# A KAM theorem for generalized Hamiltonian systems without action-angle variables

Yon Hui Jo [1], Wu Hwan Jong [2]

1: Natural Science Centre, KIM IL SUNG University, Pyongyang, Democratic People's Republic of Korea
2: Faculty of Mathematics, KIM IL SUNG University, Pyongyang, Democratic People's Republic of Korea

**Abstract**

We proved a KAM theorem on existence of invariant tori in generalized Hamiltonian systems without action-angle variables. It is a generalization of the result of de la Llave et al. [Llave, 2005] that deals with canonical Hamiltonian system.

## 1. Introduction

In this paper, we provide a KAM theorem on existence of invariant tori with a Diophantine vector for a generalized Hamiltonian system. We studied generalized Hamiltonian systems which may not be perturbations of integrable systems or may not be written in action-angle variables.

In Classical KAM theorem, persistence of invariant tori is proved for perturbations of integrable Hamiltonian systems (: nearly integrable Hamiltonian system) with action-angle variables ([Kolmogorov 1954],[Arnold 1963], [Moser 1962],[Benettin 1984],[Pöschel 2009],[Fejoz 2010],[Fejoz 2011]). Action-angle variables always exist in a neighbourhood of an invariant torus. However, a change of coordinates bringing the system to action-angle variables in general cannot be explicitly computed ([Llave, 2005]).

[Llave, 2005] obtained a KAM theorem for analytic Hamiltonian systems which are neither perturbed integrable systems nor written in action-angle variables. They obtained the theorem to solve the partial differential equation which invariant tori map satisfies by quasi-Newton method. [Haro-Llave] applied this result to the numerical computation of invariant tori. [Luque-Villanueva 2010] proved existence of low-dimensional invariant tori in Hamiltonian systems which are neither perturbed integrable systems nor written in action-angle variables. [Jong-Paek 2014] provided a KAM theorem on existence of invariant tori for differentiable Hamiltonian vector fields which are neither perturbed integrable systems nor written in action-angle variables.

[Li-Yi 2002] proved existence of invariant tori for nearly integrable generalized Hamiltonian system with action-angle variables. They applied the results to research related to the perturbation of three dimensional steady Euler fluid particle path flows. [Liu-Yihe-Huang 2005] proved existence of hyperbolic invariant tori for nearly integrable generalized Hamiltonian system with action-angle variables. [Liu-Zhu-Han 2006] proved existence of lower-dimensional invariant tori for nearly integrable generalized Hamiltonian system with action-angle variables. [Li-Yi, 2006] considered Nekhoroshev stability for nearly integrable generalized Hamiltonian system with action-angle variables using KAM theorem.

In many practical applications, one has to consider the systems that are not close to integrable system but nevertheless one has some approximately invariant tori with small enough error ― obtained, for example, by using a non-rigorous numerical method, some asymptotic expansion etc([Llave, 2005]).
On the other hand, generalized Hamiltonian systems often emerge in science and engineering, for example ecological equation of predator and prey, Lotka–Volterra model. Therefore, it will have a significance to generalize the results of [Llave, 2005] for Hamiltonian system case to generalized Hamiltonian system case.

We accordingly research a KAM theorem on existence of invariant tori for a generalized Hamiltonian system which are neither perturbed integrable systems nor written in action-angle variables. The proof for the theorem follows the proof process in [Llave, 2005] for Hamiltonian case. However there are some complexities in them because that constant structural matrix in Hamiltonian equation varies to variable structural matrix.

This paper is organized as follows: in section 2 we consider invariant manifolds given by a partial differential equations related $C^r$ semi conjugate map. The problem to find the invariant tori for generalized Hamiltonian systems is referred to find the solution of a partial differential equation (invariant tori equation) according to the result in section 2. In section 3 we consider Lagrangian property of the invariant tori and approximate Lagrangian property of the approximate invariant tori. In section 4 we construct a transformation to approximate reducible type for invariant tori equation and consider solvability of the reduced equation. In section 5 we confirm non-degenerate property of approximate solutions of invariant tori equation. Lastly in section 6 we prove existence of invariant tori of considering generalized Hamiltonian equations based on convergence of quasi-Newton method for the invariant tori equation.

## 2. Partial differential Equations given by $C^r$ semi-conjugate map and invariant manifolds

First, we define some notions and introduce some notations. Let $\mathbf{N}$, $\mathbf{Z}$, $\mathbf{R}$, $\mathbf{C}$ denote the sets of natural numbers, integers, real numbers and complex numbers, respectively. We denote set of all nonnegative integers by $\mathbf{Z}_+$. We denote Banach spaces as $\mathbf{E}$, $\mathbf{F}$, $\mathbf{G}$, $\mathbf{E}_i$, $\mathbf{F}_i$, $\mathbf{G}_i$ etc. Notation $U \overset{\circ}{\subset} X$ means that $U$ is an open set of topological space $X$.

If a subset $Y$ of topological space $X$ satisfies $Y \subset \overline{\text{int} Y}$, then $Y$ is called as standard set of $X$. Notation $Y \overset{\bullet}{\subset} X$ means that $Y$ is a standard set of topological space $X$. Let $K \subset X$ be a subset of metric space $(X, d)$. We call $diam(K) = \sup\{d(x, y) \; ; \; x, y \in K\}$ as diameter of the set $K$.

Let be $|z| = \max_{1 \le j \le m} |z_j|$ for $z \in \mathbf{R}^m$.

**Definition 2.1.** Let $\mathbf{T}^n = \mathbf{R}^n / \mathbf{Z}^n$ be an n-dimensional torus. Let $\mathbf{U}^{2n}$ be an open set of $\mathbf{R}^{2n}$ or $\mathbf{T}^n \times U$ with an open set $U$ of $\mathbf{R}^n$. Let $H : \mathbf{U}^{2n} \to \mathbf{R}$ be a real analytic function. We assume that $B(z) = (b_{ij}(z))$, $(z \in \mathbf{U}^{2n})$ is a real analytic $2n \times 2n$ matrix-valued function such that:

1) $b_{ij}(z) = -b_{ji}(z)$,

2) $\sum_{l=1}^{2n} \left( \dfrac{\partial b_{ij}}{\partial z_l}(z) b_{lk}(z) + \dfrac{\partial b_{jk}}{\partial z_l}(z) b_{li}(z) + \dfrac{\partial b_{ki}}{\partial z_l}(z) b_{lj}(z) \right) = 0$.

Ordinary differential equation

$$\frac{dz}{dt} = B(z) \nabla H(z) \tag{2.1}$$

is called as generalized Hamiltonian system with Hamiltonian $H$ and structure matrix $B$ ([Zhao 1995]). For section 2 it is good enough for $H$ and $b_{ij}$ to be $C^\infty$ function, but for 3-6 sections real analyticity is necessary. We denote non-prolongable solution satisfying initial condition $z(0) = \zeta \in \mathbf{U}^{2n}$ of differential equation (2.1) by $\varphi_t(\zeta)$, $(t \in I(\zeta))$. Here $I(\zeta)$ is the maximum existing interval of solution for the initial value problem. $I(\zeta)$ is an open interval involving 0.

Let $<\xi, \eta>$ denote the Euclidean standard scalar product on $\mathbf{R}^{2n}$. Let $\omega = \omega_z(\xi, \eta)$ $(z \in \mathbf{U}^{2n}, (\xi, \eta) \in \mathbf{R}^n \times \mathbf{R}^n)$ be a differential 2-form on $\mathbf{U}^{2n}$. If there exists a $2n \times 2n$ matrix-valued function $A(z)$, $(z \in \mathbf{U}^{2n})$ such that
$$\omega_z(\xi, \eta) = <\xi, A(z)\eta>, \ (z \in \mathbf{U}^{2n}, \ \xi, \eta \in \mathbf{R}^{2n}),$$
then $A(z)$ is called as representation matrix of differential 2-form $\omega$. We assume that matrix $B(z)$ in (2.1) is invertible for any $z \in \mathbf{U}^{2n}$. And we assume that for differential 2-form $\Omega(\xi, \eta) = <\xi, B(z)^{-1}\eta>$, $(\xi, \eta \in \mathbf{R}^{2n})$ on $\mathbf{U}^{2n}$ there exist a differential 1-form $\alpha = c(z) \cdot dz = \sum_{j=1}^{2n} c_j(z) dz_j$ such that $\Omega = d\alpha$.

**Lemma 2.1.** $\Omega$ is a symplectic form on $\mathbf{U}^{2n}$ and $(\varphi_t)^*\Omega = \Omega$ holds for any $t$.

**Proof.** Since $B(z)$ is non-degenerate skew-symmetric, $J(z) = B(z)^{-1}$ is also non-degenerate skew-symmetric. Therefore $\Omega$ is differential 2-form on $\mathbf{U}^{2n}$. Since the assumption $\Omega = d\alpha$, $\Omega$ is closed differential form so $\Omega$ is symplectic form on $\mathbf{U}^{2n}$. Since $\varphi_t$ is Poisson map ([Hairer, 2006]),
$$D\varphi_t(z)B(z)D\varphi_t(z)^T = B(\varphi_t(z))$$
holds. By taking the inverse matrix of both sides in the last relation, we have
$$[D\varphi_t(z)^T]^{-1} B(z)^{-1} D\varphi_t(z)^{-1} = B(\varphi_t(z))^{-1}.$$
And multiplying to this relation from the left by $D\varphi_t(z)^T$ and from right by $D\varphi_t(z)$, we have
$$B(z)^{-1} = D\varphi_t(z)^T B(\varphi_t(z))^{-1} D\varphi_t(z).$$
Therefore,
$$D\varphi_t(z)^T J(\varphi_t(z)) D\varphi_t(z) = J(z).$$
This means $(\varphi_t)^*\Omega = \Omega$. □

**Definition 2.2.** Given $k = (k_1, \cdots, k_n) \in \mathbf{Z}^n$, let $|k|$ denote the absolute sum norm $|k_1| + \cdots + |k_n|$ of $k$. Let $\omega = (\omega_1, \cdots, \omega_n) \in \mathbf{R}^n$. If there exist $\gamma > 0$ and $\sigma > 0$ such that
$$|k \cdot \omega| \geq \gamma |k|^{-\sigma}, \ (k \in \mathbf{Z}^n \setminus \{0\}), \tag{2.2}$$
then we call that $\omega \in \mathbf{R}^n$ is $(\gamma, \sigma)$-Diophantine. (2.2) is called as Diophantine condition. We denote the set of all $(\gamma, \sigma)$-Diophantine vectors $\omega \in \mathbf{R}^n$ by $D_n(\gamma, \sigma)$.

**Definition 2.3.** Let us $\rho > 0$ is a given positive number. Let $U_\rho = \{\theta \in \mathbf{C}^n ; \ |\mathrm{Im}\theta| < \rho\}$ and $m \in \mathbf{N}$, $V \overset{\circ}{\subset} \mathbf{R}^m$. For $U \overset{\bullet}{\subset} \mathbf{C}^n$, we let

$$\mathcal{P}(U, V) = \{K : U \to V \ ; \ K \text{ is continuous on } U, \text{ one-periodic in all its variables on } U \text{ and}$$
$$\text{real analytic on the } \operatorname{int} U \}$$

Then $\mathcal{P}(\overline{U}_\rho, \mathbf{R}^m)$ becomes the Banach space with norm $\|K\|_\rho = \sup_{\theta \in U_\rho} |K(\theta)|$ (See [Dieudonné 1960]). We denote $\mathcal{P}(\overline{U}_\rho, \mathbf{U}^{2n})$ by $\mathcal{P}_\rho$ simply. Let $\mathbf{T}_\rho^n = \{\theta \in \mathbf{C}^n / \mathbf{Z}^n \ ; \ |\operatorname{Im}\theta| \leq \rho\}$. Then the set $\mathcal{P}(\overline{U}_\rho, V)$ is identified to the set

$$\mathcal{A}(\mathbf{T}_\rho^n, V) = \{K : \mathbf{T}_\rho^n \to V \ ; \ K \text{ is continuous on } \mathbf{T}_\rho^n \text{ and real analysis on } Int\mathbf{T}_\rho^n\}.$$

Let $K_0 \in \mathcal{P}(\overline{U}_\rho, \mathbf{R}^{2n})$ be given. We fix $r > 0$ satisfying $diam(K_0(\overline{U}_\rho)) < r$ and let

$$\mathcal{B}_r = \{z \in \mathbf{R}^{2n} \ ; \ \inf_{|\operatorname{Im}\theta| \leq \rho} |z - K_0(\theta)| < r\}.$$

Then $\mathcal{B}_r$ is a neighborhood of set $K_0(\overline{U}_\rho)$ on $\mathbf{R}^{2n}$.

In this section, we clarify that two dynamical systems are $C^r$ semi-conjugate with some $C^r$ semi-conjugacy if and only if the $C^r$ semi-conjugacy satisfies a partial differential equation. In this case, the image of the $C^r$ semi-conjugacy becomes the invariant manifold of the second dynamical system. This result refers the existence of invariant manifold for dynamical system to the existence of solution for the partial differential equation.

**Definition 2.4.** Let $r \in \mathbf{N}$ or $r = \infty$. Suppose that $M$ is $C^{r+1}$ manifold, $X : M \to TM$ is $C^r$ vector field on $M$ and $x \in M$. We denote the non-prolongable solution of the initial problem

$$\dot{y} = X(y) \tag{2.3}$$
$$y(0) = x \tag{2.3}_0$$

by $\varphi_x : I_X(x) \to M$, where $I_X(x)$ is an open interval involving 0.

Let $\mathcal{D}_X = \{(t, x)\} \in \mathbf{R} \times M \ ; \ t \in I_X(x)\}$. Then $\mathcal{D}_X$ becomes an open set of $\mathbf{R} \times M$ involving $\{0\} \times M$. The map

$$\varphi : \mathcal{D}_X \to M \ ; \ \varphi(t, x) = \varphi_x(t), \quad ((t, x) \in \mathcal{D}_X) \tag{2.4}$$

is called as $C^r$ integral of vector field $X$ or $C^r$ integral of differential equation (2.3) ([Abraham-Marsden 1978]). $C^r$ integral $\varphi : \mathcal{D}_X \to M$ becomes $C^r$ map. By definition of $C^r$ integral and property of solution of autonomous system (2.3), for any $x \in M$,

$$\varphi(0, x) = x \tag{2.5}$$
$$\varphi(t, \varphi(s, x)) = \varphi(t + s, x) \quad (\forall s \in I(x), \forall t \in I(\varphi(s, x))) \tag{2.6}$$

hold. If $I_X(x) = \mathbf{R}$ holds for any $x \in M$, then the $C^r$ vector field $X$ is called as complete ([Abraham-Marsden 1978]).

If $X$ is complete, $\mathcal{D}_X = \mathbf{R} \times M$ and $C^r$ integral $\varphi : \mathbf{R} \times M \to M$ of $X$ satisfies
$$\varphi(0, x) = x$$
$$\varphi(t, \varphi(s, x)) = \varphi(t+s, x) \quad (\forall s \in \mathbf{R}, \forall t \in \mathbf{R}).$$

Then $C^r$ map $\varphi : \mathbf{R} \times M \to M$ is called $C^r$ flow of $X$. In this time, we define $\varphi_t : M \to M$ by $\varphi_t(x) = \varphi(t, x)$ for any $t \in \mathbf{R}$. Then $\varphi_t : M \to M$ becomes $C^r$ diffeomorphism and satisfies
$$\varphi_0(x) = x$$
$$(\varphi_t \circ \varphi_s)(x) = \varphi_{t+s}(x), \quad (\forall s, t \in \mathbf{R}).\ \square$$

**Definition 2.5.** Let $S \subset M$. If
1) $I(x) = \mathbf{R}$ for any $x \in S$,
2) $\varphi(t, x) \in S, \ (t \in \mathbf{R})$,

then $S$ is called as invariant set with differential equation (2.3) or with $C^r$ integral $\varphi : \mathcal{D}_X \to M$ of vector field $X$. If an invariant set $S$ with differential equation (2.3) has the structure of $C^r$ manifold, then $S$ is called as the invariant manifold with differential equation (2.3).

If $C^r$ invariant manifold $S$ with differential equation (2.3) is homeomorphism to torus $\mathbf{T}^n = \mathbf{R}^n / \mathbf{Z}^n$, then $S$ is called as invariant torus with differential equation (2.3).

**Definition 2.6.** Let $M$, $N$ be $C^{r+1}$ manifolds, $X : M \to TM$, $Y : N \to TN$ be $C^r$ vector fields, $\varphi : \mathcal{D}_X \to M$ be $C^r$ integral of $X$, $\psi : \mathcal{D}_Y \to N$ be $C^r$ integral of $Y$ and $u : M \to N$ be $C^r$ map. If relation

$$u(\varphi(t, x)) = \psi(t, u(x)), \quad (\forall x \in M, \ \forall t \in I_X(x) \cap I_Y(u(x))) \tag{2.7}$$

holds, we call that $C^r$ integrals $\varphi : \mathcal{D}_X \to M$ and $\psi : \mathcal{D}_Y \to N$ are $C^r$ semi-conjugate with map $u$. In this time, $u$ is called as $C^r$ semi-conjugacy between $\varphi$ and $\psi$.

**Definition 2.7.** Let $M$, $N$ be $C^{r+1}$ manifolds and $X : M \to TM$, $Y : N \to TN$ be complete $C^r$ vector fields. Let $\varphi : \mathbf{R} \times M \to M$, $\psi : \mathbf{R} \times N \to N$ be $C^r$ flows of $X, Y$ respectively. Let $u : M \to N$ be homeomorphism. If the relation

$$(u \circ \varphi_t)(x) = (\psi_t \circ u)(x), \quad (\forall t \in \mathbf{R}, \ \forall x \in M) \tag{2.8}$$

holds, we call that $C^r$ flows $\varphi : \mathbf{R} \times M \to M$ and $\psi : \mathbf{R} \times N \to N$ are topological conjugate with $u$. In this time, $u$ is called as topological conjugacy between $\varphi$ and $\psi$.

**Theorem 2.1.** Let $M$, $N$ be $C^{r+1}$ manifolds. Let $X : M \to TM$, $Y : N \to TN$ be $C^r$ vector fields. $C^r$ integral $\varphi : \mathcal{D}_X \to M$ of $X$ and $C^r$ integral $\psi : \mathcal{D}_Y \to N$ of $Y$ are $C^r$ semi-conjugate with $C^r$ conjugacy $u : M \to N$ if and only if $u$ satisfies partial differential equation

$$Tu(X(x)) = Y(u(x)), \ (x \in M). \tag{2.9}$$

**Proof.** Necessity: Suppose that (2.7) holds. Fix $x \in M$ arbitrarily. From (2.7) for all $t \in I_X(x) \cap I_Y(u(x))$,

$$u(\varphi_x(t)) = \psi_{u(x)}(t) \tag{2.10}$$

holds. Since differentiating both sides of (2.10) with respect to $t$, we have

$$Tu(\dot\varphi_x(t)) = \dot\psi_{u(x)}(t),$$

hence

$$Tu(X(\varphi_x(t))) = Y(\psi_{u(x)}(t)).$$

By taking $t = 0$ in the last expression we have (2.9).

Sufficiency: Let's suppose that (2.9) holds. Fix $x \in M$ arbitrarily. Let's consider the map

$$\zeta : I_X(x) \cap I_Y(u(x)) \to N \ ; \ \zeta(t) = u(\varphi_x(t)). \tag{2.11}$$

Since

$$\frac{d}{dt}\zeta(t) = \frac{d}{dt}u(\varphi_x(t)) = Tu(\dot\varphi_x(t)) = Tu(X(\varphi_x(t))) = Y(u(\varphi_x(t))) = Y(\zeta(t)),$$

$\zeta : I_X(x) \cap I_Y(u(x)) \to N$ is a solution of the differential equation

$$\dot z = Y(z) \tag{2.12}$$

satisfying initial condition $z(0) = u(x)$. On the other hand, from the definition of $\psi$,

$$\psi_{u(x)}(t), \ t \in I_X(x) \cap I_Y(u(x))$$

is also the solution of differential equation (2.12) satisfying initial condition $z(0) = u(x)$.

From uniqueness of the solution of initial problem of ordinary differential equation,

$$u(\varphi(t, \ x)) = u(\varphi_x(t)) = \zeta(t) = \psi_{u(x)}(t) = \psi(t, \ u(x))$$

holds for any $t \in I_X(x) \cap I_Y(u(x))$. □

**Corollary 2.1.** Let's denote non-prolongable solution of Hamiltonian system (2.1) satisfying initial condition $z(0) = \zeta \in \mathbf{U}^{2n}$ by $\varphi_t(\zeta), \ (t \in I(\zeta))$. For $C^\infty$ map $K : \mathbf{T}^n \to \mathbf{U}^{2n}$ we define operator $\partial_\omega$ by $\partial_\omega K(\theta) = DK(\theta)\omega, \ (\theta \in \mathbf{T}^n)$. Then $C^\infty$ map $K : \mathbf{T}^n \to \mathbf{U}^{2n}$ satisfies

$$\varphi_t(K(\theta)) = K(\theta + \omega t), \ (\forall \theta \in \mathbf{R}^n, \ \forall t \in I(K(\theta)))$$

if and only if $K$ satisfies the first order partial differential equation system

$$\partial_\omega K(\theta) = B(K(\theta))\nabla H(K(\theta)), \ (\theta \in \mathbf{T}^n).$$

**Remark 2.1.** In the case of $M = \mathbf{R}^m, \ N = \mathbf{R}^n$,

$$X(x) = \begin{pmatrix} f_1(x_1, \cdots, x_m) \\ \vdots \\ f_m(x_1, \cdots, x_m) \end{pmatrix} \text{ and } Y(y) = \begin{pmatrix} g_1(y_1, \cdots, y_n) \\ \vdots \\ g_n(y_1, \cdots, x_n) \end{pmatrix},$$

(2.7) is given as follows:

$$\sum_{j=1}^m \frac{\partial u_i}{\partial x_j}(x_1, \cdots, x_m) f_j(x_1, \cdots, x_m) = g_i(u_1(x_1, \cdots, x_m), \cdots, u_n(x_1, \cdots, x_m)), \ (i = 1, \cdots, n).$$

(2.13)

**Theorem 2.2.** Let $\mathbf{E}, \mathbf{F}$ be Banach spaces, $M$ is $C^{r+1}\mathbf{E}$ manifold, $N$ is $C^{r+1}\mathbf{F}$ manifold and $X : M \to TM$, $Y : N \to TN$ be $C^r$ vector fields. We suppose that $C^r$ injection $u : M \to N$ satisfies partial differential equation (2.9). Then

    1) $S = u(M)$ has a structure of $C^r \mathbf{E}$ manifold.

    2) $T_{u(x)}S = T_x u(T_x M)$, $(x \in M)$.

    3) If $X$ complete, then $S = u(M)$ is the invariant manifold for differential equation (2.12).

    4) If $S = u(M)$ is closed set of $N$, then $S = u(M)$ is the invariant manifold for differential equation (2.12).

    5) If $M$ is compact, then vector field $X : M \to TM$ and $Y : S \to TS$ are complete, $S = u(M)$ is the invariant manifold for differential equation (2.12) and $u : M \to S$ is homeomorphism. Moreover $C^r$ flow $\varphi : \mathbf{R} \times M \to M$ of $X : M \to TM$ and $C^r$ flow $\psi : \mathbf{R} \times S \to S$ of $Y : S \to TS$ are topological conjugate with $u$.

**Proof.** 1): Suppose $u : M \to N$ is $C^r$ injection. Now we introduce the structure of $C^r$ manifold for set $S = u(M)$. We denote the inverse map of bijection $u : M \to u(M)$ by $u^{-1} : S = u(M) \to M$. Let $\mathcal{A} = \{(U_i, \alpha_i), i \in J\}$ is the $C^r$ atlas of $M$. We put $V_i = u(U_i)$ for any $i \in J$ and define map $\beta_i : V_i \to \mathbf{E}$ as $\beta_i(y) = \alpha_i(u^{-1}(y)), (y \in V_i)$. Then we prove $\mathcal{B} = \{(V_i, \beta_i), i \in J\}$ is $C^r$ atlas of set $S = u(M)$. First of all, we have

$$\bigcup_{i \in J} V_i = \bigcup_{i \in J} u(U_i) = u(\bigcup_{i \in J} U_i) = u(M) = S$$

obviously. For any $i, j \in J$

$$\beta_i(V_i \cap V_j) = (\alpha_i \circ u^{-1})(u(U_i) \cap u(U_j)) = \alpha_i(U_i \cap U_j)$$

is open set of Banach space $\mathbf{E}$. Now we prove that coordinate transformations are $C^r$. Since

$$\beta_j \beta_i^{-1} = \alpha_j u^{-1}(\alpha_i u^{-1})^{-1} = \alpha_j u^{-1} u \alpha_i^{-1} = \alpha_j \alpha_i^{-1} : \alpha_i(U_i \cap U_j) \to \alpha_j(U_i \cap U_j),$$

$\beta_j \beta_i^{-1} : \beta_i(V_i \cap V_j) \to \beta_j(V_i \cap V_j)$ is $C^r$. Therefore $\mathcal{B} = \{(V_i, \beta_i), i \in J\}$ is $C^r$ atlas of set $S = u(M)$. Hence set $S = u(M)$ has the structure of $C^r \mathbf{E}$ manifold.

    2): Next we prove that

$$T_y S = T_x u(T_x M) \tag{2.14}$$

for $x \in M$, $y = u(x) \in S = u(M)$. First we prove $T_y S \subset T_x u(T_x M)$. Let $w \in T_y S$. Then there exists open interval $I \subset \mathbf{R} : 0 \in I$ and $C^r$ map $c : I \to S = u(M)$ with $c(0) = y$ such that $w = [c]_y$. Here $[c]_y$ is the equivalent class of $C^1$ curves on $M$ which are tangent to $c$ at $y$ (See [Abraham-Marsden 1978]). Now let $x \in M$ is an element satisfies $y = u(x)$ and define $C^1$ map $\tilde{c} : I \to M$ as $\tilde{c} = u^{-1} \circ c$. Then we have $\tilde{c}(0) = u^{-1}(c(0)) = u^{-1}(u(x)) = x$. Here we put $v = [\tilde{c}]_x (\in T_x M)$. In this time we have

$$T_x u(v) = T_x u([\tilde{c}]_x) = [u \circ \tilde{c}]_{u(x)} = [c]_y = w.$$

Therefore $w \in T_x u(T_x M)$ holds. Hence we have proved $T_y S \subset T_x u(T_x M)$.

Next we prove $T_x u(T_x M) \subset T_y S$. Let take $v \in T_x M$ arbitrarily. Then there exists open interval $I \subset \mathbf{R} : 0 \in I$ and $C^1$ map $\tilde{c} : I \to M$ with $\tilde{c}(0) = x$ such that $v = [\tilde{c}]_x$. Let's consider $C^1$ map $c = u \circ \tilde{c} : I \to S = u(M)$. Then $c(0) = u \circ \tilde{c}(0) = u(x) = y$. We put $w = [c]_y \in T_y S$. Then since

$$T_x u(v) = T_x u([\tilde{c}]_x) = [u \circ \tilde{c}]_{u(x)} = [c]_y = w,$$

we have $T_x u(T_x M) \subset T_y S$. Therefore (2.14) is proved.

3): Let's take $y \in S$ arbitrarily. Then there is a unique $x \in M$ such that $y = u(x)$. Since $X$ is complete, we have

$$I_X(x) \cap I_Y(u(x)) = \mathbf{R} \cap I_Y(u(x)) = I_Y(u(x)) = I_Y(y).$$

Then since

$$\psi(t, u(x)) = u(\varphi(t, x)) \in u(M) = S$$

for any $t \in I_X(x) \cap I_Y(u(x)) = I_Y(y)$ from theorem 2.1, $S$ is the invariant manifold of (2.12).

4): Let's suppose that (2.9) holds for $C^r$ injection $u : M \to N$. $S = u(M)$ is not submanifold of $N$ generally but $S$ has the structure of $C^r$ manifold itself as we see in 1) and the tangent space of $S$ at $y \in S$ is given by (2.14). Let's suppose that $y = u(x)$ for some $x \in M$. Since we have

$$Y(y) = Y(u(x)) = Tu(X(x)) = T_x u(X(x)) \in T_x u(T_x M) = T_y S$$

from (2.9), $\tilde{Y} : S \to TS$ ; $\tilde{Y}(y) = Y(y)$ $(y \in S)$ defines a $C^r$ vector field on $C^r$ manifold $S$.

We denote the $C^r$ integral of differential equation

$$\dot{z} = \tilde{Y}(z), \ (z \in S) \tag{2.15}$$

with vector field $\tilde{Y} : S \to TS$ by $\tilde{\psi} : \mathcal{D}_{\tilde{Y}} \to S$. Let's take $y \in S$ and $t \in I_{\tilde{Y}}(y)$ arbitrarily. Then since

$$\frac{d}{dt}\tilde{\psi}_y(t) = \tilde{Y}(\tilde{\psi}_y(t)) = Y(\tilde{\psi}_y(t)),$$

$\tilde{\psi}_y(t), t \in I_{\tilde{Y}}(y)$ is the solution of (2.12) satisfies $z(0) = y$. Therefore we have $I_{\tilde{Y}}(y) \subset I_Y(y)$ and $\tilde{\psi}_y(t) = \psi_y(t), (\forall t \in I_{\tilde{Y}}(y))$.

Now we prove $I_{\tilde{Y}}(y) = I_Y(y)$. If we deny this fact, then there exists $s \in I_Y(y)$ such that $t < s, \ (\forall t \in I_{\tilde{Y}}(y))$ or $s' \in I_Y(y)$ such that $s' < t, \ (\forall t \in I_{\tilde{Y}}(y))$. We consider the first case. The argument for the second case is as same as the first case. Let $s_0 = \inf\{s \in I_Y(y) \ ; \ t < s, \ (\forall t \in I_{\tilde{Y}}(y))\}$. Then $s_0 = \sup I_{\tilde{Y}}(y)$. There exists limit $z_0 = \lim_{t \to s_0} \tilde{\psi}_y(t) = \lim_{t \to s_0} \psi_y(t) \in N$. Then since $S$ is closed set of $N$ by the assumption, $z_0 \in S$ holds. There exists solution $\xi(t)$, $t \in (-\delta, \delta)$ of (2.12) with $z(0) = z_0$. If we define

$$\eta(t) = \begin{cases} \tilde{\psi}_y(t), & (t \in I_{\tilde{Y}}(y)) \\ \xi(t - s_0), & t \in [s_0, s_0 + \delta) \end{cases},$$

then $\eta(t)$, $t \in I_{\tilde{Y}}(y) \cup [s_0, s_0 + \delta)$ is the solution of (2.12) with $z(0) = y$. This contradicts that $\psi_y(t)$, $t \in I_{\tilde{Y}}(y)$ is non-prolongable solution of (2.12) with $z(0) = y$. Therefore we have $I_{\tilde{Y}}(y) = I_Y(y)$. Sine $\psi_y(t) = \tilde{\psi}_y(t) \in S$ for any $t \in I_Y(y) = I_{\tilde{Y}}(y)$, $S$ is invariant manifold of (2.12).

5) Suppose that $M$ is compact. Then $S = u(M)$ is also compact and vector field $X : M \to TM$ and $Y : S \to TS$ are complete. From 3) $S = u(M)$ is the invariant manifold of (2.12). Since bijective continuous map from compact topological space to Hausdorff topological space is homeomorphism, $u : M \to S$ is homeomorphism. Then from theorem 2.1

$$u(\theta + \omega t) = (u \circ \varphi_t)(\theta) = u(\varphi(t, \theta)) = \psi(t, u(\theta)) = \psi_t(u(\theta)) = (\psi_t \circ u)(\theta)$$

holds for any $x \in M$, $t \in \mathbf{R}$. □

Now consider the generalized Hamiltonian system

$$(2.1): \frac{dz}{dt} = B(z)\nabla H(z).$$

For $\omega \in \mathbf{R}^n$ let us define $\partial_\omega u(\theta) = Du(\theta)\omega$.

**Corollary 2.2.** Suppose that $\omega \in \mathbf{R}^n$ is a Diophantine vector and $C^\infty$ injection $u : \mathbf{T}^n \to \mathbf{U}^{2n}$ satisfies the partial differential equation

$$\partial_\omega u(\theta) = B(u(\theta))\nabla H(u(\theta)). \tag{2.16}$$

Then

1) $\mathcal{T} = u(\mathbf{T}^n)$ is invariant torus of (2.1).

2) $u : \mathbf{T}^n \to \mathcal{T}$ is homeomorphism.

3) The flow of (2.1) restricted on $\mathcal{T}$ is defined on $\mathbf{R} \times \mathcal{T}$ and it is topologically conjugate for homeomorphism $u : \mathbf{T}^n \to \mathcal{T}$ with The Kronecker's flow $\varphi(t, \theta) = \theta + \omega t$, $(t \in \mathbf{R}, \theta \in \mathbf{T}^n)$.

Therefore the restricted flow of (2.1) on $\mathcal{T}$ is quasiperiodic.

**Proof.** Let define $C^\infty$ vector field $X$ on $\mathbf{T}^n$ by

$$X : \mathbf{T}^n \to \mathbf{R}^n \ ; \ X(\theta) = \omega. \tag{2.17}$$

Then $C^\infty$ flow of $X$ coincides with the Kronecker's flow. We define that $C^\infty$ vector field $Y$ on $\mathbf{U}^{2n}$ by

$$Y : \mathbf{U}^{2n} \to \mathbf{R}^{2n} \ ; \ Y(z) = B(z)\nabla H(z). \tag{2.18}$$

We denote $C^\infty$ integral of $Y$ by $\psi(t, z)$, $(z \in \mathbf{U}^{2n}, t \in I_Y(z))$. Then (2.9) is written by (2.16).

Therefore from 2.2 $\mathcal{T} = u(\mathbf{T}^n)$ is invariant manifold of (2.1) and $u : \mathbf{T}^n \to \mathcal{T}$ is homeomorphism. Hence $\mathcal{T}$ is invariant torus of (2.1). From theorem 2.2 $\psi$ restricted on $\mathcal{T}$ is defined on $\mathbf{R} \times \mathcal{T}$ and $\psi \mid \mathcal{T}$ is topologically conjugate to the Kronecker's flow $\varphi(t, \theta) = \theta + \omega t$, $(t \in \mathbf{R}, \theta \in \mathbf{T}^n)$ with $u : \mathbf{T}^n \to \mathcal{T}$. Thus $u(\theta + \omega t) = \psi_t(u(\theta))$, $(t \in \mathbf{R}, \theta \in \mathbf{T}^n)$. □

From Corollary 2.1 the problem to obtain invariant torus for generalized Hamiltonian system (2.1) refers to the problem to obtain solution of nonlinear partial differential equation (2.16).

**Definition 2.8.** Based on Corollary 2.1 we call (2.16) as invariant torus equation. Let's consider map $K \in \mathcal{P}(\overline{U}_\rho, \mathbf{U}^{2n})$. Then we call function $e(\theta) = \partial_\omega K(\theta) - B(K(\theta)) \nabla H(K(\theta)), \theta \in \overline{U}_\rho$ as error of map $K$ for equation (2.16), real number $\varepsilon = \|e\|_\rho = \sup\{\|e(\theta)\| ; \theta \in \overline{U}_\rho\}$ as error-size of map $K$ for equation (2.16). And we call $K$ as approximate solution of equation (2.16) with error $\varepsilon$. We call approximate solution of (2.16) with zero error as solution of (2.16). If $K \in \mathcal{P}(\overline{U}_\rho, \mathbf{U}^{2n})$ is solution of (2.16), then

$$\partial_\omega K(\theta) = B(K(\theta)) \nabla H(K(\theta))$$

for any $\theta \in \overline{U}_\rho$.

## 3. Lagrangian invariant torus

In this section we show that the image of the solution of the first order nonlinear partial differential equation (2.16), driven by the generalized Hamiltonian dynamical system (2.1), becomes Lagrangian invariant torus with symplectic structure with (2.1) and the image of an approximate solution becomes the approximate Lagrangian invariant torus. This fact plays an essential role to get the existence of the solution applying Newton method to the first order nonlinear partial differential equation (2.16).

Let $m, n$ be natural numbers. Let us denote $\mathbf{R}$ or $\mathbf{C}$ by $\mathbf{K}$. Let us introduce maximum norm for an element $x \in \mathbf{K}^n$ of finite dimensional Euclidean Space $\mathbf{K}^n$ as $|x| = \max_{1 \le i \le n} |x_i|$. Taking $m \times n$ matrix $A$ as the element of $\mathbf{K}^{mn}$, the maximum norm is given as

$$|A| = \max\{|a_{ij}| ; 1 \le i \le m, i \le j \le n\}.$$

On the other hand, we define the *absolute sum norm* of $x \in \mathbf{K}^n$ as $|x|_1 = |x_1| + \cdots + |x_n|$. Then

$$\frac{1}{n}|x|_1 \le |x| \le |x|_1, \quad (\forall x \in \mathbf{K}^n)$$

holds. We define the *operator norm* of $m \times n$ matrix $A = (a_{ij}) \in \mathbf{K}^{mn}$ as

$$\|A\| = \sup_{x \in \mathbf{K}^n \setminus \{0\}} \frac{|Ax|}{|x|} = \sup_{|x|=1} |Ax|.$$

**Proposition 3.1.** Let us $A$ is $m \times n$ matrix, $B$ is $n \times l$ matrix, and $x \in \mathbf{K}^n$. Then the followings hold:

1) $|AB| \le n |A| \cdot |B|$
2) $|Ax| \le n |A| \cdot |x|$

3) $\|A\| \leq n|A|$

4) $|Ax| \leq \|A\| \cdot |x|$

5) $\|AB\| \leq \|A\| \cdot \|B\|$

6) $|A| \leq \|A\|$

7) $\|A^T\| \leq m|A|$

8) $\|A^T\| \leq m\|A\|$

**Proof.** 1): $|\sum_{k=1}^{n} a_{ik} b_{kj}| \leq \sum_{k=1}^{n} |a_{ik}||b_{kj}| \leq \sum_{k=1}^{n} |A||B| = n|A||B|$.

2): This is a special case of 1).

3): From 1)
$$\|A\| = \sup_{|x|=1} |Ax| \leq \sup_{|x|=1} (n|A|\|x|) = n|A|.$$

4), 5): From operator norm of matrix, it is directly obtained.

6): Let
$$A = \begin{pmatrix} a_{11} & a_{12} & \cdots & a_{1n} \\ a_{21} & a_{22} & \cdots & a_{2n} \\ \vdots & \vdots & & \vdots \\ a_{m1} & a_{m2} & \cdots & a_{mn} \end{pmatrix} = \begin{pmatrix} a^1 & a^2 & \cdots & a^n \end{pmatrix}.$$

If we apply 4) for $e_i = (0, \cdots, 1, \cdots, 0)^T$, $(i = 1, \cdots, n)$, then
$|a^i| = |Ae_i| \leq \|A\| \cdot |e_i| = \|A\|$. From $|a^i| \leq \|A\|$ we have $|A| = \max_{1 \leq i \leq n} |a^i| \leq \|A\|$. Moreover
$|A| = \max_{1 \leq i \leq n} |a^i| \leq \|A\|$.

7): From 3) and the definition of maximum norm, $\|A^T\| \leq m|A^T| = m|A|$.

8): From 7) and 6), $\|A^T\| \leq m|A| \leq m\|A\|$. □

**Proposition 3.2.** Let us $m, n, l \in \mathbf{N}$ and $A_{ij}$, $(1 \leq i \leq m, 1 \leq j \leq n)$ are $l \times l$ matrices. We define $ml \times nl$ matrix $A$ as
$$A = \begin{pmatrix} A_{11} & \cdots & A_{1n} \\ \vdots & & \vdots \\ A_{m1} & \cdots & A_{mn} \end{pmatrix}.$$

Then $\|A\| \leq \max_{1 \leq i \leq m} \sum_{j=1}^{n} \|A_{ij}\|$.

**Proof.** Let $x = \begin{pmatrix} x_1 \\ \vdots \\ x_n \end{pmatrix} \in \mathbf{R}^{nl}$, $x_1, \cdots, x_n \in \mathbf{R}^l$. Then

$$\|A\| = \sup_{\|x\|=1} \|Ax\| \leq \sup_{\|x\|=1} \max_{1 \leq i \leq m} \left| \sum_{j=1}^{n} A_{ij} x_j \right| \leq \max_{1 \leq i \leq m} \sum_{j=1}^{n} \sup_{\|x\|=1} |A_{ij} x_j| = \max_{1 \leq i \leq m} \sum_{j=1}^{n} \|A_{ij}\|. \square$$

**Proposition 3.3.** (Cauchy's estimate) Let $z_0 \in \mathbf{C}$ and complex function $f(z)$ is analytic in $|z - z_0| < R$. We suppose that for $a_n \in \mathbf{C}$, $(n = 0, 1, 2, \cdots)$ $f(z) = \sum_{n=0}^{\infty} a_n(z-z_0)^n$. We put $M(r) = \max\{|f(\zeta)|\,;\ |\zeta - z_0| = r\}$ for any $0 < r < R$. Then we have

$$|a_n| \leq \frac{M(r)}{r^n},\ (m = 0, 1, 2, \cdots). \tag{5.2}$$

**Proof.** See [Ahlfors 1966].

**Corollary 3.1.** Let $\rho > 0$ and $U_\rho = \{\theta \in \mathbf{C}\,;\ |\mathrm{Im}\,\theta| < \rho\}$. Let $\|f\|_\rho = \sup_{\theta \in U_\rho} |f(\theta)|$ for bounded analytic function $f: U_\rho \to \mathbf{C}$. Then for any $0 < \delta < \rho$,

$$\|f'\|_{\rho-\delta} \leq \delta^{-1} \|f\|_\rho$$

holds.

**Proof.** Fix $\theta_0 \in U_{\rho-\delta}$ arbitrarily. Taking $\theta$ arbitrarily from closed ball $B(\theta_0, \delta)$ with radius $\delta$, since $|\mathrm{Im}\,\theta| \leq |\mathrm{Im}(\theta - \theta_0)| + |\mathrm{Im}\,\theta_0| < \rho - \delta + \delta = \rho$, $B(\theta_0, \delta) \subset U_\rho$ holds. From Cauchy's estimate,

$$|f'(\theta)| \leq \delta^{-1} \max\{|f(\theta)|,\ |\theta - \theta_0| = \delta\} \leq \delta^{-1} \sup\{|f(\theta)|,\ \theta \in U_\rho\} = \delta^{-1} \|f\|_\rho,$$

moreover $\|f'\|_\rho \leq \delta^{-1} \|f\|_\rho$. □

**Corollary 3.2.** Let $n$, $m$ be natural numbers and $\rho > 0$, we let

$$U_\rho = \{\theta \in \mathbf{C}^n\,;\ |\mathrm{Im}\,\theta| = \max_{1 \leq j \leq n} |\mathrm{Im}\,\theta_j| < \rho\}.$$

For bounded analytic function $f: U_\rho \to \mathbf{C}^m$, we let

$$\|f\|_\rho = \sup_{\theta \in U_\rho} |f(\theta)| = \sup_{\theta \in U_\rho} \max_{1 \leq j \leq m} |f_j(\theta)|.$$

For $\theta \in U_\rho$, we let

$$Df(\theta) = \begin{pmatrix} \frac{\partial f_1}{\partial \theta_1}(\theta) & \cdots & \frac{\partial f_1}{\partial \theta_n}(\theta) \\ \vdots & & \vdots \\ \frac{\partial f_m}{\partial \theta_1}(\theta) & \cdots & \frac{\partial f_m}{\partial \theta_n}(\theta) \end{pmatrix}.$$

Let $\|Df(\theta)\|$ be a operator norm of matrix $Df(\theta)$ and $\|Df\|_\rho = \sup_{\theta \in U_\rho} \|Df(\theta)\|$.

Then for arbitrary $0 < \delta < \rho$

$$\|Df\|_{\rho-\delta} \leq n\delta^{-1} \|f\|_\rho$$

holds.

**Proof.** For any $\theta_0 \in U_{\rho-\delta}$, $B(\theta_0, \delta) \subset U_\rho$ holds. Fix $1 \leq i \leq m$, $1 \leq j \leq n$ arbitrarily. From Cauchy's estimate

$$\left|\frac{\partial f_i}{\partial \theta_j}(\theta)\right| \leq \delta^{-1} \max\{|f_j(\theta)|,\ |\theta_j - \theta_{0j}| = \delta\} \leq \delta^{-1} \|f\|_\rho$$

and from 3) of proposition 3.2 and Corollary 3.1, we have

$$\| Df(\theta) \| \leq n | Df(\theta) | = n \max \left\{ \left| \frac{\partial f_i}{\partial \theta_j}(\theta) \right| ; \ 1 \leq i \leq n, \ 1 \leq j \leq m \right\} \leq n\delta^{-1} \| f \|_\rho.$$

Hence we have $\| Df \|_{\rho-\delta} \leq n\delta^{-1} \| f \|_\rho$. □

**Lemma 3.1.** Let $\alpha_z(dz) = c(z) \cdot dz = \sum_{j=1}^{2n} c_j(z) dz_j$ be differential 1-form, $\Omega$ be differential 2-form given by $\Omega = d\alpha$ on $\mathbf{U}^{2n}$ and $K : \mathbf{T}^n \to \mathbf{U}^{2n}$ be a $C^\infty$ map. Let $a(\theta) = DK(\theta)^T c(K(\theta))$ and $L(\theta) = Da(\theta)^T - Da(\theta)$. Then the followings hold;

1) The representation matrix of differential 2-form $\Omega$ is given by $Dc(\theta)^T - Dc(\theta)$.

2) Pull-back $K^*\alpha$ of differential 1-form by $\alpha$ is given by $K^*\alpha = \sum_{i=1}^n a_i(\theta) d\theta_i$.

3) The representation matrix of the pull-back $K^*\Omega$ of differential 1-form $\Omega$ by $K$ is given $L(\theta) = Da(\theta)^T - Da(\theta)$.

**Definition 3.1.** ([Ito 1990]) Let $(M, \omega)$ be a $2n$-dimensional symplectic manifold and $N$ an $n$-dimensional submanifold of $M$. $N$ is called Lagrangian if $\omega(\xi, \eta) = 0, (\forall \xi, \eta \in T_p N)$ for any $p \in N$. □

If $n$-dimensional torus $\mathcal{T} = K(\mathbf{T}^n) \subset \mathbf{U}^{2n}$ is invariant with the generalized Hamiltonian system (2.1) and $\omega(\xi, \eta) = 0, (\forall \xi, \eta \in T_z\mathbf{U}^{2n}, \forall z \in \mathcal{T})$, then $\mathcal{T}$ is called as Lagrangian invariant torus of the generalized Hamiltonian system (2.1).

On the other hand let $L(z), z \in \mathcal{T}$ be the representation matrix of $\omega | \mathcal{T}$. If for any $\varepsilon > 0$ and any $d \in (0, \rho)$ there exist some $\delta > 0$ such that if $\| e \|_\rho < \delta$ then $\| L \|_{\rho-d} < \varepsilon$ holds, $\mathcal{T}$ is called by approximate Lagrangian invariant torus.

**Theorem 3.1.** If $C^r$ embedding $K : \mathbf{T}^n \to \mathbf{U}^{2n}$ is the solution of (2.16) with Diophantine vector $\omega \in \mathbf{R}^n$, $\mathcal{T} = K(\mathbf{T}^n)$ is Lagrangian invariant torus for the generalized Hamiltonian system (2.1).

**Proof.** From corollary 2.1, $\mathcal{T} = K(\mathbf{T}^n)$ is the invariant torus of the generalized Hamiltonian system (2.1). First let us show the representation matrix $L(\theta)$ of $K^*\Omega$ is constant function. For $v \in \mathbf{R}^n$, let us define map $T_v : \mathbf{T}^n \to \mathbf{T}^n$ by $T_v(\theta) = \theta + v$.

From corollary 2.1 we have

$$\varphi_t \circ K = K \circ T_{\omega t}, \ (\forall t \in \mathbf{R}). \tag{3.1}$$

Then from lemma 2.1 and (3.1), we obtain

$$K^*\Omega = K^*(\varphi_t^*\Omega) = (\varphi_t \circ K)^*\Omega = (K \circ T_{\omega t})^*\Omega. \tag{3.2}$$

Let us let $\pi_t : \mathbf{T}^n \to \mathbf{T}^n$; $\pi_t(\theta) = \theta + \omega t, (t \in \mathbf{R})$ and consider measure preserving dynamical system $(\mathbf{T}^n, \mu, \pi)$ of Kroneker flow. Since $\omega \in \mathbf{R}^n$ is rational independent, measure preserving dynamical system $(\mathbf{T}^n, \mu, \pi)$ is ergodic (See [Arnold- Avez 1968]). For any $\xi, \eta \in \mathbf{R}^n$ we let $f(\theta) = (K^*\Omega)_\theta(\xi, \eta), (\theta \in \mathbf{T}^n)$. Then since

$$((K \circ T_{\omega t})^*\Omega)_\theta(\xi, \eta) = (K^*\Omega)_{\theta+\omega t}(\xi, \eta) = f(\theta + \omega t) = f \circ T_{\omega t}(\theta) = f \circ \pi_t(\theta),$$

we obtain

$$f = f \circ \pi_t, \quad (\forall t \in \mathbf{R}) . \tag{3.3}$$

Since $(\mathbf{T}^n, \mu, \pi)$ is ergodic, $\pi$-invariant function $f$ is constant almost everywhere ([Arnold-Avez 1968]). Moreover since $f$ is continuous, $f$ is constant function. Therefore $L(\theta)$ is constant function.

From lemma 3.1 the representation matrix $L(\theta)$ of $K^*\Omega$ is written by $L(\theta) = Dc(\theta)^T - Dc(\theta)$. And from $\int_{\mathbf{T}^n} Dc(\theta)d\theta = O$, we have $\int_{\mathbf{T}^n} L(\theta)d\theta = O$. Since $L(\theta)$ is constant function, we obtain $L(\theta) \equiv O$. Thus $K^*\Omega = 0$ holds.

Now let us take $\Omega$ as symplectic structure on $\mathcal{T}$ and $K$ is $C^\infty$ diffeomorphism from $\mathbf{T}^n$ to $\mathcal{T} = K(\mathbf{T}^n)$. Taking any $z \in \mathcal{T} = K(\mathbf{T}^n)$ and any $\xi, \eta \in T_z\mathcal{T}$, there is a unique $\theta \in \mathbf{T}^n$ such that $z = K(\theta)$ and there are unique $\xi', \eta' \in T_\theta\mathbf{T}^n$ such that $DK(\theta)\xi' = \xi$, $DK(\theta)\eta' = \eta$. Then since

$$\Omega(\xi, \eta) = \Omega(DK(\theta)\xi', DK(\theta)\eta') = K^*\Omega(\xi', \eta') = 0,$$

$n$-dimensional torus $\mathcal{T} = K(\mathbf{T}^n)$ is Lagrangian invariant torus of (2.1). □

**Lemma 3.2.** Given $\omega \in D_n(\gamma, \sigma)$, $h \in \mathcal{P}(\overline{U}_\rho, \mathbf{R})$, we assume that $<h> = \int_{\mathbf{T}^n} h(\theta)d\theta$ is zero. Then there is a unique real analytic solution $v : U_\rho \to \mathbf{R}^{2n}$ of the first order partial differential equation

$$\partial_\omega v(\theta) = h(\theta) \tag{3.4}$$

such that one-periodic in all its variables and its average is 0. Moreover the following holds;

$$\| v \|_{\rho-\delta} \le \mu \gamma^{-1} \delta^{-\sigma} \| h \|_\rho , \tag{3.5}$$

where $\mu$ is constant depending only on $n$ and $\sigma$.

**Proof.** See [Rüssmann 1975]. By [Rüssmann 1975], in the case of setting $n$-torus as $\mathbf{T}^n = \mathbf{R}^n / 2\pi\mathbf{Z}^n$, $\mu$ is given by $\mu = 2^{n/2-\sigma-1} 3^{n/2+1} \sigma^{1/2} \pi \sqrt{\Gamma(2\sigma)}$.

**Definition 3.2.** Let us assume $\omega \in D_n(\gamma, \sigma)$. Let $K : U_\rho \to \mathbf{U}^{2n}$ be a real analytic map. The error function of $K$ is defined by

$$e(\theta) = B(K(\theta))\nabla H(K(\theta)) - \partial_\omega K(\theta), \quad (\theta \in U_\rho) . \tag{3.6}$$

**Theorem 3.2.** Let $\alpha_z = c(z) \cdot dz = \sum_{j=1}^{2n} c_j(z)dz_j$, $\Omega = d\alpha$. Then the followings hold:

1) The each components of vector $DK(\theta)^T B(K(\theta))^{-1} e(\theta)$ have zero average:

$$< DK(\theta)^T B(K(\theta))^{-1} e(\theta) > = 0 .$$

2) If $e = 0$ thus if $K$ is solution of (2.4), then there exist $b : \mathbf{T}^n \to \mathbf{R}$ and $a_0 \in \mathbf{R}^n$ such that

$$K^*\alpha = db + a_0 d\theta . \tag{3.7}$$

3) There exists $\tilde{b} : \mathbf{T}^n \to \mathbf{R}$, $g : \mathbf{T}^n \to \mathbf{R}^n$ and $a_0 \in \mathbf{R}^n$ such that

$$K^*\alpha = d\tilde{b} + a_0 d\theta + \sum_{j=1}^{n} g_j(\theta)d\theta_j . \tag{3.8}$$

Moreover

$$\partial_\omega g(\theta) = DK(\theta)^T B(K(\theta))^{-1} e(\theta) . \tag{3.9}$$

4) There exists a constant $c_1 \geq \max\{1, \mu\}$ depending only on $n, \sigma, \|DK\|_\rho, \|B^{-1}\|_{\mathcal{B}_r}$ which is denoted by polynomial with positive coefficients with $\|DK\|_\rho, \|B^{-1}\|_{\mathcal{B}_r}$ such that the representation matrix $L(\theta)$ of $K^*\Omega$ satisfies

$$\|L\|_{\rho-2\delta} \leq c_1 \gamma^{-1} \delta^{-(\sigma+1)} \|e\|_\rho. \tag{3.10}$$

Especially, if $e = 0$, then $L = 0$.

**Proof.** 1): From 2) of lemma 3.1, the pull-back of $\alpha_z = c(z) \cdot dz = \sum_{j=1}^{2n} c_j(z) dz_j$ under $K$ is given

$$(K^*\alpha)_\theta(d\theta) = \sum_{i=1}^{n} a_i(\theta) d\theta_i, \quad a_i(\theta) = (\frac{\partial K(\theta)}{\partial \theta_i})^T c(K(\theta)).$$

And From 1) of lemma 3.1, the representation matrix of $\Omega = d\alpha = d(\sum_{i=1}^{2n} c_i(z) dz_i)$ is given by $J(z) = Dc(z)^T - Dc(z)$. Now let us calculate $\partial_\omega a_j(\theta)$.

$$\partial_\omega a_j(\theta) = (\partial_\omega \frac{\partial K(\theta)}{\partial \theta_j})^T c(K(\theta)) + (\frac{\partial K(\theta)}{\partial \theta_j})^T \partial_\omega c(K(\theta)) =$$

$$= (\frac{\partial}{\partial \theta_j} \partial_\omega K(\theta))^T c(K(\theta)) + (\frac{\partial K(\theta)}{\partial \theta_j})^T Dc(K(\theta)) \partial_\omega K(\theta) =$$

$$= \{\frac{\partial}{\partial \theta_j}[B(K(\theta))\nabla H(K(\theta))]\}^T c(K(\theta)) - \{\frac{\partial}{\partial \theta_j} e(\theta)\}^T c(K(\theta)) +$$

$$+ (\frac{\partial K(\theta)}{\partial \theta_j})^T Dc(K(\theta))[B(K(\theta))\nabla H(K(\theta)) - e(\theta)] =$$

$$= \frac{\partial}{\partial \theta_j}\{[B(K(\theta))\nabla H(K(\theta))]^T c(K(\theta))\} - [B(K(\theta))\nabla H(K(\theta))]^T \frac{\partial}{\partial \theta_j} c(K(\theta)) -$$

$$- \frac{\partial}{\partial \theta_j}[e(\theta)^T c(K(\theta))] + e(\theta)^T \frac{\partial}{\partial \theta_j} c(K(\theta)) +$$

$$+ (\frac{\partial K(\theta)}{\partial \theta_j})^T Dc(K(\theta))[B(K(\theta))\nabla H(K(\theta)) - e(\theta)] =$$

$$= \frac{\partial}{\partial \theta_j}\{[B(K(\theta))\nabla H(K(\theta))]^T c(K(\theta)) - e(\theta)^T c(K(\theta))\} +$$

$$+ (\frac{\partial K(\theta)}{\partial \theta_j})^T Dc(K(\theta))[B(K(\theta))\nabla H(K(\theta))] -$$

$$- [B(K(\theta))\nabla H(K(\theta))]^T Dc(K(\theta))\frac{\partial K(\theta)}{\partial \theta_j} +$$

$$+ e(\theta)^T Dc(K(\theta))\frac{\partial K(\theta)}{\partial \theta_j} - (\frac{\partial K(\theta)}{\partial \theta_j})^T Dc(K(\theta))e(\theta)$$

$$\tag{3.11}$$

Since
$$[B(K(\theta))\nabla H(K(\theta))]^T Dc(K(\theta))\frac{\partial K(\theta)}{\partial \theta_j} =$$
$$= \left\{[B(K(\theta))\nabla H(K(\theta))]^T Dc(K(\theta))\frac{\partial K(\theta)}{\partial \theta_j}\right\}^T = \frac{\partial K(\theta)}{\partial \theta_j}^T Dc(K(\theta))^T [B(K(\theta))\nabla H(K(\theta))],$$

we have
$$(\frac{\partial K(\theta)}{\partial \theta_j})^T Dc(K(\theta))[B(K(\theta))\nabla H(K(\theta))] - [B(K(\theta))\nabla H(K(\theta))]^T Dc(K(\theta))\frac{\partial K(\theta)}{\partial \theta_j} =$$
$$= (\frac{\partial K(\theta)}{\partial \theta_j})^T [Dc(K(\theta)) - Dc(K(\theta))^T] B(K(\theta))\nabla H(K(\theta)).$$

(3.10)

In the same way, we obtain
$$e(\theta)^T Dc(K(\theta))\frac{\partial K(\theta)}{\partial \theta_j} - (\frac{\partial K(\theta)}{\partial \theta_j})^T Dc(K(\theta))e(\theta) = (\frac{\partial K(\theta)}{\partial \theta_j})^T [Dc(K(\theta))^T - Dc(K(\theta))]e(\theta). \quad (3.13)$$

From (3.11), (3.10) and (3.13), we have
$$\partial_\omega a_j(\theta) = \frac{\partial}{\partial \theta_j}\{[B(K(\theta))\nabla H(K(\theta))]^T c(K(\theta)) - e(\theta)^T c(K(\theta))\} -$$
$$- (\frac{\partial K(\theta)}{\partial \theta_j})^T [Dc(K(\theta))^T - Dc(K(\theta))]B(K(\theta))\nabla H(K(\theta)) +$$
$$+ (\frac{\partial K(\theta)}{\partial \theta_j})^T [Dc(K(\theta))^T - Dc(K(\theta))]e(\theta)$$
$$= \frac{\partial}{\partial \theta_j}\{[B(K(\theta))\nabla H(K(\theta))]^T c(K(\theta)) - e(\theta)^T c(K(\theta))\} -$$
$$- (\frac{\partial K(\theta)}{\partial \theta_j})^T B(K(\theta))^{-1} B(K(\theta))\nabla H(K(\theta)) + (\frac{\partial K(\theta)}{\partial \theta_j})^T B(K(\theta))^{-1} e(\theta) =$$
$$= \frac{\partial}{\partial \theta_j}\{[B(K(\theta))\nabla H(K(\theta))]^T c(K(\theta)) - e(\theta)^T c(K(\theta))\}$$
$$- (\frac{\partial K(\theta)}{\partial \theta_j})^T \nabla H(K(\theta)) + (\frac{\partial K(\theta)}{\partial \theta_j})^T B(K(\theta))^{-1} e(\theta) =$$
$$= \frac{\partial}{\partial \theta_j}\{[B(K(\theta))\nabla H(K(\theta))]^T c(K(\theta)) - H(K(\theta)) - e(\theta)^T c(K(\theta))\} +$$
$$+ (\frac{\partial K(\theta)}{\partial \theta_j})^T B(K(\theta))^{-1} e(\theta).$$

Thus we obtain
$$\partial_\omega a_j(\theta) = \frac{\partial}{\partial \theta_j}\{[B(K(\theta))\nabla H(K(\theta))]^T c(K(\theta)) - H(K(\theta)) - e(\theta)^T c(K(\theta))\} +$$
$$+ (\frac{\partial K(\theta)}{\partial \theta_j})^T B(K(\theta))^{-1} e(\theta).$$

(3.14)

Since the average of left-hand side of (3.14) and the average of the first term equal to zero, we have

$$< DK(\theta)^T B(K(\theta))^{-1} e(\theta) > = 0 \qquad (3.15)$$

On the other hand, from lemma 3.2, for any $j = 1, \cdots, n$, there exists real analytic function with zero average such that

$$\partial_\omega g = DK(\theta)^T B(K(\theta))^{-1} e(\theta) \qquad (3.16)$$

holds.

2): In the first item of the right-hand side of (3.14) we let

$$h_1 = [B(K(\theta))\nabla H(K(\theta))]^T c(K(\theta)) - H(K(\theta))$$

and

$$h_2 = e(\theta)^T c(K(\theta)).$$

Then by lemma 3.2, there exist real analytic functions $b_1, b_2 : \mathbf{T}^n \to \mathbf{R}$ with the zero average such that

$$\partial_\omega b_1 = h_1 - <h_1>, \quad \partial_\omega b_2 = h_2 - <h_2>.$$

If $K$ is the solution of (2.16) and $e(\theta) \equiv 0$, then $g \equiv 0$ and $b_2 \equiv 0$. Therefore

$$\partial_\omega a_j(\theta) = \frac{\partial h_1}{\partial \theta_j} = \frac{\partial}{\partial \theta_j}\{\partial_\omega b_1 - <h_1>\} = \partial_\omega \frac{\partial b_1}{\partial \theta_j}.$$

Thus we obtain

$$\partial_\omega \left( a_j(\theta) - \frac{\partial b_1}{\partial \theta_j} \right) = 0.$$

Considering $<\frac{\partial b_1}{\partial \theta_j}> = 0$, from lemma 3.2 we have

$$a_j(\theta) - \frac{\partial h_1}{\partial \theta_j} = <a_j(\theta) - \frac{\partial h_1}{\partial \theta_j}> = <a_j(\theta)>.$$

Therefore

$$K^*\alpha = \sum_{j=1}^n a_j(\theta)d\theta_j = \sum_{j=1}^n (\frac{\partial b_1}{\partial \theta_j} + <a_j(\theta)>)d\theta_j = \sum_{j=1}^n \frac{\partial b_1}{\partial \theta_j}d\theta_j + \sum_{j=1}^n <a_j(\theta)> d\theta_j.$$

Taking

$$b = b_1, \quad a_0 = (<a_1(\theta)>, \cdots, <a_n(\theta)>)^T,$$

we have $K^*\alpha = db + a_0 \cdot d\theta$. Thus 2) is proved.

3): Since

$$\partial_\omega a_j(\theta) = \frac{\partial}{\partial \theta_j}\{h_1 - h_2\} + \partial_\omega g_j = \frac{\partial}{\partial \theta_j}\{\partial_\omega b_1 + <h_1> - [\partial_\omega b_2 + <h_2>]\} + \partial_\omega g_j =$$

$$= \partial_\omega \{\frac{\partial b_1}{\partial \theta_j} - \frac{\partial b_2}{\partial \theta_j} + g_j\},$$

we have

$$\partial_\omega [a_j(\theta) - \frac{\partial b_1}{\partial \theta_j} + \frac{\partial b_2}{\partial \theta_j} - g_j(\theta)] = 0.$$

From $<\frac{\partial b_1}{\partial \theta_j}> = 0$, $<\frac{\partial b_2}{\partial \theta_j}> = 0$, $<g_j> = 0$, we have

$$a_j(\theta) - \frac{\partial b_1}{\partial \theta_j} + \frac{\partial b_2}{\partial \theta_j} - g_j(\theta) = <a_j(\theta)>,$$

$$a_j(\theta) = \frac{\partial (b_1 - b_2)}{\partial \theta_j} + <a_j(\theta)> + g_j(\theta).$$

By substituting above expression into $K^*\alpha = \sum_{j=1}^n a_j(\theta) d\theta_j$, we have

$$K^*\alpha = \sum_{j=1}^n a_j(\theta) d\theta_j = \sum_{j=1}^n \frac{\partial}{\partial \theta_j}(b_1 - b_2) d\theta_j + \sum_{j=1}^n <a_j(\theta)> d\theta_j + \sum_{j=1}^n g_j(\theta) d\theta_j.$$

Thus taking $\tilde{b} = b_1 - b_2$, $a_0 = (<a_1>, \cdots, <a_n>)^T$,

$$K^*\alpha = d\tilde{b} + a_0 \cdot d\theta + g(\theta) \cdot d\theta.$$

4): From 3),

$$K^*\Omega = dK^*\alpha = d(g(\theta) \cdot d\theta)$$

holds and from lemma 3.1 we get $L(\theta) = Dg(\theta)^T - Dg(\theta)$. Therefore from 8) of proposition 3.2, we have

$$\| L(\theta) \| \leq \| Dg(\theta)^T \| + \| Dg(\theta) \| \leq n \| Dg(\theta) \| + \| Dg(\theta) \| \leq 2n \| Dg(\theta) \|.$$

From above expression, Cauchy's estimate (corollary 3.2), (3.16) and lemma 3.2, there exists constant $c_1 \geq \max\{1, \mu\}$ depending only on $n, \sigma, \| DK \|_\rho, \| B^{-1} \|_{\mathcal{B}_r}$ which is denoted by polynomial with positive coefficients with $\| DK \|_\rho, \| B^{-1} \|_{\mathcal{B}_r}$ such that

$$\| L \|_{\rho - \delta} \leq c_1 \gamma^{-1} \delta^{-(\sigma+1)} \| e \|_\rho. \quad \square$$

## 4. Variable Transformation to approximate reducible type

In this section we reduce the invariant torus equation (2.16): $\partial_\omega K(\theta) = B(K(\theta)) \nabla H(K(\theta))$ to an approximately solvable form by a variable transformation. Let $n \in \mathbf{N} : n \geq 2$, $\gamma \in (0, 1)$ and $\omega \in D_n(\gamma, \sigma)$. Let's define $F : \mathcal{P}_\rho = \mathcal{P}(\overline{U}_\rho, \mathbf{U}^{2n}) \to \mathcal{P}(U_\rho, \mathbf{R}^{2n})$ by

$$(F(K))(\theta) = B(K(\theta)) \nabla H(K(\theta)) - \partial_\omega K(\theta). \tag{4.1}$$

The problem to solve (2.16) is equivalent to the problem to solve nonlinear equation with unknown $K$:

$$F(K) = 0. \tag{4.2}$$

In order to solve (4.2) by quasi-Newton method, let's drive the linearized equation

$$DF(K)\Delta = -F(K). \tag{4.3}$$

We define $F_1 : \mathcal{P}_\rho \to \mathcal{P}(U_\rho, \mathbf{R}^{2n})$, $F_2 : \mathcal{P}_\rho \to \mathcal{P}(U_\rho, \mathbf{R}^{2n})$ by

$$(F_1(K))(\theta) = B(K(\theta))\nabla H(K(\theta)),$$
$$(F_2(K))(\theta) = \partial_\omega K(\theta)$$

respectively. And we define continuous bilinear map $T: L(\mathbf{R}^{2n}, \mathbf{R}^{2n}) \times \mathbf{R}^{2n} \to \mathbf{R}^{2n}$ by $T(S, z) = Sz$. Then since

$$B(K(\theta))\nabla H(K(\theta)) = (T \circ (B \circ K, \ \nabla H \circ K))(\theta), \qquad (4.4)$$

we have

$$F_1(K) = T \circ (B \circ K, \ \nabla H \circ K). \qquad (4.5)$$

Now we define $f \in C(\mathcal{P}_\rho, \ \mathcal{P}(U_\rho, \ L(\mathbf{R}^{2n}, \mathbf{R}^{2n})))$, $g \in C(\mathcal{P}_\rho, \ \mathcal{P}(U_\rho, \ \mathbf{R}^{2n}))$ by

$$f(K) = B \circ K, \quad g(K) = \nabla H \circ K \qquad (4.6)$$

respectively. For $\xi \in C(\mathcal{P}_\rho, \ \mathcal{P}(U_\rho, \ L(\mathbf{R}^{2n}, \mathbf{R}^{2n})))$, $\eta \in C(\mathcal{P}_\rho, \ \mathcal{P}(U_\rho, \ \mathbf{R}^{2n}))$ we define $\tilde{T}(\xi, \ \eta) \in C(\mathcal{P}_\rho, \ \mathcal{P}(U_\rho, \ \mathbf{R}^{2n}))$ by $\tilde{T}(\xi, \ \eta) = T \circ (\xi, \ \eta)$. Then

$$\tilde{T}: C(\mathcal{P}_\rho, \ \mathcal{P}(U_\rho, \ L(\mathbf{R}^{2n}, \mathbf{R}^{2n}))) \times C(\mathcal{P}_\rho, \ \mathcal{P}(U_\rho, \ \mathbf{R}^{2n})) \to C(\mathcal{P}_\rho, \ \mathcal{P}(U_\rho, \ \mathbf{R}^{2n}))$$

is continuous bilinear map and thus

$$D\tilde{T}(\xi, \ \eta)(d\xi, \ d\eta) = \tilde{T}(\xi, \ d\eta) + \tilde{T}(d\xi, \ \eta),$$
$$(\forall d\xi \in C(\mathcal{P}_\rho, \ \mathcal{P}(U_\rho, \ L(\mathbf{R}^{2n}, \mathbf{R}^{2n}))), \ \forall d\eta \in C(\mathcal{P}_\rho, \ \mathcal{P}(U_\rho, \ \mathbf{R}^{2n}))). \quad (4.7)$$

Since

$$F_1(K) = T \circ (B \circ K, \ \nabla H \circ K) = T \circ (f(K), \ g(K)) =$$
$$= \tilde{T}(f(K), \ g(K)) = [\tilde{T} \circ (f, \ g)](K),$$

we have

$$F_1 = \tilde{T} \circ (f, \ g). \qquad (4.8)$$

From the differentiation formula for composition map,

$$DF_1(K) = D\tilde{T}(f(K), \ g(K)) \circ (Df(K), \ Dg(K)) \qquad (4.9)$$

holds. Using (4.7) we have

$$DF_1(K)\Delta = D\tilde{T}(f(K), \ g(K))(Df(K)\Delta, \ Dg(K)\Delta) =$$
$$= \tilde{T}(f(K), \ Dg(K)\Delta) + \tilde{T}(Df(K)\Delta, \ g(K)) = \qquad (4.10)$$
$$= T \circ (f(K), \ Dg(K)\Delta) + T \circ (Df(K)\Delta, \ g(K)).$$

For maps $S: U_\rho \to L(\mathbf{R}^{2n}, \mathbf{R}^{2n})$ and $u: U_\rho \to \mathbf{R}^{2n}$, we define map $S.u: U_\rho \to \mathbf{R}^{2n}$ by $S.u(\theta) = S(\theta)u(\theta)$. Since $Df(K)\Delta = (DB \circ K).\Delta$, we have

$$(Df(K)\Delta)(\theta) = ((DB \circ K).\Delta)(\theta) = DB(K(\theta))\Delta(\theta) \qquad (4.11)$$

from composition theorem ([Irwin 1980]). In the same way

$$(Dg(K)\Delta)(\theta) = ((D\nabla H \circ K).\Delta)(\theta) = D\nabla H(K(\theta))\Delta(\theta) \qquad (4.12)$$

holds. If we remind (4.11) and (4.12) for (4.10),

$$[DF_1(K)\Delta](\theta) =$$
$$= T(f(K)(\theta), \ Dg(K)\Delta(\theta)) + T(Df(K)\Delta(\theta), \ g(K)(\theta)) =$$
$$= T(B(K(\theta)), \ D\nabla H(K(\theta))\Delta(\theta)) +$$
$$+ T(DB(K(\theta))\Delta(\theta), \ \nabla H(K(\theta))) =$$
$$= B(K(\theta))D\nabla H(K(\theta))\Delta(\theta) + DB(K(\theta))\Delta(\theta)\nabla H(K(\theta)).$$

Thus
$$(DF_1(K)\Delta)(\theta) = DB(K(\theta))\Delta(\theta)\nabla H(K(\theta)) + B(K(\theta))D\nabla H(K(\theta))\Delta(\theta). \quad (4.13)$$
On the other hand, since $F_2(K) = \partial_\omega K$ is continuous linear map with $K$, we obtain
$$(DF_2(K)\Delta)(\theta) = \partial_\omega \Delta(\theta). \quad (4.14)$$
Therefore from (4.13) and (4.14), the linearized equation (4.3) is given by
$$DB(K(\theta))\Delta(\theta)\nabla H(K(\theta)) + B(K(\theta))D\nabla H(K(\theta))\Delta(\theta) - \partial_\omega \Delta(\theta) =$$
$$= \partial_\omega K(\theta) - B(K(\theta))\nabla H(K(\theta)). \quad (4.15)$$
Then, since
$$\left[DB(K(\theta))\Delta(\theta)\nabla H(K(\theta))\right]_i = \sum_{j=1}^{2n}\sum_{k=1}^{2n}\frac{\partial b_{ij}}{\partial z_k}(K(\theta))\Delta_k(\theta)\frac{\partial H}{\partial z_j}(K(\theta)) =$$
$$= \sum_{k=1}^{2n}\left[\sum_{j=1}^{2n}\frac{\partial b_{ij}}{\partial z_k}(K(\theta))\frac{\partial H}{\partial z_j}(K(\theta))\right]\Delta_k(\theta),$$
taking $2n \times 2n$ matrix-valued function $\Phi(z)$ whose $(i, k)$-element is given by
$$\Phi_{ik}(z) = \sum_{j=1}^{2n}\frac{\partial b_{ij}}{\partial z_k}(z)\frac{\partial H}{\partial z_j}(z), \quad z \in \mathbf{U}^{2n}, \quad (4.16)$$
we have
$$DB(z)h\nabla H(z) = \Phi(z)h, \quad (h \in \mathbf{R}^{2n}), \quad (4.17)$$
$$DB(K(\theta))\Delta(\theta)\nabla H(K(\theta)) = \Phi(K(\theta))\Delta(\theta).$$
Therefore, if we define
$$\Psi(z) = D\nabla H(z), \quad z \in \mathbf{U}^{2n}, \quad (4.18)$$
$$A(\theta) = \Phi(K(\theta)) + B(K(\theta))\Psi(K(\theta)), \quad (4.19)$$
$$e(\theta) = B(K(\theta))\nabla H(K(\theta)) - \partial_\omega K(\theta), \quad (\theta \in U_\rho), \quad (4.20)$$
then the linearized equation (4.3) for quasi-Newton method for (4.2) is written as
$$A(\theta)\Delta(\theta) - \partial_\omega \Delta(\theta) = -e(\theta). \quad (4.21)$$
Let's remark that
$$A(\theta)h = DB(K(\theta))h\nabla H(K(\theta)) + B(K(\theta))D\nabla H(K(\theta))h. \quad (4.19)'$$
**Definition 4.1.** Suppose that $K \in \mathcal{P}_\rho$ satisfies the followings:

(ND1) $Y(\theta) = K(\theta)^T K(\theta), \; (\theta \in U_\rho)$ (4.22)

is non-degenerate matrix. We denote the inverse matrix of $Y(\theta)$ by $N(\theta)$.

(ND2) The average $<S^0>$ of $n \times n$ matrix-valued function
$$S^0(\theta) := N(\theta)DK(\theta)^T\{A(\theta)B(K(\theta)) - B(K(\theta))A(\theta) -$$
$$- DB(K(\theta))B(K(\theta))\nabla H(K(\theta)) +$$
$$+ B(K(\theta))DK(\theta)N(\theta)DK(\theta)^T[A(\theta) + A(\theta)^T]\}DK(\theta)N(\theta)$$
(4.23)

is non-degenerate matrix. Then we call that $K \in \mathcal{P}_\rho$ is non-degenerate.

We denote the set of non-degenerate $K \in \mathcal{P}_\rho$ by $\mathcal{ND}(\rho)$.

**Lemma 4.1.** For

$$f(z) := B(z)\nabla H(z) \tag{4.24}$$

we define
$$a(z) = Df(z). \tag{4.25}$$

Then
$$a(z)B(z) + B(z)a(z)^T = DB(z)B(z)\nabla H(z) \tag{4.26}$$

holds.

**Proof.** We denote the non-prolongable solution of (2.1): $\dfrac{dz}{dt} = B(z)\nabla H(z)$ satisfying $y(0) = z$ by $\varphi_t(z), t \in I(x)$. From

$$\frac{d}{dt}\varphi_t(z) = f(\varphi_t(z)), \tag{4.27}$$

$$\varphi_0(z) = z \tag{4.27$_0$}$$

we have
$$\left.\frac{d}{dt}\varphi_t(z)\right|_{t=0} = f(z) = B(z)\nabla H(z). \tag{4.28}$$

From the differentiability theorem of the solution of ordinary differential equation with respect to initial condition, we have

$$\frac{d}{dt}D\varphi_t(z) = Df(\varphi_t(z))D\varphi_t(z) = a(\varphi_t(z))D\varphi_t(z), \tag{4.29}$$

$$D\varphi_t(z)|_{t=0} = I_{2n} \quad \text{(where } I_{2n} \text{ is } 2n \text{ dimensional unit matrix )}. \tag{4.29$_0$}$$

Since $\varphi_t$ is Poisson map(See [Hairer, 2006]),

$$D\varphi_t(z)B(z)D\varphi_t(z)^T = B(\varphi_t(z)) \tag{4.30}$$

holds. Differentiating the both sides of (4.30) with respect to $t$, we have

$$\frac{d}{dt}[D\varphi_t(z)]B(z)D\varphi_t(z)^T + D\varphi_t(z)B(z)\frac{d}{dt}[D\varphi_t(z)]^T = DB(\varphi_t(z))\frac{d}{dt}[\varphi_t(z)].$$

From (4.27) and (4.29), we get

$$[a(\varphi_t(z))D\varphi_t(z)]B(z)D\varphi_t(z)^T + D\varphi_t(z)B(z)[a(\varphi_t(z))D\varphi_t(z)]^T =$$
$$= DB(\varphi_t(z))B(\varphi_t(z))\nabla H(\varphi_t(z)).$$

Taking $t = 0$ in the last relation, from (4.27)$_0$ and (4.29)$_0$ we have

$$a(z)B(z) + B(z)a(z)^T = DB(z)B(z)\nabla H(z). \square$$

**Remark 4.1.**

1) $a(z)h = B(z)D\nabla H(z)h + DB(z)h\nabla H(z), \ (h \in \mathbf{R}^{2n})$

2) $a(K(\theta))h = B(K(\theta))D\nabla H(K(\theta))h + DB(K(\theta))h\nabla H(K(\theta)) = A(\theta)h, \ (h \in \mathbf{R}^{2n}, K \in \mathcal{ND}(\rho))$

3) $A(\theta)B(K(\theta)) + B(K(\theta))A(\theta)^T = DB(K(\theta))B(K(\theta))\nabla H(K(\theta)) \tag{4.31}$

hold.

In fact, if we define $T: L(\mathbf{R}^{2n}, \mathbf{R}^{2n}) \times \mathbf{R}^{2n} \to \mathbf{R}^{2n}; T(S, h) = Sh$, then we have
$$f(z) = B(z)\nabla H(z) = T(B(z), \nabla H(z)).$$

Applying the differentiation formula for continuous bilinear map, we have

$$a(z)h = Df(z)h = [DT(B(z), \nabla H(z)) \circ (DB(z), D\nabla H(z))]h =$$
$$= DT(B(z), \nabla H(z))(DB(z)h, D\nabla H(z)h) =$$
$$= T(DB(z)h, \nabla H(z)) + T(B(z), D\nabla H(z)h) =$$
$$= DB(z)h\nabla H(z) + B(z)D\nabla H(z)h.$$

Substituting $h = K(\theta)$, $\theta \in U_\rho$ into the both sides of the last relation, we get

$$a(z)h = Df(z)h = [DT(B(z), \nabla H(z)) \circ (DB(z), D\nabla H(z))]h =$$
$$= DT(B(z), \nabla H(z))(DB(z)h, D\nabla H(z)h) =$$
$$= T(DB(z)h, \nabla H(z)) + T(B(z), D\nabla H(z)h) =$$
$$= DB(z)h\nabla H(z) + B(z)D\nabla H(z)h.$$

Substituting $z = K(\theta)$ into the both sides of (4.26): $a(z)B(z) + B(z)a(z)^T = DB(z)B(z)\nabla H(z)$ and taking $a(K(\theta)) = A(\theta)$ yield

$$A(\theta)B(K(\theta)) + B(K(\theta))A(\theta)^T = DB(K(\theta))B(K(\theta))\nabla H(K(\theta)). \quad \square$$

Let us fix $K \in \mathcal{ND}(\rho)$ in the following discussions.

**Lemma 4.2.** If $K$ is the solution of (2.16), then $\left\{ \dfrac{\partial K(\theta)}{\partial \theta_j}, \ B(K(\theta))\dfrac{\partial K(\theta)}{\partial \theta_j} \right\}_{j=1,\cdots,n}$ forms the basis of $T_{K(\theta)}\mathbf{U}^{2n} \cong \mathbf{R}^{2n}$.

**Proof.** The representation matrix of $K^*\Omega$ is given by $L(\theta) = DK(\theta)^T B(K(\theta))^{-1} DK(\theta)$ and there exists a constant $c_1 \geq \max\{1, \ \mu\}$ depending on $n, \sigma$ such that the representation matrix $L(\theta)$ of $K^*\Omega$ satisfies (3.10): $\|L\|_{\rho-2\delta} \leq c_1 \gamma^{-1}\delta^{-(\sigma+1)}\|e\|_\rho$. Thus if $K(\theta)$ is a solution of (2.16), then $DK(\theta)^T B(K(\theta))^{-1} DK(\theta) = 0$. For any $i, j = 1, \cdots, n$, vectors $A_i = \dfrac{\partial K(\theta)}{\partial \theta_i} \in \mathbf{R}^{2n}$ is orthogonal to $C_j = B(K(\theta))^{-1}\dfrac{\partial K(\theta)}{\partial \theta_j} \in \mathbf{R}^{2n}$. Let $E \subset \mathbf{R}^{2n}$ is the linear subspace of $\mathbf{R}^{2n}$ that $A_1, \cdots, A_n$ span. Then $C_j \in E^\perp$ ($j = 1, \cdots, n$). From the condition $\det DK(\theta)^T DK(\theta) \neq 0$, $\{A_j\}_{j=1}^n$ are linearly independent and since $B(K(\theta)) : \mathbf{R}^{2n} \to \mathbf{R}^{2n}$ is linear isomorphism, $\{C_j\}_{j=1}^n$ are also linearly independent. Moreover from $\dim E = n$, we have $\dim E^\perp = n$. Hence $\{C_j\}_{j=1}^n$ forms a basis of $E^\perp$.

Therefore, $\left\{ \dfrac{\partial K(\theta)}{\partial \theta_j}, \ B(K(\theta))^{-1}\dfrac{\partial K(\theta)}{\partial \theta_j} \right\}_{j=1,\cdots,n}$ forms the basis of $\mathbf{R}^{2n}$.

Since $B(K(\theta)) : \mathbf{R}^{2n} \to \mathbf{R}^{2n}$ is the linear isomorphism, $\left\{ B(K(\theta))\dfrac{\partial K(\theta)}{\partial \theta_j}, \ \dfrac{\partial K(\theta)}{\partial \theta_j} \right\}_{j=1,\cdots,n}$ is also the basis of $\mathbf{R}^{2n}$. Let denote $n$-dimensional unit matrix by $I$ and $2n$-dimensional unit matrix by $I_{2n}$. Let

$$V(\theta) := \begin{pmatrix} 0 & I \\ -I & -N(\theta)^T DK(\theta)^T B(K(\theta))DK(\theta)N(\theta) \end{pmatrix}$$

$$R(\theta) = \begin{pmatrix} L(\theta) & 0 \\ 0 & 0 \end{pmatrix}$$

$$M(\theta) := (DK(\theta) \ B(K(\theta))DK(\theta)N(\theta)). \tag{4.32}$$

Then the following lemma holds.

**Lemma 4.3.** We assume $0 < \delta < \rho/2$.

1) Let $c_1$ be the constant in theorem 3.2. Then there exists a constant $c_2 \geq c_1$, depending only on $n$, $\sigma$, $\|DK\|_\rho$, $\|B^{-1}\|_{\mathcal{B}_r}$, $\|B\|_{\mathcal{B}_r}$, which is denoted by a polynomial with positive coefficients with $\|DK\|_\rho$, $\|B^{-1}\|_{\mathcal{B}_r}$, $\|B\|_{\mathcal{B}_r}$ such that: if $c_2 \gamma^{-1} \delta^{-(\sigma+1)} \|e\|_\rho \leq \frac{1}{2}$, then $M(\theta)$ is invertible and the inverse is given by

$$M(\theta)^{-1} = V(\theta)^{-1} M(\theta)^T B(K(\theta))^{-1} + M_e(\theta), \tag{4.33}$$

where

$$M_e(\theta) = -[I_{2n} + V(\theta)^{-1} R(\theta)]^{-1} V(\theta)^{-1} R(\theta) V(\theta)^{-1} M(\theta)^T B(K(\theta))^{-1}. \tag{4.34}$$

2) There exists a constant $c_3 \geq c_2$, depending only on $n$, $\sigma$, $\|DK\|_\rho$, $\|B^{-1}\|_{\mathcal{B}_r}$, $\|B\|_{\mathcal{B}_r}$, which is denoted by polynomial with positive coefficients with $\|DK\|_\rho$, $\|B^{-1}\|_{\mathcal{B}_r}$, $\|B\|_{\mathcal{B}_r}$ such that

$$\|M_e\|_{\rho-2\delta} \leq c_3 \gamma^{-1} \delta^{-(\sigma+1)} \|e\|_\rho. \tag{4.35}$$

Especially if $e = 0$ then $M_e = 0$.

**Proof.** Since $B(K(\theta))^T = -B(K(\theta))$, from $b_{ij}(z) = -b_{ji}(z)$ we have

$$B(K(\theta))^T B(K(\theta))^{-1} = -I. \tag{4.36}$$

Then

$$\begin{aligned}
M(\theta)^T & B(K(\theta))^{-1} M(\theta) = \\
&= \begin{pmatrix} DK(\theta)^T \\ N(\theta)^T DK(\theta)^T B(K(\theta))^T \end{pmatrix} B(K(\theta))^{-1} (DK(\theta) \ B(K(\theta))DK(\theta)N(\theta)) \\
&= \begin{pmatrix} DK(\theta)^T B(K(\theta))^{-1} \\ -N(\theta)^T DK(\theta)^T \end{pmatrix} (DK(\theta) \ B(K(\theta))DK(\theta)N(\theta)) \\
&= \begin{pmatrix} L(\theta) & I \\ -I & -N(\theta)^T DK(\theta)^T B(K(\theta))DK(\theta)N(\theta) \end{pmatrix} \\
&= \begin{pmatrix} L(\theta) & 0 \\ 0 & 0 \end{pmatrix} + \begin{pmatrix} 0 & I \\ -I & -N(\theta)^T DK(\theta)^T B(K(\theta))DK(\theta)N(\theta) \end{pmatrix} \\
&= R(\theta) + V(\theta).
\end{aligned} \tag{4.37}$$

Then $V(\theta)$ is non-degenerate and

$$V(\theta)^{-1} = \begin{pmatrix} -N(\theta)^T DK(\theta)^T B(K(\theta))DK(\theta)N(\theta) & -I \\ I & 0 \end{pmatrix}. \tag{4.38}$$

From definition of $R(\theta)$ and (3.10),

$$\|R(\theta)\|_{\rho-2\delta} \le \|L(\theta)\|_{\rho-2\delta} \le c_1\gamma^{-1}\delta^{-(\sigma+1)}\|e\|_\rho \tag{4.39}$$

and there exists a constant $c_2 \ge c_1$, depending only on $n$, $\sigma$, $\|DK\|_\rho$, $\|B^{-1}\|_{\mathcal{B}_r}$, $\|B\|_{\mathcal{B}_r}$, which is denoted by polynomial with positive coefficients with $\|DK\|_\rho$, $\|B^{-1}\|_{\mathcal{B}_r}$, $\|B\|_{\mathcal{B}_r}$ such that

$$\|V(\theta)^{-1}R(\theta)\|_{\rho-2\delta} \le c_2\gamma^{-1}\delta^{-(\sigma+1)}\|e\|_\rho.$$

From the assumption of the theorem,

$$\|V(\theta)^{-1}R(\theta)\|_{\rho-2\delta} \le c_2\gamma^{-1}\delta^{-(\sigma+1)}\|e\|_\rho \le \frac{1}{2}$$

and from Neumann's series theorem $(V(\theta)^{-1}R(\theta)+I_{2n})$ is non-degenerate and

$$\|(I_{2n}+V(\theta)^{-1}R(\theta))^{-1}\|_{\rho-2\delta} \le \frac{1}{1-\|V(\theta)^{-1}R(\theta)\|_{\rho-2\delta}} \le 2. \tag{4.40}$$

Multiplying (4.37) from the left by $[V(\theta)(V(\theta)^{-1}R(\theta)+I_{2n})]^{-1}$, we have

$$[V(\theta)(V(\theta)^{-1}R(\theta)+I_{2n})]^{-1}M(\theta)^T B(K(\theta))^{-1}M(\theta) = I_{2n}.$$

Considering

$$(I_{2n}+V(\theta)^{-1}R(\theta))^{-1} = I_{2n} - (I_{2n}+V(\theta)^{-1}R(\theta))^{-1}V(\theta)^{-1}R(\theta),$$

we get

$$\begin{aligned}
&[V(\theta)(V(\theta)^{-1}R(\theta)+I_{2n})]^{-1}M(\theta)^T B(K(\theta))^{-1} = \\
&= [V(\theta)(V(\theta)^{-1}R(\theta)+I_{2n})]^{-1}M(\theta)^T B(K(\theta))^{-1} = \\
&= (I_{2n}+V(\theta)^{-1}R(\theta))^{-1}V(\theta)^{-1}M(\theta)^T B(K(\theta))^{-1} = \\
&= [I_{2n} - (I_{2n}+V(\theta)^{-1}R(\theta))^{-1}V(\theta)^{-1}R(\theta)]V(\theta)^{-1}M(\theta)^T B(K(\theta))^{-1} = \\
&= V(\theta)^{-1}M(\theta)^T B(K(\theta))^{-1} - \\
&\quad - (I_{2n}+V(\theta)^{-1}R(\theta))^{-1}V(\theta)^{-1}R(\theta)V(\theta)^{-1}M(\theta)^T B(K(\theta))^{-1} = \\
&= V(\theta)^{-1}M(\theta)^T B(K(\theta))^{-1} + M_e(\theta).
\end{aligned}$$

Therefore from the relation

$$[V(\theta)(V(\theta)^{-1}R(\theta)+I_{2n})]^{-1}M(\theta)^T B(K(\theta))^{-1}M(\theta) = I_{2n}$$

$M(\theta)$ is invertible and

$$M(\theta)^{-1} = V(\theta)^{-1}M(\theta)^T B(K(\theta))^{-1} + M_e(\theta).$$

Then, there exists a constant $c_3 \ge c_2$, depending only on $n$, $\sigma$, $\|DK\|_\rho$, $\|N\|_\rho$, $\|B\|_{\mathcal{B}_r}$, $\|B^{-1}\|_{\mathcal{B}_r}$,

which is denoted by polynomial with positive coefficients with $\|DK\|_\rho, \|N\|_\rho, \|B\|_{\mathcal{B}_r} \|B^{-1}\|_{\mathcal{B}_r}$ such that

$$\|M_e\|_{\rho-2\delta} \le c_3 \gamma^{-1} \delta^{-(\sigma+1)} \|e\|_\rho . \square$$

**Lemma 4.4.** If $K(\theta)$ is the approximate solution of (2.16), then
$$A(\theta)DK(\theta) - \partial_\omega DK(\theta) = De(\theta). \tag{4.41}$$

Especially if $K(\theta)$ is the solution of (2.16), then
$$A(\theta)DK(\theta) - \partial_\omega DK(\theta) = De(\theta). \tag{4.41}$$

**Proof.** Let $T: L(\mathbf{R}^{2n}, \mathbf{R}^{2n}) \times \mathbf{R}^{2n} \to \mathbf{R}^{2n}$ be the continuous bilinear map given by
$$T(S, z) = Sz.$$

Then differentiating the right-hand side of
$$(4.4): B(K(\theta))\nabla H(K(\theta)) = (T \circ (B \circ K, \ \nabla H \circ K))(\theta),$$

from the differentiation formula for composite map
$$D(T \circ (B \circ K, \ \nabla H \circ K))(\theta) =$$
$$= DT(B(K(\theta)), \ \nabla H(K(\theta))) \circ (DB(K(\theta)) \circ DK(\theta), \ D\nabla H(K(\theta)) \circ DK(\theta)).$$

Therefore by the differentiation formula for continuous bilinear map:
$$DT(S, \ z)(dS, \ dz) = T(S, \ dz) + T(dS, \ z), \ (dS \in L(\mathbf{R}^{2n}, \mathbf{R}^{2n}), dz \in \mathbf{R}^{2n}),$$

we have
$$D(T \circ (B \circ K, \ \nabla H \circ K))(\theta)h =$$
$$= \{DT(B(K(\theta)), \ \nabla H(K(\theta)))\}(DB(K(\theta))DK(\theta)h, \ D\nabla H(K(\theta))DK(\theta)h) =$$
$$= T(B(K(\theta)), \ D\nabla H(K(\theta))DK(\theta)h) + T(DB(K(\theta))DK(\theta)h, \ \nabla H(K(\theta))) =$$
$$= B(K(\theta))D\nabla H(K(\theta))DK(\theta)h + DB(K(\theta))DK(\theta)h \nabla H(K(\theta)) =$$
$$= B(K(\theta))\Psi(K(\theta))DK(\theta)h + \Phi(K(\theta))DK(\theta)h$$

for any $h \in \mathbf{R}^n$. Therefore (4.17),(4.18),(4.19) give us
$$D(T \circ (B \circ K, \ \nabla H \circ K))(\theta) = [B(K(\theta))\Psi(K(\theta)) + \Phi(K(\theta))]DK(\theta) = A(\theta)DK(\theta).$$

Differentiating the both sides of
$$B(K(\theta))\nabla H(K(\theta)) - \partial_\omega K(\theta) = e(\theta)$$

yields
$$A(\theta)DK(\theta) - \partial_\omega DK(\theta) = De(\theta). \ \square$$

**Lemma 4.5.** For the solution $K$ of (2.16), defining $M$ by (4.32), we have
$$A(\theta)M(\theta) - \partial_\omega M(\theta) = M(\theta) \begin{pmatrix} 0 & S^0(\theta) \\ 0 & 0 \end{pmatrix}. \tag{4.43}$$

**Proof.** If we take variable transformation based on lemma 4.3 $\Delta(\theta) = M(\theta)\xi(\theta)$ to the linearized equation (4.21): $A(\theta)\Delta(\theta) - \partial_\omega \Delta(\theta) = -e(\theta)$, then we get
$$[A(\theta)M(\theta) - \partial_\omega M(\theta)]\xi(\theta) - M(\theta)\partial_\omega \xi(\theta) = -e(\theta). \tag{4.44}$$

Now let us calculate $A(\theta)M(\theta) - \partial_\omega M(\theta)$. Since $K(\theta)$ is the exact solution of (2.16), $A(\theta)DK(\theta) - \partial_\omega DK(\theta) = 0$ holes by lemma 4.4. On the other hand,

$$A(\theta)M(\theta) - \partial_\omega M(\theta) =$$
$$= (A(\theta)DK(\theta) \quad A(\theta)B(K(\theta))DK(\theta)N(\theta)) - (\partial_\omega DK(\theta) \quad \partial_\omega [B(K(\theta))DK(\theta)N(\theta)]) =$$
$$= (A(\theta)DK(\theta) - \partial_\omega DK(\theta), \quad A(\theta)B(K(\theta))DK(\theta)N(\theta) - \partial_\omega [B(K(\theta))DK(\theta)N(\theta)]).$$

Hence the first $n$ columns of $A(\theta)M(\theta) - \partial_\omega M(\theta)$ are all zero. Let us denote the last $n$ columns of $A(\theta)M(\theta) - \partial_\omega M(\theta)$ by $W_1(\theta)$. Thus

$$W_1(\theta) = A(\theta)B(K(\theta))DK(\theta)N(\theta) - \partial_\omega [B(K(\theta))DK(\theta)N(\theta)]. \tag{4.45}$$

On the other hand
$$\partial_\omega [B(K(\theta))DK(\theta)N(\theta)] =$$
$$= [\partial_\omega B(K(\theta))]DK(\theta)N(\theta) + B(K(\theta))\partial_\theta DK(\theta)N(\theta) + B(K(\theta))DK(\theta)\partial_\omega N(\theta).$$
$$\tag{4.46}$$

By taking into account (2.16): $\partial_\omega K(\theta) = B(K(\theta))\nabla H(K(\theta))$,

$$\partial_\omega (B \circ K)(\theta) = DB(K(\theta))\partial_\omega K(\theta) = DB(K(\theta))B(K(\theta))\nabla H(K(\theta)). \tag{4.47}$$

Taking $\partial_\omega$ on the both sides of $I_n = N(\theta)DK(\theta)^T DK(\theta)$ and using (4.42), we have

$$0 = \partial_\omega [N(\theta)DK(\theta)^T DK(\theta)] =$$
$$= \partial_\omega N(\theta)[DK(\theta)^T DK(\theta)] + N(\theta)[\partial_\omega DK(\theta)]^T DK(\theta) + N(\theta)DK(\theta)^T \partial_\omega DK(\theta) =$$
$$= \partial_\omega N(\theta)[DK(\theta)^T DK(\theta)] + N(\theta)DK(\theta)^T A(\theta)^T DK(\theta) +$$
$$+ N(\theta)DK(\theta)^T A(\theta)DK(\theta).$$

Multiplying the both sides of the above relation by $N(\theta) = (DK(\theta)^T DK(\theta))^{-1}$ from the right, we have

$$0 = \partial_\omega N(\theta) + N(\theta)DK(\theta)^T A(\theta)^T DK(\theta)N(\theta) +$$
$$+ N(\theta)DK(\theta)^T A(\theta)DK(\theta)N(\theta) =$$
$$= \partial_\omega N(\theta) + N(\theta)DK(\theta)^T [A(\theta) + A(\theta)^T]DK(\theta)N(\theta).$$

Thus we get
$$\partial_\omega N(\theta) = -N(\theta)DK(\theta)^T [A(\theta) + A(\theta)^T]DK(\theta)N(\theta). \tag{4.48}$$

Substituting (4.46), (4.47), (4.48) into (4.45) gives us
$$W_1(\theta) = A(\theta)B(K(\theta))DK(\theta)N(\theta)) - \partial_\omega [B(K(\theta))DK(\theta)N(\theta)] =$$
$$= A(\theta)B(K(\theta))DK(\theta)N(\theta)) - [\partial_\omega B(K(\theta))]DK(\theta)N(\theta) -$$
$$- B(K(\theta))\partial_\theta DK(\theta)N(\theta) - B(K(\theta))DK(\theta)\partial_\omega N(\theta) =$$
$$= A(\theta)B(K(\theta))DK(\theta)N(\theta)) - DB((K(\theta))B(K(\theta))\nabla H(K(\theta))DK(\theta)N(\theta) -$$
$$- B(K(\theta))A(\theta)DK(\theta)N(\theta) +$$
$$+ B(K(\theta))DK(\theta)N(\theta)DK(\theta)^T [A(\theta) + A(\theta)^T]DK(\theta)N(\theta) =$$
$$= [A(\theta)B(K(\theta)) - B(K(\theta))A(\theta)]DK(\theta)N(\theta)) -$$
$$- DB((K(\theta))B(K(\theta))\nabla H(K(\theta))DK(\theta)N(\theta) +$$
$$+ B(K(\theta))DK(\theta)N(\theta)DK(\theta)^T [A(\theta) + A(\theta)^T]DK(\theta)N(\theta)$$

and hence
$$W_1(\theta) = [A(\theta)B(K(\theta)) - B(K(\theta))A(\theta)]DK(\theta)N(\theta)) -$$

$$-DB((K(\theta))B(K(\theta))\nabla H(K(\theta))DK(\theta)N(\theta) +$$
$$+ B(K(\theta))DK(\theta)N(\theta)DK(\theta)^T[A(\theta)+A(\theta)^T]DK(\theta)N(\theta).$$
(4.49)

Multiplying relation (4.49) from the left by $DK(\theta)^T B(K(\theta))^{-1}$ and using lemma 4.1 and remark 4.1, we have

$$DK(\theta)^T B(K(\theta))^{-1}W_1(\theta) =$$
$$= DK(\theta)^T B(K(\theta))^{-1}A(\theta)B(K(\theta))DK(\theta)N(\theta) -$$
$$- DK(\theta)^T A(\theta)DK(\theta)N(\theta) -$$
$$- DK(\theta)^T B(K(\theta))^{-1}DB(K(\theta))B(K(\theta))\nabla H(K(\theta))DK(\theta)N(\theta) +$$
$$+ DK(\theta)^T[A(\theta)^T + A(\theta)]DK(\theta)N(\theta)$$
$$= DK(\theta)^T B(K(\theta))^{-1}A(\theta)B(K(\theta))DK(\theta)N(\theta) -$$
$$- DK(\theta)^T B(K(\theta))^{-1}DB(K(\theta))B(K(\theta))\nabla H(K(\theta))DK(\theta)N(\theta) +$$
$$+ DK(\theta)^T A(\theta)^T DK(\theta)N(\theta) =$$
$$= DK(\theta)^T B(K(\theta))^{-1}[A(\theta)B(K(\theta)) - DB(K(\theta))B(K(\theta))\nabla H(K(\theta)) +$$
$$+ B(K(\theta))A(\theta)^T]DK(\theta)N(\theta) = 0.$$

Thus we get
$$DK(\theta)^T B(K(\theta))^{-1}W_1(\theta) = 0. \tag{4.50}$$

Since $\left\{\dfrac{\partial K(\theta)}{\partial \theta_j}, \; B(K(\theta))\dfrac{\partial K(\theta)}{\partial \theta_j}\right\}_{j=1,\cdots,n}$ is the basis of $\mathbf{R}^{2n}$ by lemma 4.2, there are unique $n\times n$ matrices $S(\theta), T(\theta)$ such that

$$W_1(\theta) = A(\theta)M(\theta) - \partial_\omega M(\theta) = DK(\theta)S(\theta) + B(K(\theta))DK(\theta)N(\theta)T(\theta). \tag{4.51}$$

Multiplying the both sides of (4.51) by $DK(\theta)^T B(K(\theta))^{-1}$ and using (4.50),

$$DK(\theta)^T B(K(\theta))^{-1}W_1 = DK(\theta)^T B(K(\theta))^{-1}DK(\theta)S(\theta) + T(\theta).$$

And since $DK(\theta)^T B(K(\theta))^{-1}DK(\theta) = 0$ by Lagrangian property of $K$, we obtain
$$T(\theta) = 0, \tag{4.52}$$
$$W_1(\theta) = DK(\theta)S(\theta). \tag{4.53}$$

On the other hand, multiplying the both sides of (4.51) from the left by $N(\theta)DK(\theta)^T$ and using (4.52), (4.49), we have

$$S(\theta) = N(\theta)DK(\theta)^T W_1 =$$
$$= N(\theta)DK(\theta)^T\{A(\theta)B(K(\theta)) - B(K(\theta))A(\theta) - DB(K(\theta))B(K(\theta))\nabla H(K(\theta)) +$$
$$+ B(K(\theta))DK(\theta)N(\theta)DK(\theta)^T[A(\theta)^T + A(\theta)]\}DK(\theta)N(\theta) = S^0(\theta).$$

Since
$$M(\theta)\begin{pmatrix} 0 & S^0(\theta) \\ 0 & 0 \end{pmatrix} = (DK(\theta) \;\; B(K(\theta))DK(\theta)N(\theta))\begin{pmatrix} 0 & S^0(\theta) \\ 0 & 0 \end{pmatrix} =$$
$$= (0 \;\; DK(\theta)S^0(\theta)) = (0 \;\; W_1) = A(\theta)M(\theta) - \partial_\omega M(\theta)$$

by (4.53), we obtain the assertion of the lemma. □

**Lemma 4.6.** We assume that $K$ is the approximate solution of (2.16) and $c_2\gamma^{-1}\delta^{-(\sigma+1)}\|e\|_\rho \leq \frac{1}{2}$. Let

$$E(\theta) = A(\theta)M(\theta) - \partial_\omega M(\theta) - M(\theta)\begin{pmatrix} 0 & S^0(\theta) \\ 0 & 0 \end{pmatrix}. \quad (4.54)$$

Then there exist a constant $c_4 \geq c_3$, depending only on $n, \sigma, \|N\|_\rho, \|DK\|_\rho, \|B\|_{C^1,\mathcal{B}_r}, \|B^{-1}\|_{\mathcal{B}_r}$, $|H|_{C^2,\mathcal{B}_r}$, which is denoted by polynomial with positive coefficients with $\|DK\|_\rho, \|N\|_\rho, \|B\|_{\mathcal{B}_r}$, $\|B^{-1}\|_{\mathcal{B}_r}$ such that

$$\left\|M(\theta)^{-1}E(\theta)\right\|_{\rho-2\delta} \leq c_4\gamma^{-1}\delta^{-(\sigma+1)}\|e\|_\rho. \quad (4.55)$$

**Proof.** In order to obtain estimate (4.55), let us calculate $E(\theta)$. Since $K(\theta)$ is the approximation solution of (2.16),

$$B(K(\theta))\nabla H(K(\theta)) - \partial_\omega K(\theta) = e(\theta).$$

Differentiating the both sides of the above relation with respect to $\theta$, from lemma 4.4

$$A(K(\theta))DK(\theta) - \partial_\omega DK(\theta) = De(\theta).$$

Thus

$$A(\theta)M(\theta) - \partial_\omega M(\theta) = A(\theta)(DK(\theta) \quad B(K(\theta))DK(\theta)N(\theta)) - $$
$$- \partial_\omega(DK(\theta) \quad B(K(\theta))DK(\theta)N(\theta)) = $$
$$= (A(\theta)DK(\theta) - \partial_\omega DK(\theta) \quad A(\theta)B(K(\theta))DK(\theta)N(\theta) - \partial_\omega[B(K(\theta))DK(\theta)N(\theta)])$$
$$= (De(\theta) \quad A(\theta)B(K(\theta))DK(\theta)N(\theta) - \partial_\omega[B(K(\theta))DK(\theta)N(\theta)])$$
$$= (De(\theta) \quad W_1(\theta)),$$

where

$$W_1(\theta) = $$
$$= A(\theta)B(K(\theta))DK(\theta)N(\theta) - \partial_\omega[B(K(\theta))DK(\theta)N(\theta)] = $$
$$= A(\theta)B(K(\theta))DK(\theta)N(\theta) - [\partial_\omega B(K(\theta))]DK(\theta)N(\theta) - $$
$$- B(K(\theta))[\partial_\omega DK(\theta)]N(\theta) - B(K(\theta))DK(\theta)\partial_\omega N(\theta) = \quad (4.56)$$
$$= A(\theta)B(K(\theta))DK(\theta)N(\theta) - [\partial_\omega B(K(\theta))]DK(\theta)N(\theta) - $$
$$- B(K(\theta))A(\theta)DK(\theta)N(\theta) + B(K(\theta))De(\theta)N(\theta) - $$
$$- B(K(\theta))DK(\theta)\partial_\omega N(\theta).$$

Taking $\partial_\omega$ on the both sides of $I_n = N(\theta)DK(\theta)^T DK(\theta)$, we have

$$0 = \partial_\omega[N(\theta)DK(\theta)^T DK(\theta)] = $$
$$= \partial_\omega N(\theta)[DK(\theta)^T DK(\theta)] + N(\theta)[\partial_\omega DK(\theta)]^T DK(\theta) + N(\theta)DK(\theta)^T \partial_\omega DK(\theta) = $$
$$= \partial_\omega N(\theta)[DK(\theta)^T DK(\theta)] + N(\theta)[A(\theta)DK(\theta) - De(\theta)]^T DK(\theta) + $$
$$+ N(\theta)DK(\theta)^T[A(\theta)DK(\theta) - De(\theta)]$$

and multiplying the above relation from the right by $N(\theta) = (DK(\theta)^T DK(\theta))^{-1}$, we have

$$0 = \partial_\omega N(\theta) + N(\theta)DK(\theta)^T A(\theta)^T DK(\theta)N(\theta) - N(\theta)De(\theta)^T DK(\theta)N(\theta) +$$
$$+ N(\theta)DK(\theta)^T A(\theta)DK(\theta)N(\theta) - N(\theta)DK(\theta)^T De(\theta)N(\theta) =$$
$$= \partial_\omega N(\theta) + N(\theta)DK(\theta)^T [A(\theta) + A(\theta)^T]DK(\theta)N(\theta) -$$
$$- N(\theta)[DK(\theta)^T De(\theta) + De(\theta)^T DK(\theta)]N(\theta).$$

Thus we obtain
$$\partial_\omega N(\theta) = -N(\theta)DK(\theta)^T[A(\theta)+A(\theta)^T]DK(\theta)N(\theta) + \qquad (4.57)$$
$$+ N(\theta)[DK(\theta)^T De(\theta) + De(\theta)^T DK(\theta)]N(\theta).$$

Using (4.57) and
$$\partial_\omega(B(K(\theta))) = DB(K(\theta))\partial_\omega K(\theta) = DB(K(\theta))[B(K(\theta))\nabla H(K(\theta)) - e(\theta)],$$
we get
$$W_1(\theta) =$$
$$= A(\theta)B(K(\theta))DK(\theta)N(\theta) - B(K(\theta))A(\theta)DK(\theta)N(\theta) +$$
$$+ B(K(\theta))De(\theta)N(\theta) -$$
$$- DB(K(\theta))[B(K(\theta))\nabla H(K(\theta)) - e(\theta)]DK(\theta)N(\theta) +$$
$$+ B(K(\theta))DK(\theta)N(\theta)DK(\theta)^T[A(\theta) + A(\theta)^T]DK(\theta)N(\theta) -$$
$$- B(K(\theta))DK(\theta)N(\theta)[DK(\theta)^T De(\theta) + De(\theta)^T DK(\theta)]N(\theta) =$$
$$= A(\theta)B(K(\theta))DK(\theta)N(\theta) - B(K(\theta))A(\theta)DK(\theta)N(\theta) -$$
$$- DB(K(\theta))B(K(\theta))\nabla H(K(\theta))DK(\theta)N(\theta) +$$
$$+ B(K(\theta))DK(\theta)N(\theta)DK(\theta)^T[A(\theta) + A(\theta)^T]DK(\theta)N(\theta) +$$
$$+ B(K(\theta))De(\theta)N(\theta) + DB(K(\theta))e(\theta)DK(\theta)N(\theta) -$$
$$- B(K(\theta))DK(\theta)N(\theta)[DK(\theta)^T De(\theta) + De(\theta)^T DK(\theta)]N(\theta).$$
(4.58)

Multiplying the both sides of (4.58) by $DK(\theta)^T B(K(\theta))^{-1}$ and using 3) of remark 4.1, we obtain

$$DK(\theta)^T B(K(\theta))^{-1} W_1(\theta) =$$
$$= DK(\theta)^T B(K(\theta))^{-1} A(\theta) B(K(\theta)) DK(\theta) N(\theta) - DK(\theta)^T A(\theta) DK(\theta) N(\theta) -$$
$$- DK(\theta)^T B(K(\theta))^{-1} DB(K(\theta)) B(K(\theta)) \nabla H(K(\theta)) DK(\theta) N(\theta) +$$
$$+ DK(\theta)^T [A(\theta) + A(\theta)^T] DK(\theta) N(\theta) + DK(\theta)^T De(\theta) N(\theta) +$$
$$+ DK(\theta)^T B(K(\theta))^{-1} DB(K(\theta)) e(\theta) DK(\theta) N(\theta) -$$
$$- [DK(\theta)^T De(\theta) + De(\theta)^T DK(\theta)] N(\theta) =$$
$$= DK(\theta)^T B(K(\theta))^{-1} A(\theta) B(K(\theta)) DK(\theta) N(\theta) -$$
$$- DK(\theta)^T B(K(\theta))^{-1} DB(K(\theta)) B(K(\theta)) \nabla H(K(\theta)) DK(\theta) N(\theta) +$$
$$+ DK(\theta)^T A(\theta)^T DK(\theta) N(\theta) +$$
$$+ DK(\theta)^T B(K(\theta))^{-1} DB(K(\theta)) e(\theta) DK(\theta) N(\theta) - De(\theta)^T DK(\theta) N(\theta) =$$
$$= DK(\theta)^T B(K(\theta))^{-1} [A(\theta) B(K(\theta)) - DB(K(\theta)) B(K(\theta)) \nabla H(K(\theta)) +$$
$$+ B(K(\theta)) A(\theta)^T] DK(\theta) N(\theta) + DK(\theta)^T B(K(\theta))^{-1} DB(K(\theta)) e(\theta) DK(\theta) N(\theta) -$$
$$- De(\theta)^T DK(\theta) N(\theta) =$$
$$= DK(\theta)^T B(K(\theta))^{-1} DB(K(\theta)) e(\theta) DK(\theta) N(\theta) - De(\theta)^T DK(\theta) N(\theta).$$

Thus,
$$DK(\theta)^T B(K(\theta))^{-1} W_1(\theta) = DK(\theta)^T B(K(\theta))^{-1} DB(K(\theta)) e(\theta) DK(\theta) N(\theta) - De(\theta)^T DK(\theta) N(\theta). \qquad (4.59)$$

Since $\left\{ \dfrac{\partial K(\theta)}{\partial \theta_j},\ B(K(\theta)) \dfrac{\partial K(\theta)}{\partial \theta_j} \right\}_{j=1,\cdots,n}$ is the basis of $\mathbf{R}^{2n}$ when $K(\theta)$ is the solution of (2.16), if the error $e(\theta)$ is small enough, then $\left\{ \dfrac{\partial K(\theta)}{\partial \theta_j},\ B(K(\theta)) \dfrac{\partial K(\theta)}{\partial \theta_j} \right\}_{j=1,\cdots,n}$ also forms the basis of $\mathbf{R}^{2n}$.

Therefore there exists $n \times n$ matrices $S(\theta)$ and $N(\theta)$ uniquely such that
$$W_1(\theta) = DK(\theta) S(\theta) + B(K(\theta)) DK(\theta) N(\theta) T(\theta). \qquad (4.60)$$

Multiplying (4.60) from the left by $DK(\theta)^T B(K(\theta))^{-1}$,
$$DK(\theta)^T B(K(\theta))^{-1} W_1(\theta) = DK(\theta)^T B(K(\theta))^{-1} DK(\theta) S(\theta) + T(\theta) = L(\theta) S(\theta) + T(\theta). \qquad (4.61)$$

If we remind (4.59) for (4.61), then
$$T(\theta) = DK(\theta)^T B(K(\theta))^{-1} W_1(\theta) - L(\theta) S(\theta) =$$
$$= DK(\theta)^T B(K(\theta))^{-1} DB(K(\theta)) e(\theta) DK(\theta) N(\theta) - De(\theta)^T DK(\theta) N(\theta) - L(\theta) S(\theta). \qquad (4.62)$$

Multiplying relation (4.60) from the left by $N(\theta) DK(\theta)^T$,

$$N(\theta)DK(\theta)^T W_1(\theta) =$$
$$= S(\theta) + N(\theta)DK(\theta)^T B(K(\theta))DK(\theta)N(\theta)T(\theta) =$$
$$= N(\theta)DK(\theta)^T B(K(\theta))DK(\theta)N(\theta)DK(\theta)^T B(K(\theta))^{-1} \cdot$$
$$\cdot DB(K(\theta))e(\theta)DK(\theta)N(\theta) +$$
$$+ S(\theta) - N(\theta)DK(\theta)^T B(K(\theta))DK(\theta)N(\theta)De(\theta)^T DK(\theta)N(\theta) -$$
$$- N(\theta)DK(\theta)^T B(K(\theta))DK(\theta)N(\theta)L(\theta)S(\theta).$$

(4.63)

Multiplying relation (4.58) from the left by $N(\theta)DK(\theta)^T$, we have

$$N(\theta)DK(\theta)^T W_1(\theta) =$$
$$= N(\theta)DK(\theta)^T \{A(\theta)B(K(\theta)) - B(K(\theta))A(\theta) - DB(K(\theta))B(K(\theta))\nabla H(K(\theta)) +$$
$$+ B(K(\theta))DK(\theta)N(\theta)DK(\theta)^T [A(\theta) + A(\theta)^T]\}DK(\theta)N(\theta) +$$
$$+ N(\theta)DK(\theta)^T \{B(K(\theta))De(\theta) + DB(K(\theta))e(\theta)DK(\theta) -$$
$$- B(K(\theta))DK(\theta)N(\theta)[DK(\theta)^T De(\theta) + De(\theta)^T DK(\theta)]\}N(\theta).$$

(4.64)

If we compare (4.63) with (4.64) when $e = 0$, we have

$$S(\theta) = N(\theta)DK(\theta)^T \{A(\theta)B(K(\theta)) - B(K(\theta))A(\theta) -$$
$$- DB(K(\theta))B(K(\theta))\nabla H(K(\theta)) +$$
$$+ B(K(\theta))DK(\theta)N(\theta)DK(\theta)^T [A(\theta) + A(\theta)^T]\}DK(\theta)N(\theta) =$$
$$= S^0(\theta).$$

On the other hand, since

$$E(\theta) = A(\theta)M(\theta) - \partial_\omega M(\theta) - M(\theta)\begin{pmatrix} 0 & S^0(\theta) \\ 0 & 0 \end{pmatrix} = (De(\theta), \ W_1(\theta) - DK(\theta)S^0(\theta)), \qquad (4.65)$$

from (4.58) we have

$$W_1(\theta) - DK(\theta)S^0(\theta) =$$
$$= -DK(\theta)S^0(\theta) +$$
$$+ A(\theta)B(K(\theta))DK(\theta)N(\theta) - B(K(\theta))A(\theta)DK(\theta)N(\theta) -$$
$$- DB(K(\theta))B(K(\theta))\nabla H(K(\theta))DK(\theta)N(\theta) +$$
$$+ B(K(\theta))DK(\theta)N(\theta)DK(\theta)^T[A(\theta) + A(\theta)^T]DK(\theta)N(\theta) +$$
$$+ B(K(\theta))De(\theta)N(\theta) + DB(K(\theta))e(\theta)DK(\theta)N(\theta) -$$
$$- B(K(\theta))DK(\theta)N(\theta)[DK(\theta)^T De(\theta) + De(\theta)^T DK(\theta)]N(\theta).$$
(4.66)

Now, let us estimate $\|A\|_\rho$. Fix $\xi \in \mathbf{R}^{2n}$ arbitrarily such that $\|\xi\|=1$. Then from

$$\|A(\theta)\xi\| = \|DB(K(\theta))\xi \nabla H(K(\theta)) + B(K(\theta))D\nabla H(K(\theta))\xi\| \leq$$
$$\leq \|DB(K(\theta))\| \cdot \|\nabla H(K(\theta))\| + \|B(K(\theta))\| \cdot \|D\nabla H(K(\theta))\| \leq$$
$$\leq \|B\|_{C^1, \mathcal{B}_r} |H|_{C^1, \mathcal{B}_r} + \|B\|_{\mathcal{B}_r} |H|_{C^2, \mathcal{B}_r} \leq$$
$$\leq 2\|B\|_{C^1, \mathcal{B}_r} |H|_{C^2, \mathcal{B}_r},$$

we have

$$\|A\|_\rho \leq 2\|B\|_{C^1, \mathcal{B}_r} |H|_{C^2, \mathcal{B}_r}. \qquad (4.67)$$

Here let us estimate $\|S^0(\theta)\|_\rho$. Since

$$S^0(\theta) := N(\theta)DK(\theta)^T\{A(\theta)B(K(\theta)) - B(K(\theta))A(\theta) -$$
$$- DB(K(\theta))B(K(\theta))\nabla H(K(\theta)) +$$
$$+ B(K(\theta))DK(\theta)N(\theta)DK(\theta)^T[A(\theta) + A(\theta)^T]\}DK(\theta)N(\theta),$$

we obtain

$$\|S^0(\theta)\|_\rho \leq$$
$$\leq \|N(\theta)\|_\rho^2\, 2n\, \|DK(\theta)\|_\rho^2 \{2\|A\|_\rho \|B\|_{\mathcal{B}_r} +$$
$$+ \|B\|_{C^1, \mathcal{B}_r}^2 |H|_{C^1, \mathcal{B}_r} + \|B\|_{\mathcal{B}_r} \|N(\theta)\|_\rho\, 2n\, \|DK(\theta)\|_\rho^2\, 4n\, \|A\|_\rho\} \leq$$
$$\leq \|N(\theta)\|_\rho^2\, 2n\, \|DK(\theta)\|_\rho^2 \{4\|B\|_{C^1, \mathcal{B}_r}^2 |H|_{C^2, \mathcal{B}_r} +$$
$$+ \|B\|_{C^1, \mathcal{B}_r}^2 |H|_{C^1, \mathcal{B}_r} +$$
$$+ 16n^2 \|N(\theta)\|_\rho \|DK(\theta)\|_\rho^2 \|B\|_{C^1, \mathcal{B}_r}^2 |H|_{C^2, \mathcal{B}_r}\} \leq$$
$$\leq (\|N(\theta)\|_\rho + 1)^2\, 2n(\|DK(\theta)\|_\rho + 1)^2 \cdot$$
$$\cdot \{(16n^2 + 5)(\|N(\theta)\|_\rho + 1)(\|DK(\theta)\|_\rho + 1)^2 \|B\|_{C^1, \mathcal{B}_r}^2 |H|_{C^2, \mathcal{B}_r}\} \leq$$
$$\leq 2n(16n^2 + 5)(\|N(\theta)\|_\rho + 1)^3(\|DK(\theta)\|_\rho + 1)^4 \|B\|_{C^1, \mathcal{B}_r}^2 |H|_{C^2, \mathcal{B}_r}.$$

Therefore we get

$$\|S^0(\theta)\|_\rho \leq 2n(16n^2 + 5)(\|N(\theta)\|_\rho + 1)^3(\|DK(\theta)\|_\rho + 1)^4 \|B\|_{C^1, \mathcal{B}_r}^2 |H|_{C^2, \mathcal{B}_r}. \quad (4.68)$$

From (4.33): $M(\theta)^{-1} = V(\theta)^{-1}M(\theta)^T B(K(\theta))^{-1} + M_e(\theta)$, we get

$$\left\|M(\theta)^{-1}E(\theta)\right\|_{\rho-2\delta} = \left\|V(\theta)^{-1}M(\theta)^T B(K(\theta))^{-1}E(\theta) + M_e(\theta)E(\theta)\right\|_{\rho-2\delta} \le$$
$$\le \left\|V(\theta)^{-1}M(\theta)^T B(K(\theta))^{-1}E(\theta)\right\|_{\rho-2\delta} + \left\|M_e(\theta)E(\theta)\right\|_{\rho-2\delta} \le$$
$$\le \|V(\theta)^{-1}\|_{\rho-2\delta}\|M(\theta)^T B(K(\theta))^{-1}E(\theta)\|_{\rho-2\delta} + \left\|M_e(\theta)E(\theta)\right\|_{\rho-2\delta}$$
(4.69)

On the other hand,
$$M(\theta)^T B(K(\theta))^{-1}E(\theta) =$$
$$= \begin{pmatrix} DK(\theta)^T \\ N(\theta)^T DK(\theta)^T B(K(\theta))^T \end{pmatrix} B(K(\theta))^{-1}(De(\theta),\ W_1(\theta) - DK(\theta)S^0(\theta)) =$$
$$= \begin{pmatrix} DK(\theta)^T B(K(\theta))^{-1}De(\theta) & DK(\theta)^T B(K(\theta))^{-1}(W_1(\theta) - DK(\theta)S^0(\theta)) \\ N(\theta)^T DK(\theta)^T B(K(\theta))^T B(K(\theta))^{-1}De(\theta) & N(\theta)^T DK(\theta)^T B(K(\theta))^T B(K(\theta))^{-1}(W_1(\theta) - DK(\theta)S^0(\theta)) \end{pmatrix} =$$
$$= \begin{pmatrix} DK(\theta)^T B(K(\theta))^{-1}De(\theta) & DK(\theta)^T B(K(\theta))^{-1}(W_1(\theta) - DK(\theta)S^0(\theta)) \\ -N(\theta)DK(\theta)^T De(\theta) & -N(\theta)DK(\theta)^T(W_1(\theta) - DK(\theta)S^0(\theta)) \end{pmatrix}$$
(4.70)

holds. Here we get the estimate
$$\| DK(\theta)^T B(K(\theta))^{-1}(W_1(\theta) - DK(\theta)S^0(\theta))\|_{\rho-2\delta} \le$$
$$=\| DK(\theta)^T B(K(\theta))^{-1}W_1(\theta)\|_{\rho-2\delta} + \| DK(\theta)^T B(K(\theta))^{-1}DK(\theta)S^0(\theta)\|_{\rho-2\delta} =$$
$$=\| DK(\theta)^T B(K(\theta))^{-1}DB(K(\theta))e(\theta)DK(\theta)N(\theta)\|_{\rho-2\delta} +$$
$$+ \| De(\theta)^T DK(\theta)N(\theta)\|_{\rho-2\delta} + \|L(\theta)S^0(\theta)\|_{\rho-2\delta}.$$

And then from (4.63) and $S(\theta) = S^0(\theta)$, we have
$$N(\theta)DK(\theta)^T W_1(\theta) - S^0(\theta) =$$
$$= N(\theta)DK(\theta)^T B(K(\theta))DK(\theta)N(\theta)DK(\theta)^T B(K(\theta))^{-1} \cdot$$
$$\cdot DB(K(\theta))e(\theta)DK(\theta)N(\theta) +$$
$$- N(\theta)DK(\theta)^T B(K(\theta))DK(\theta)N(\theta)De(\theta)^T DK(\theta)N(\theta) -$$
$$- N(\theta)DK(\theta)^T B(K(\theta))DK(\theta)N(\theta)L(\theta)S(\theta).$$

Since
$$E(\theta) = (De(\theta),\ W_1(\theta) - DK(\theta)S^0(\theta))$$

and
$$M_e(\theta) = -[I_{2n} + V(\theta)^{-1}R(\theta)]^{-1}V(\theta)^{-1}R(\theta)V(\theta)^{-1}M(\theta)^T B(K(\theta))^{-1},$$

we have
$$M_e(\theta)E(\theta) = (M_e(\theta)De(\theta),\ M_e(\theta)[W_1(\theta) - DK(\theta)S^0(\theta)]).$$

In the last relation,
$$M_e(\theta)[W_1(\theta) - DK(\theta)S^0(\theta)] =$$
$$= -[I_{2n} + V(\theta)^{-1}R(\theta)]^{-1}V(\theta)^{-1}R(\theta)V(\theta)^{-1}M(\theta)^T B(K(\theta))^{-1}[W_1(\theta) - DK(\theta)S^0(\theta)]$$

holds. And we have

$$M(\theta)^T B(K(\theta))^{-1}[W_1(\theta) - DK(\theta)S^0(\theta)] =$$

$$= \begin{pmatrix} DK(\theta)^T \\ N(\theta)^T DK(\theta)^T B(K(\theta))^T \end{pmatrix} B(K(\theta))^{-1}[W_1(\theta) - DK(\theta)S^0(\theta)] =$$

$$= \begin{pmatrix} DK(\theta)^T B(K(\theta))^{-1}[W_1(\theta) - DK(\theta)S^0(\theta)] \\ N(\theta)^T DK(\theta)^T B(K(\theta))^T B(K(\theta))^{-1}[W_1(\theta) - DK(\theta)S^0(\theta)] \end{pmatrix} =$$

$$= \begin{pmatrix} DK(\theta)^T B(K(\theta))^{-1}[W_1(\theta) - DK(\theta)S^0(\theta)] \\ -N(\theta)DK(\theta)^T [W_1(\theta) - DK(\theta)S^0(\theta)] \end{pmatrix}.$$

Considering the above relations in (4.69), we obtain (4.55). □

**Theorem 4.1.** Let $K$ be an approximate solution of (2.16) such that $c_2 \gamma^{-1} \delta^{-(\sigma+1)} \|e\|_\rho \leq \dfrac{1}{2}$. Then the followings hold;

1) By the variable transformation $\Delta(\theta) = M(\theta)\xi(\theta)$, the linearized equation (4.21) is transformed to the following form.

$$\left[ \begin{pmatrix} 0 & S^0(\theta) \\ 0 & 0 \end{pmatrix} + C(\theta) \right] \xi(\theta) - \partial_\omega \xi(\theta) = p(\theta) + w(\theta), \qquad (4.71)$$

where

$$C(\theta) = M(\theta)^{-1} E(\theta)$$

$$p(\theta) = -\begin{pmatrix} N(\theta)^T DK(\theta)^T (I - B(K(\theta))DK(\theta)N(\theta)DK(\theta)B(K(\theta))^{-1})e(\theta) \\ DK(\theta)^T B(K(\theta))^{-1} e(\theta) \end{pmatrix}$$

$$w(\theta) = M_e(\theta)e(\theta).$$

Moreover,

$$\begin{aligned} \|C\|_{\rho-2\delta} &\leq c_4 \gamma^{-1} \delta^{-(\sigma+1)} \|e\|_\rho, \\ \|w\|_{\rho-\delta} &\leq c_4 \gamma^{-1} \delta^{-(\sigma+1)} \|e\|_\rho^2 \end{aligned} \qquad (4.72)$$

hold, where $c_4$ is the constant in lemma 4.6.

2) If we take $C(\theta) = 0$, $w(\theta) = 0$ in (4.71), then (4.71) is reduced to

$$\partial_\omega \xi_x(\theta) = S^0(\theta)\xi_y(\theta) - p_x(\theta) \qquad (4.73)\text{-}1$$

$$\partial_\omega \xi_y(\theta) = -p_y(\theta) . \qquad (4.73)\text{-}2$$

Moreover for any $0 < \delta < \dfrac{\rho}{2}$, there exists a solution $\xi = (\xi_x, \xi_y) \in \mathcal{P}(U_{\rho-\delta}, \mathbf{R}^{2n})$ of (4.73) such that there exist a constant $c_5$, depending only on

$$n, \sigma, \|N\|_\rho, \|DK\|_\rho, |<S>^{-1}|, \|B\|_{C^1, \mathcal{B}_r}, \|B^{-1}\|_{\mathcal{B}_r}, |H|_{C^2, \mathcal{B}_r},$$

satisfying $c_5 \geq c_4$ such that

$$\begin{aligned} \|\xi_y\|_{\rho-\delta} &\leq c_5 \gamma^{-1} \delta^{-\sigma} \|e\|_\rho, \\ \|\xi_x\|_{\rho-2\delta} &\leq c_5 \gamma^{-2} \delta^{-2\sigma} \|e\|_\rho. \end{aligned} \qquad (4.74)$$

**Proof.** If we take the variable transformation $\Delta(\theta) = M(\theta)\xi(\theta)$ for the linearized equation

(4.21): $A(\theta)\Delta(\theta) - \partial_\omega \Delta(\theta) = -e(\theta)$, then we obtain
$$[A(\theta)M(\theta) - \partial_\omega M(\theta)]\xi(\theta) - M(\theta)\partial_\omega \xi(\theta) = -e(\theta).$$

Therefore we have
$$M(\theta)^{-1}[A(\theta)M(\theta) - \partial_\omega M(\theta)]\xi(\theta) - \partial_\omega \xi(\theta) = -M(\theta)^{-1}e(\theta). \tag{4.75}$$

From (4.54), we get
$$M(\theta)^{-1}(A(\theta)M(\theta) - \partial_\omega M(\theta)) =$$
$$= M(\theta)^{-1}\left[M(\theta)\begin{pmatrix} 0 & S^0(\theta) \\ 0 & 0 \end{pmatrix} + E(\theta)\right] = \begin{pmatrix} 0 & S^0(\theta) \\ 0 & 0 \end{pmatrix} + M(\theta)^{-1}E(\theta). \tag{4.76}$$

Here we let $C(\theta) = M(\theta)^{-1}E(\theta)$. If we let
$$p(\theta) = -V(\theta)^{-1}M(\theta)^T B(K(\theta))^{-1}e(\theta), \quad w(\theta) = M_e(\theta)e(\theta), \tag{4.77}$$

then we have
$$-M(\theta)^{-1}e(\theta) = p(\theta) + w(\theta). \tag{4.78}$$

By considering (4.76), (4.77) in (4.75), we obtain (4.71). We obtain the estimate for $C(\theta)$ from lemma 4.6 and the estimate for $w(\theta)$ from lemma 4.3.

2): If we take $C(\theta) = 0$, $w(\theta) = 0$ in (4.71), then (4.71) is reduced to
$$\begin{pmatrix} 0 & S^0(\theta) \\ 0 & 0 \end{pmatrix}\xi(\theta) - \partial_\omega \xi(\theta) = p(\theta). \tag{4.79}$$

If we let in (4.79)
$$\xi(\theta) = \begin{pmatrix} \xi_x(\theta) \\ \xi_y(\theta) \end{pmatrix}, \quad p(\theta) = \begin{pmatrix} p_x(\theta) \\ p_y(\theta) \end{pmatrix},$$

then we have
$$\begin{pmatrix} 0 & S^0(\theta) \\ 0 & 0 \end{pmatrix}\begin{pmatrix} \xi_x(\theta) \\ \xi_y(\theta) \end{pmatrix} - \begin{pmatrix} \partial_\omega \xi_x(\theta) \\ \partial_\omega \xi_y(\theta) \end{pmatrix} = \begin{pmatrix} p_x(\theta) \\ p_y(\theta) \end{pmatrix}$$

and hence (4.79) is denoted as system of equations with $\xi_x$ and $\xi_y$ such as (4.73).

From (4.32), (4.37) we have
$$p(\theta) = -V(\theta)^{-1}M(\theta)^T B(K(\theta))^{-1}e(\theta) =$$
$$= -\begin{pmatrix} -N(\theta)^T DK(\theta)^T B(K(\theta))DK(\theta)N(\theta) & -I \\ I & 0 \end{pmatrix}\begin{pmatrix} DK(\theta)^T \\ N(\theta)^T DK(\theta)^T B(K(\theta))^T \end{pmatrix}B(K(\theta))^{-1}e(\theta) =$$
$$= -\begin{pmatrix} -N(\theta)^T DK(\theta)^T B(K(\theta))DK(\theta)N(\theta) & -I \\ I & 0 \end{pmatrix}\begin{pmatrix} DK(\theta)^T B(K(\theta))^{-1}e(\theta) \\ -N(\theta)^T DK(\theta)^T e(\theta) \end{pmatrix} =$$
$$= \begin{pmatrix} N(\theta)^T DK(\theta)^T [B(K(\theta))DK(\theta)N(\theta)DK(\theta)^T B(K(\theta))^{-1} - I_{2n}]e(\theta) \\ -DK(\theta)^T B(K(\theta))^{-1}e(\theta) \end{pmatrix}.$$

Therefore we obtain
$$p_x(\theta) = N(\theta)^T DK(\theta)^T [B(K(\theta))DK(\theta)N(\theta)DK(\theta)^T B(K(\theta))^{-1} - I_{2n}]e(\theta)$$
$$p_y(\theta) = -DK(\theta)^T B(K(\theta))^{-1}e(\theta).$$

Since $<DK(\theta)^T B(K(\theta))^{-1}e(\theta)>=0$ by theorem 3.2 and $\omega$ satisfies the Diophantine condition, using lemma 3.2, (4.73)-2: $\partial_\omega \xi_y(\theta) = -p_y(\theta)$ has unique analytic solution $\tilde{\xi}_y(\theta)$, $\theta \in U_\rho$ with zero average and

$$\|\tilde{\xi}_y\|_{\rho-\delta} \leq \mu \gamma^{-1} \delta^{-\sigma} \|p_y\|_\rho \tag{4.80}$$

holds. Let $\bar{\xi}_y \in \mathbf{R}^n$ be a function that $\xi_y(\theta) = \tilde{\xi}_y(\theta) + \bar{\xi}_y$ satisfies

$$<S^0(\theta)\xi_y - p_x(\theta))>=0.$$

We let

$$\bar{\xi}_y = <S^0(\theta)>^{-1}<p_x(\theta)> - <S^0(\theta)>^{-1}<S_0(\theta)\tilde{\xi}_y(\theta)>. \tag{4.81}$$

Since $<S^0(\theta)\xi_y - p_x(\theta)>=0$, from lemma 3.2

$$(4.73)\text{-}1: \partial_\omega \xi_x = S^0(\theta)\xi_y - p_x(\theta)$$

has unique analytic solution $\xi_x(\theta)$, $\theta \in U_\rho$ with zero average. From the denotation of $p_x(\theta)$, there exists a constant $c_5' \geq c_4$ depending only on

$n, \sigma, \|N\|_\rho, \|DK\|_\rho, \|B^{-1}\|_{B_r}, |B\|_{C^1, B_r}, H|_{C^2, B_r}$ such that

$$\|p_x\|_\rho \leq (2\pi)^{-n} c_5' \|e\|_\rho. \tag{4.82}$$

Now let us estimate $\xi_x$, $\xi_y$. First we estimate $|\bar{\xi}_y|$. From (4.82)

$$\left|<S^0(\theta)>^{-1}<p_x(\theta)>\right| \leq \left|<S^0(\theta)>^{-1}\right| \cdot \int_{\mathbf{T}^n} |p_x(\theta)| d\theta \leq c_5' \left|<S^0(\theta)>^{-1}\right| \cdot \|e\|_\rho. \tag{4.83}$$

From the definition (4.23) of $S^0(\theta)$, (4.80) and $p_y(\theta) = DK(\theta)^T B(K(\theta))^{-1} e(\theta)$, there exist a constant $c_5'' \geq c_4$ depending only on $n, \sigma, \|N\|_\rho, \|DK\|_\rho, \|B^{-1}\|_{B_r}, \|B\|_{C^1, B_r}, H|_{C^2, B_r}$ such that

$$|<S^0(\theta)>^{-1}<S^0(\theta)\tilde{\xi}_y(\theta)>| \leq |<S^0(\theta)>^{-1}| \int_{\mathbf{T}^n} |S^0(\theta)\tilde{\xi}_y(\theta)| d\theta \leq$$

$$\leq |<S^0(\theta)>^{-1}| c_5'' \mu \gamma^{-1} \delta^{-\sigma} \|e\|_\rho. \tag{4.84}$$

Then from (4.81),(4.83) and (4.84) we have

$$|\bar{\xi}_y| \leq c_5' \left|<S^0(\theta)>^{-1}\right| \cdot \|e\|_\rho + |<S^0(\theta)>^{-1}| c_5'' \mu \gamma^{-1} \delta^{-\sigma} \|e\|_\rho =$$
$$= |<S^0(\theta)>^{-1}| \cdot [c_5' + c_5'' \mu \gamma^{-1} \delta^{-\sigma}] \|e\|_\rho \leq \tag{4.85}$$
$$\leq |<S^0(\theta)>^{-1}| \cdot [c_5' + c_5'' \mu] \gamma^{-1} \delta^{-\sigma} \|e\|_\rho.$$

Therefore

$$\|\xi_y\|_{\rho-\delta} \leq \|\tilde{\xi}_y\|_{\rho-\delta} + |\bar{\xi}_y| \leq$$
$$\leq \|DK\|_\rho \|B^{-1}\|_{B_r} \gamma^{-1} \delta^{-\sigma} \|e\|_\rho + |<S^0(\theta)>^{-1}| \cdot [c_5' + c_5'' \mu] \gamma^{-1} \delta^{-\sigma} \|e\|_\rho \leq$$
$$\leq (\|DK\|_\rho + 1)(\|B^{-1}\|_{B_r} + 1)(|<S^0(\theta)>^{-1}| + 1)[c_5' + c_5'' \mu] \gamma^{-1} \delta^{-\sigma} \|e\|_\rho \leq$$
$$\leq c_5''' \gamma^{-1} \delta^{-\sigma} \|e\|_\rho,$$

where

$$c_5''' = (\|DK\|_\rho + 1)(\|B^{-1}\|_{\mathcal{B}_r} + 1)(|<S^0(\theta)>^{-1}| + 1)[c_5' + c_5''\mu]. \qquad (4.86)$$

Then $c_5''' \geq c_4$ holds. Now, we estimate $\|\xi_x\|_{\rho-2\delta}$. By using lemma 3.2,

$$\|\xi_x\|_{\rho-2\delta} \leq \mu\gamma^{-1}\delta^{-\sigma} \|S^0(\theta)\xi_y - p_x(\theta)\|_{\rho-\delta} \leq$$

$$\leq \mu\gamma^{-1}\delta^{-\sigma} \|S^0(\theta)\xi_y\|_{\rho-\delta} + \|p_x\|_{\rho-\delta} \leq$$

$$\leq \mu\gamma^{-1}\delta^{-\sigma} \|S^0(\theta)\|_{\rho-\delta}\|\xi_y\|_{\rho-\delta} + (2\pi)^{-n} c_5' \|e\|_{\rho-\delta} \leq$$

$$\leq \mu\gamma^{-1}\delta^{-\sigma} \|S^0(\theta)\|_{\rho-\delta} c_5'''\gamma^{-1}\delta^{-\sigma} \|e\|_{\rho-\delta} + (2\pi)^{-n} c_5' \|e\|_{\rho-\delta} \leq$$

$$\leq c_5\gamma^{-2}\delta^{-2\sigma} \|e\|_\rho,$$

where $c_5$ is a constant depending only on $n$, $\sigma$, $\|N\|_\rho$, $\|DK\|_\rho$, $|<S^0>^{-1}|$, $\|B^{-1}\|_{\mathcal{B}_r}$, $\|B\|_{C^1,\mathcal{B}_r}$, $H|_{C^2,\mathcal{B}_r}$, satisfying $c_5 \geq c_5'''$. □

**Lemma 4.7.** We assume that $c_4\gamma^{-1}\delta^{-(\sigma+1)} \|e\|_\rho \leq 1/2$ holds. Then

1) If $\xi(\theta)$, $\theta \in U_{\rho-\delta}$ is a solution of the reduced linearized equation (4.73) that the existence is guaranteed by theorem 4.1, then $\Delta(\theta) = M(\theta)\xi(\theta)$, $\theta \in U_{\rho-\delta}$ satisfies

$$A(\theta)\Delta(\theta) - \partial_\omega \Delta(\theta) = -e(\theta) + M(\theta)[C(\theta)\xi(\theta) - w(\theta)]. \qquad (4.87)$$

2) Moreover, there exist constants $c_6$, $c_7$ depending only on $n$, $\sigma$, $\|N\|_\rho$, $\|DK\|_\rho$, $|<S^0>^{-1}|$, $\|B\|_{C^1,\mathsf{B}_r}$, $\|B^{-1}\|_{\mathsf{B}_r}$, $|H|_{C^2,\mathsf{B}_r}$, satisfying $c_7 \geq c_6 \geq c_5$ such that

$$\|\Delta\|_{\rho-2\delta} \leq c_6\gamma^{-2}\delta^{-2\sigma} \|e\|_\rho \qquad (4.88)$$

and

$$\|DF(K)\Delta + e\|_{\rho-2\delta} \leq c_7\gamma^{-3}\delta^{-(3\sigma+1)} \|e\|_\rho^2. \qquad (4.89)$$

**Proof.** The condition $c_4\gamma^{-1}\delta^{-(\sigma+1)} \|e\|_\rho \leq 1/2$ guarantees the existence of $M(\theta)^{-1}$, thus the existence of the variable transformation $\Delta(\theta) = M(\theta)\xi(\theta)$. Let $\xi(\theta) = \begin{pmatrix} \xi_x(\theta) \\ \xi_y(\theta) \end{pmatrix}$, $\theta \in U_{\rho-\delta}$ be the solution of the reduced linearized equation (4.79) that its existence is guaranteed by theorem 4.1. Let us transform (4.79) into the form

$$\left[\begin{pmatrix} 0 & S^0(\theta) \\ 0 & 0 \end{pmatrix} + C(\theta)\right]\xi(\theta) - \partial_\omega \xi(\theta) = p(\theta) + w(\theta) + [C(\theta)\xi(\theta) - w(\theta)]. \qquad (4.90)$$

Since

$$\begin{pmatrix} 0 & S^0(\theta) \\ 0 & 0 \end{pmatrix} + C(\theta) = M(\theta)^{-1}[A(\theta)M(\theta) - \partial_\omega M(\theta)],$$

(4.90) is written by

$$M(\theta)^{-1}[A(\theta)M(\theta) - \partial_\omega M(\theta)]\xi(\theta) - \partial_\omega \xi(\theta) = p(\theta) + w(\theta) + [C(\theta)\xi(\theta) - w(\theta)].$$

Since from (4.37) we have $M(\theta)[p(\theta) + w(\theta)] = -e(\theta)$, the last relation is written by

$$A(\theta)M(\theta)\xi(\theta) - \partial_\omega M(\theta)\xi(\theta) - M(\theta)\partial_\omega \xi(\theta) =$$
$$= M(\theta)[p(\theta) + w(\theta)] + M(\theta)[C(\theta)\xi(\theta) - w(\theta)] =$$
$$= -e(\theta) + M(\theta)[C(\theta)\xi(\theta) - w(\theta)].$$

Therefore (4.87) holds.

Now we confirm (4.88). In fact since $\Delta(\theta) = M(\theta)\xi(\theta)$ and $M(\theta) = (DK(\theta) \quad B(K(\theta))DK(\theta)N(\theta))$, we have

$$\|\Delta\|_{\rho-2\delta} \leq (\|DK\|_\rho + 1)(\|N\|_\rho + 1)(\|B\|_{\mathcal{B}_r} + 1)\|\xi\|_{\rho-2\delta}. \tag{4.91}$$

From theorem 4.1, $\xi(\theta) = \begin{pmatrix} \xi_x(\theta) \\ \xi_y(\theta) \end{pmatrix}$ satisfies

$$\|\xi\|_{\rho-2\delta} \leq \max\{\|\xi_x\|_{\rho-2\delta},\ \|\xi_y\|_{\rho-\delta}\} \leq \max\{c_5\gamma^{-2}\delta^{-2\sigma}\|e\|_\rho,\ c_5\gamma^{-1}\delta^{-\sigma}\|e\|_\rho\} \leq c_5\gamma^{-2}\delta^{-2\sigma}\|e\|_\rho.$$

Therefore if we let $c_6 = (\|DK\|_\rho + 1)(\|N\|_\rho + 1)(\|B\|_{\mathcal{B}_r} + 1)c_5$, then

$$\|\Delta\|_{\rho-2\delta} \leq c_6\gamma^{-2}\delta^{-2\sigma}\|e\|_\rho.$$

On the other hand, since from (4.87) we have
$$DF(K)\Delta(\theta) + e(\theta) = A(\theta)\Delta(\theta) - \partial_\omega \Delta(\theta) + e(\theta) = M(\theta)[C(\theta)\xi(\theta) - w(\theta)],$$

we get

$$\|DF(K)\Delta + e\|_{\rho-2\delta} \leq \|M\|_{\rho-2\delta}[\|C\|_{\rho-2\delta} \cdot \|\xi\|_{\rho-2\delta} + \|w\|_{\rho-2\delta}].$$

Then from theorem 4.1 we have

$$\|C\|_{\rho-2\delta} \leq c_4\gamma^{-1}\delta^{-(\sigma+1)}\|e\|_\rho, \qquad \|\xi_y\|_{\rho-\delta} \leq c_5\gamma^{-1}\delta^{-\sigma}\|e\|_\rho,$$
$$\|w\|_{\rho-\delta} \leq c_4\gamma^{-1}\delta^{-(\sigma+1)}\|e\|_\rho^2, \qquad \|\xi_x\|_{\rho-2\delta} \leq c_5\gamma^{-2}\delta^{-2\sigma}\|e\|_\rho.$$

Therefore we have

$$\|DF(K)\Delta + e\| \leq$$
$$\leq \|M\|_\rho [c_4\gamma^{-1}\delta^{-(\sigma+1)}\|e\|_\rho \cdot c_5\gamma^{-2}\delta^{-2\sigma}\|e\| + c_5\gamma^{-1}\delta^{-(\sigma+1)}\|e\|_\rho^2] =$$
$$= \|M\|_\rho [c_4 c_5 \gamma^{-3}\delta^{-(3\sigma+1)}\|e\|_\rho^2 + c_5\gamma^{-1}\delta^{-(\sigma+1)}\|e\|_\rho^2] \leq$$
$$\leq (\|DK\|_\rho + 1)(\|N\|_\rho + 1)(\|B\|_{\mathcal{B}_r} + 1)c_4 c_5 \gamma^{-3}\delta^{-(3\sigma+1)}\|e\|_\rho^2.$$

Hence if we let $c_7 = (\|DK\|_\rho + 1)(\|N\|_\rho + 1)(\|B\|_{\mathcal{B}_r} + 1)c_4 c_6$, then

$$\|DF(K)\Delta + e\|_{\rho-2\delta} \leq c_7\gamma^{-3}\delta^{-(3\sigma+1)}\|e\|_\rho^2. \quad \square$$

## 5. Confirmation of non-degenerate property of approximate solutions

Given $\rho_0 > 0$, we choose $0 < \delta_0 < \min\{1, \rho_0/12\}$. And put

$$\delta_m = \delta_0 2^{-m} \ (m \geq 1), \quad \rho_m = \rho_{m-1} - 3\delta_{m-1} \ (m \geq 1).$$

Then

$$\rho_m = \rho_{m-1} - 3\delta_{m-1} = \rho_0 - 6\delta_0(1 - 2^{-m}) > 0.$$

We put

$$\rho_\infty := \lim_{m\to\infty} \rho_m = \rho_0 - 6\delta_0 > 0.$$

Assume that $K_0 \in \mathcal{ND}(\rho_0)$ is an approximate solution of (2.16). We put $e_0 = F(K_0)$. Let's suppose that $\Delta K_0(\theta)$, $\theta \in U_{\rho_0-\delta_0}$ is the approximate solution of

$$DF(K_0)\Delta K_0 = -e_0 \tag{5.1}$$

given by theorem 4.1 and put $K_1(\theta) = K_0(\theta) + \Delta K_0(\theta)$, $\theta \in U_{\rho_0-\delta_0}$. Let's contract the domain of $K_1$: $K_1(\theta)$, $\theta \in U_{\rho_1} = U_{\rho_0-3\delta_0}$. Let's assume that $K_{m-1}(\theta)$, $\theta \in U_{\rho_{m-1}}$ is already defined for $m \geq 1$. Then we define

$$e_{m-1} = F(K_{m-1}). \tag{5.2}$$

Let's suppose that $\Delta K_{m-1}(\theta)$, $\theta \in U_{\rho_{m-1}-\delta_{m-1}}$ is the approximate solution of

$$DF(K_{m-1})\Delta K_{m-1} = -e_{m-1} \tag{5.3}$$

based on $K_{m-1}(\theta)$ given by theorem 4.1. We define

$$K_m(\theta) = K_{m-1}(\theta) + \Delta K_{m-1}(\theta), \quad \theta \in U_{\rho_m}. \tag{5.4}$$

For $m \geq 0$ we define

$$N_m(\theta) = (DK_m(\theta)^T DK_m(\theta))^{-1} \tag{5.5}$$

in the case of that $DK_m(\theta)^T DK_m(\theta)$ is invertible. And we define

$$A_m(\theta) = \Phi(K_m(\theta)) + B(K_m(\theta))\Psi(K_m(\theta)) \tag{5.6}$$

and

$$\begin{aligned}S_m^0(\theta) = &N_m(\theta)DK_m(\theta)^T \{A_m(\theta)B(K_m(\theta)) - B(K_m(\theta))A_m(\theta) - \\ &- DB(K_m(\theta))B(K_m(\theta))\nabla H(K_m(\theta)) + \\ &+ B(K_m(\theta))DK_m(\theta)N_m(\theta)DK_m(\theta)^T[A_m(\theta)+A_m(\theta)^T]\}DK_m(\theta)N_m(\theta).\end{aligned} \tag{5.7}$$

Now for $1 \leq j \leq 7$ we denote constant $c_j$, that we replace $K$ in $c_j$ with $K_{m-1}$, by $c_{m-1}^{(j)}$. For example constant

$$c_4(n, \ \sigma, \ \|N\|_\rho, \ \|DK\|_\rho, \ \|B^{-1}\|_{\mathcal{B}_r}, \ \|B\|_{C^1,\mathcal{B}_r}, \ H|_{C^2,\mathcal{B}_r})$$

in lemma 4.6 is replaced with

$$c_{m-1}^{(4)} = c_4(n, \ \sigma, \ \|N_{m-1}\|_{\rho_{m-1}}, \ \|DK_{m-1}\|_{\rho_{m-1}}, \ \|B^{-1}\|_{\mathcal{B}_r}, \ \|B\|_{C^1,\mathcal{B}_r}, \ H|_{C^2,\mathcal{B}_r}).$$

**Lemma 5.1.** Assume that for $m \geq 1$ $K_{m-1} \in \mathcal{ND}(\rho_{m-1})$ is approximate solution of (2.16) and

$$r_{m-1} := \|K_{m-1} - K_0\|_{\rho_{m-1}} < r. \tag{5.8}$$

If $e_{m-1}$ defined by (5.2) satisfies condition $c_{m-1}^{(4)}\gamma^{-1}\delta_{m-1}^{-(\sigma+1)}\|e_{m-1}\|_{\rho_{m-1}} \leq 1/2$ (then lemma 4.7 stands), then there exists function $\Delta K_{m-1} \in \mathcal{P}_{\rho_{m-1}-\delta_{m-1}}$ such that

$$\|\Delta K_{m-1}\|_{\rho_{m-1}-2\delta_{m-1}} \leq c_{m-1}^{(8)}\gamma^{-2}\delta_{m-1}^{-2\sigma}\|e_{m-1}\|_{\rho_{m-1}} \tag{5.9}$$

$$\|D\Delta K_{m-1}\|_{\rho_{m-1}-3\delta_{m-1}} \leq c_{m-1}^{(8)}\gamma^{-2}\delta_{m-1}^{-(2\sigma+1)}\|e_{m-1}\|_{\rho_{m-1}}. \tag{5.10}$$

Here $c_{m-1}^{(8)}$ is constant depending on

$$n, \sigma, \|N_{m-1}\|_{\rho_{m-1}}, \|DK_{m-1}\|_{\rho_{m-1}}, |<S_{m-1}^{0}>^{-1}|, \|B\|_{C^1,\mathcal{B}_r}, \|B^{-1}\|_{\mathcal{B}_r}, |H|_{C^2,\mathcal{B}_r}$$

that satisfies $c_{m-1}^{(8)} \geq c_{m-1}^{(7)}$. Moreover if

$$r_{m-1} + c_{m-1}^{(8)}\gamma^{-2}\delta_{m-1}^{-2\sigma}\|e_{m-1}\|_{\rho_{m-1}} < r \tag{5.11}$$

holds, then $K_m(\theta) = K_{m-1}(\theta) + \Delta K_{m-1}(\theta), \; (\theta \in U_{\rho_m})$ satisfies

$$K_m(U_{\rho_m}) \subset \mathcal{B}_r = \{z \in \mathbf{C}^{2n} \; ; \; \inf_{|\mathrm{Im}\theta| \leq \rho_0}|z - K_0(\theta)| < r\} \tag{5.12}$$

and there exists constant $c_{m-1}^{(9)}$ depending on

$$n, \sigma, \|N_{m-1}\|_{\rho_{m-1}}, \|DK_{m-1}\|_{\rho_{m-1}}, |<S_{m-1}^{0}>^{-1}|, \|B\|_{C^2,\mathcal{B}_r}, \|B^{-1}\|_{\mathcal{B}_r}, |H|_{C^3,\mathcal{B}_r}$$

satisfying $c_{m-1}^{(9)} \geq c_{m-1}^{(8)}$ such that

$$\|e_m\|_{\rho_m} \leq c_{m-1}^{(9)}\gamma^{-4}\delta_{m-1}^{-4\sigma}\|e_{m-1}\|_{\rho_{m-1}}^2. \tag{5.13}$$

**Proof.** Let's suppose that $\Delta K_{m-1}(\theta) = M_{m-1}(\theta)\xi_{m-1}(\theta), \; \theta \in \mathcal{P}_{\rho_{m-1}-\delta_{m-1}}$ is the approximate solution of $DF(K_{m-1})\Delta K_{m-1} = -e_{m-1}$ given by theorem 4.1. First we prove (5.9), (5.10). From the estimate $\|\Delta\|_{\rho-2\delta} \leq c_6\gamma^{-2}\delta^{-2\sigma}\|e\|_\rho$ in lemma 4.7 we obtain

$$\|\Delta K_{m-1}\|_{\rho_{m-1}-2\delta_{m-1}} \leq c_{m-1}^{(6)}\gamma^{-2}\delta_{m-1}^{-2\sigma}\|e_{m-1}\|_{\rho_{m-1}}.$$

Now we put $S(\theta, \delta) = \{\varphi \in \mathbf{C}^n ; |\varphi - \theta| \leq \delta\}$ for $\theta \in \mathbf{C}^n$ and positive number $\delta$. Let's fix $\theta \in U_{\rho_{m-1}-3\delta_{m-1}}$ arbitrarily. Then $S(\theta, \delta_{m-1}) \subset U_{\rho_{m-1}-2\delta_{m-1}}$. Therefore from Cauchy' estimate we have

$$|D(\Delta K_{m-1})(\theta)| \leq n\delta_{m-1}^{-1}\max\{|\Delta K_{m-1}(\theta)| \; ; \; |\theta| \leq \delta_{m-1}\} \leq$$
$$\leq n\delta_{m-1}^{-1}\max\{|\Delta K_{m-1}(\theta)| \; ; \; \theta \in U_{\rho_{m-1}-2\delta_{m-1}}\} =$$
$$\leq n\delta_{m-1}^{-1}\|\Delta K_{m-1}\|_{\rho_{m-1}-2\delta_{m-1}} \leq nc_{m-1}^{(6)}\gamma^{-2}\delta_{m-1}^{-2\sigma}\|e_{m-1}\|_{\rho_{m-1}}.$$

Then for $c_{m-1}^{(8)} = nc_{m-1}^{(7)}$

$$\| D\Delta K_{m-1} \|_{\rho_{m-1}-3\delta_{m-1}} \leq c_{m-1}^{(8)}\gamma^{-2}\delta_{m-1}^{-(2\sigma+1)} \| e_{m-1} \|_{\rho_{m-1}}$$

holds. Because $K_m(\theta) = K_{m-1}(\theta) + \Delta K_{m-1}(\theta),\ \theta \in U_{\rho_m}$, we have

$$\| K_m - K_0 \|_{\rho_{m-1}-2\delta_{m-1}} \leq \| K_m - K_{m-1} \|_{\rho_{m-1}-2\delta_{m-1}} + \| K_{m-1} - K_0 \|_{\rho_{m-1}-2\delta_{m-1}} \leq$$
$$\leq c_{m-1}^{(6)}\gamma^{-2}\delta_{m-1}^{-2\sigma} \| e_{m-1} \|_{\rho_{m-1}} + r_{m-1} < r$$

from (5.11). Hence for any $\theta$ satisfying $|\mathrm{Im}\theta| \leq \rho_{m-1} - 2\delta_{m-1}$,

$$\inf_{|\mathrm{Im}\theta'| \leq \rho_0} | K_m(\theta) - K_0(\theta') | \leq | K_m(\theta) - K_0(\theta) | \leq \| K_m - K_0 \|_{\rho_{m-1}-2\delta_{m-1}} < r$$

thus $K_m(U_{\rho_m}) \subset \mathcal{B}_r = \{z \in \mathbf{C}^{2n}\ ;\ \inf_{|\mathrm{Im}\theta'| \leq \rho_0} |z - K_0(\theta')| < r\}$ holds.

If we define

$$R(K, v) = F(v) - F(K) - DF(K)(v - K), \tag{5.14}$$

then

$$F(K_m) = F(K_{m-1}) + DF(K_{m-1})(K_m - K_{m-1}) + R(K_{m-1},\ K_m).$$

Therefore

$$e_m(\theta) = e_{m-1}(\theta) + DF(K_{m-1})\Delta K_{m-1}(\theta) + R(K_{m-1},\ K_m)(\theta).$$

Accordingly

$$\|e_m\|_{\rho_m} \leq \|e_{m-1} + DF(K_{m-1})\Delta K_{m-1}\|_{\rho_m} + \|R(K_{m-1},\ K_m)\|_{\rho_m}.$$

Because of $\rho_m = \rho_{m-1} - 3\delta_{m-1}$, we have

$$\|e_{m-1} + DF(K_{m-1})\Delta K_{m-1}\|_{\rho_m} = \|e_{m-1} + DF(K_{m-1})\Delta K_{m-1}\|_{\rho_{m-1}-3\delta_{m-1}} \leq$$
$$\leq \|e_{m-1} + DF(K_{m-1})\Delta K_{m-1}\|_{\rho_{m-1}-2\delta_{m-1}} \leq c_{m-1}^{(7)}\gamma^{-3}\delta_{m-1}^{-(3\sigma+1)} \| e_{m-1} \|_{\rho_{m-1}}^2$$

from (4.89).

Let's estimate $\| D^2F(K_{m-1} + t\Delta K_{m-1}) \|$. If we put

$$f(K) = (B \circ K)(\nabla H \circ K),\ g(K) = \partial_\omega K,$$

then

$$F(K) \coloneqq (B \circ K)(\nabla H \circ K) - \partial_\omega K = f(K) - g(K)$$

and $D^2 g(K) = 0$ because $g(K)$ is continuous linear map for $K$. Therefore we have

$D^2F(K) = D^2 f(K)$. From composition theorem ([Irwin 1980])

$$[Dg_*(f)h](x) = [(Dg \circ f)\cdot h](x) = Dg(f(x))h(x) = \frac{d}{dt}g(f(x) + th(x))\bigg|_{t=0}.$$

And therefore

$$[D^2 f(K)(\xi, \eta)](\theta) =$$
$$= D^2 B(K(\theta))\xi(\theta)\eta(\theta)\nabla H(K(\theta)) + DB(K(\theta))\xi(\theta)D\nabla H(K(\theta))\eta(\theta) +$$
$$+ DB(K(\theta))\eta(\theta)D\nabla H(K(\theta))\xi(\theta) + B(K(\theta))D^2\nabla H(K(\theta))\xi(\theta)\eta(\theta).$$

Accordingly

$$\|D^2 f(K)\| := \sup\{\|D^2 f(K)(\xi, \eta)\|_{\rho_0}; \|\xi\|_{\rho_0} = 1, \|\eta\|_{\rho_0} = 1\} =$$
$$= \sup\{\sup_{\theta \in U_{\rho_0}} \|[D^2 f(K)(\xi, \eta)](\theta)\|; \|\xi\|_{\rho_0} = 1, \|\eta\|_{\rho_0} = 1\} \leq$$
$$\leq 4\|B\|_{C^2, \mathcal{B}_r} |H|_{C^3, \mathcal{B}_r}$$

holds, where $\|D^2 f(K)\|$ denotes the norm of continuous bilinear map ([Dieudonné 1960]). Therefore we have for any $t \in [0, 1]$

$$\|D^2 F(K_{m-1} + t\Delta K_{m-1})\| \leq 4\|B\|_{C^2, \mathcal{B}_r} \|H\|_{C^3, \mathcal{B}_r}$$

holds.

Next let's estimate $\|R(K_{m-1}, K_m)\|_{\rho_m}$. By the integral-typed Taylor's formula ([Lang 2002])

$$F(K_m) = F(K_{m-1}) + DF(K_{m-1})\Delta K_{m-1} + \int_0^1 D^2 F(K_{m-1} + t\Delta K_{m-1})(1-t)(\Delta K_{m-1})^2 dt,$$

we have

$$\|R(K_{m-1}, K_m)\|_{\rho_m} = \left\|\int_0^1 D^2 F(K_{m-1} + t\Delta K_{m-1})(1-t)(\Delta K_{m-1})^2 dt\right\|_{\rho_m} =$$
$$\leq 4\|B\|_{C^2, \mathcal{B}_r} |H|_{C^3, \mathcal{B}_r} \|\Delta K_{m-1}\|^2_{\rho_m} \leq$$
$$\leq 4\|B\|_{C^2, \mathcal{B}_r} |H|_{C^3, \mathcal{B}_r} \|\Delta K_{m-1}\|^2_{\rho_{m-1} - 2\delta_{m-1}} \leq$$
$$\leq 4\|B\|_{C^2, \mathcal{B}_r} |H|_{C^3, \mathcal{B}_r} c_{m-1}^{(6)\,2} \gamma^{-4} \delta_{m-1}^{-4\sigma} \|e_{m-1}\|^2_{\rho_{m-1}}.$$

Accordingly if we put $c_{m-1}^{(9)} = (c_{m-1}^{(8)} + 4\|B\|_{C^2, \mathcal{B}_r} \|H\|_{C^3, \mathcal{B}_r} c_{m-1}^{(6)\,2})$, then

$c_{m-1}^{(9)} \geq c_{m-1}^{(8)}$ and

$$\|e_m\|_{\rho_m} \leq \|F(K_{m-1}) + DF(K_{m-1})\Delta K_{m-1}\|_{\rho_m} + \|R(K_{m-1}, K_m)\|_{\rho_m} \leq$$
$$\leq c_{m-1}^{(8)} \gamma^{-3} \delta_{m-1}^{-(3\sigma+1)} \|e_{m-1}\|^2_{\rho_{m-1}} + 4\|B\|_{C^2, \mathcal{B}_r} \|H\|_{C^3, \mathcal{B}_r} c_{m-1}^{(6)\,2} \gamma^{-4} \delta_{m-1}^{-4\sigma} \|e_{m-1}\|^2_{\rho_{m-1}} \leq$$
$$\leq c_{m-1}^{(9)} \gamma^{-4} \delta_{m-1}^{-4\sigma} \|e_{m-1}\|^2_{\rho_{m-1}}$$

holds. □

**Lemma 5.2.** Let's put

$$c_{m-1}^{(10)} = 4(\|N_{m-1}\|_{\rho_{m-1}} + 1)^2 [2\|DK_{m-1}\|_{\rho_{m-1}} c_{m-1}^{(9)} + c_{m-1}^{(9)}].$$

Then $c_{m-1}^{(10)} \geq c_{m-1}^{(9)}$ and $c_{m-1}^{(10)}$ depend on

$$n, \sigma, \|N_{m-1}\|_{\rho_{m-1}}, \|DK_{m-1}\|_{\rho_{m-1}}, |<S_{m-1}^0>^{-1}|, \|B\|_{C^2, \mathcal{B}_r}, \|B^{-1}\|_{\mathcal{B}_r}, \|H\|_{C^3, \mathcal{B}_r}.$$

Assume that

a) $K_{m-1} \in \mathcal{ND}(\rho_{m-1})$ is approximate solution of (2.16) and $r_{m-1} := \|K_{m-1} - K_0\|_{\rho_{m-1}} < r$ holds for $m \geq 1$.

b) $e_{m-1}$ satisfies $c_{m-1}^{(10)} \gamma^{-2} \delta_{m-1}^{-(2\sigma+1)} \|e_{m-1}\|_{\rho_{m-1}} \leq 1/2$.

c) (5.11): $r_{m-1} + c_{m-1}^{(8)} \gamma^{-2} \delta_{m-1}^{-2\sigma} \|e_{m-1}\|_{\rho_{m-1}} < r$ holds.

Then we obtain

1) $DK_m(\theta)^T DK_m(\theta)$ is invertible and its inverse map $N_m$ satisfies

$$\|N_m\|_{\rho_m} \leq \|N_{m-1}\|_{\rho_{m-1}} + c_{m-1}^{(10)} \gamma^{-2} \delta_{\rho_{m-1}}^{-2(\sigma+1)} \|e_{m-1}\|_{\rho_{m-1}}. \tag{5.15}$$

2) There exists constant $c_{m-1}^{(11)}$ depending on

$$n, \sigma, \|N_{m-1}\|_{\rho_{m-1}}, \|DK_{m-1}\|_{\rho_{m-1}}, |<S_{m-1}^0>^{-1}|, \|B\|_{C^2, \mathcal{B}_r}, \|B^{-1}\|_{\mathcal{B}_r}, \|H\|_{C^3, \mathcal{B}_r}$$

that satisfies $c_{m-1}^{(11)} \geq c_{m-1}^{(10)}$ such that

$$\|\Delta S_{m-1}^0\|_{\rho_m} \leq c_{m-1}^{(11)} \gamma^{-2} \delta_{m-1}^{-2(\sigma+1)} \|e_{m-1}\|_{\rho_{m-1}}.$$

3) Let $c_{m-1}^{(12)} = 2(|<S_{m-1}^0>^{-1}|+1)^2 c_{m-1}^{(11)}$. If $e_{m-1}$ satisfies $c_{m-1}^{(12)} \gamma^{-2} \delta_{m-1}^{-(2\sigma+1)} \|e_{m-1}\|_{\rho_{m-1}} \leq 1/2$, then $<S_m^0>$ is invertible. Moreover

$$|<S_m^0>^{-1}| < |<S_{m-1}^0>^{-1}| + c_{m-1}^{(12)} \gamma^{-2} \delta_{\rho_{m-1}}^{-2(\sigma+1)} \|e_{m-1}\|_{\rho_{m-1}} \tag{5.16}$$

holds.

**Proof.** First, estimate norm $\|Z_{m-1}\|_{\rho_m}$ of

$$Z_{m-1}(\theta) = DK_m^T DK_m - DK_{m-1}^T DK_{m-1}, \quad \theta \in U_{\rho_m} \tag{5.17}$$

for $m \geq 1$. Since

$$Z_{m-1} = DK_m^T DK_m - DK_{m-1}^T DK_{m-1} =$$
$$= [DK_{m-1} + D\Delta K_{m-1}]^T [DK_{m-1} + D\Delta K_{m-1}] - DK_{m-1}^T DK_{m-1} =$$
$$= D\Delta K_{m-1}^T D\Delta K_{m-1} + DK_{m-1}^T D\Delta K_{m-1} + D\Delta K_{m-1}^T DK_{m-1},$$

if we consider the assumption b), then from lemma 5.1 we obtain

$$\|Z_{m-1}\|_{\rho_m} \le$$
$$\le 2\|DK_{m-1}\|_{\rho_{m-1}-2\delta_{m-1}} \|D\Delta K_{m-1}\|_{\rho_{m-1}-3\delta_{m-1}} + \|D\Delta K_{m-1}\|_{\rho_{m-1}-3\delta_{m-1}} \|D\Delta K_{m-1}\|_{\rho_{m-1}-3\delta_{m-1}} \le$$
$$\le 2\|DK_{m-1}\|_{\rho_{m-1}-2\delta_{m-1}} c_{m-1}^{(9)} \gamma^{-2} \delta_{m-1}^{-2(\sigma+1)} \|e_{m-1}\|_{\rho_{m-1}} + \left[ c_{m-1}^{(9)} \gamma^{-2} \delta_{m-1}^{-(2\sigma+1)} \|e_{m-1}\|_{\rho_{m-1}} \right]^2 \le$$
$$\le 2\|DK_{m-1}\|_{\rho_{m-1}} c_{m-1}^{(9)} \gamma^{-2} \delta_{m-1}^{-2(\sigma+1)} \|e_{m-1}\|_{\rho_{m-1}} + \frac{1}{2} c_{m-1}^{(9)} \gamma^{-2} \delta_{m-1}^{-(2\sigma+1)} \|e_{m-1}\|_{\rho_{m-1}} =$$
$$= [2\|DK_{m-1}\|_{\rho_{m-1}} c_{m-1}^{(9)} + \frac{1}{2} c_{m-1}^{(9)}] \gamma^{-2} \delta_{m-1}^{-2(\sigma+1)} \|e_{m-1}\|_{\rho_{m-1}}.$$

Then since

$$\|N_{m-1}Z_{m-1}\|_{\rho_m} \le \|N_{m-1}\|_{\rho_{m-1}} \|Z_{m-1}\|_{\rho_m} \le$$
$$\le \|N_{m-1}\|_{\rho_{m-1}} [2\|DK_{m-1}\|_{\rho_{m-1}} c_{m-1}^{(9)} + \frac{1}{2} c_{m-1}^{(9)}] \gamma^{-2} \delta_{m-1}^{-2(\sigma+1)} \|e_{m-1}\|_{\rho_{m-1}} \le$$
$$\le (\|N_{m-1}\|_{\rho_{m-1}} + 1)[2\|DK_{m-1}\|_{\rho_{m-1}} c_{m-1}^{(9)} + c_{m-1}^{(9)}] \gamma^{-2} \delta_{m-1}^{-2(\sigma+1)} \|e_{m-1}\|_{\rho_{m-1}} \le$$
$$\le c_{m-1}^{(10)} \gamma^{-2} \delta_{m-1}^{-2(\sigma+1)} \|e_{m-1}\|_{\rho_{m-1}} \le \frac{1}{2},$$

(5.18)

from Neumann's series theorem matrix $(I_n + N_{m-1}Z_{m-1})$ is invertible and

$$\|(I_n + N_{m-1}Z_{m-1})^{-1}\|_{\rho_m} \le (1 - \|N_{m-1}Z_{m-1}\|_{\rho_m})^{-1} \le (1 - \|N_{m-1}\|_{\rho_m} \cdot \|Z_{m-1}\|_{\rho_m})^{-1} \le 2.$$

Hence $DK_m^T DK_m = DK_{m-1}^T DK_{m-1} + Z_{m-1} = DK_{m-1}^T DK_{m-1}(I_n + N_{m-1}Z_{m-1})$ is invertible and

$$N_m = (I_n + N_{m-1}Z_{m-1})^{-1} N_{m-1},$$

$$\|N_m\|_{\rho_m} \le \|(I_n + N_{m-1}Z_{m-1})^{-1}\|_{\rho_m} \cdot \|N_{m-1}\|_{\rho_m} \le (1 - \|N_{m-1}\|_{\rho_m} \cdot \|Z_{m-1}\|_{\rho_m})^{-1} \|N_{m-1}\|_{\rho_m}.$$

Since

$$(1 - \|N_{m-1}\|_{\rho_m} \cdot \|Z_{m-1}\|_{\rho_m})^{-1} = 1 + (1 - \|N_{m-1}\|_{\rho_m} \cdot \|Z_{m-1}\|_{\rho_m})^{-1} \|N_{m-1}\|_{\rho_m} \cdot \|Z_{m-1}\|_{\rho_m} \le$$
$$\le 1 + 2\|N_{m-1}\|_{\rho_m} \cdot \|Z_{m-1}\|_{\rho_m},$$

we have

$$\|N_m\|_{\rho_m} \le (1 - \|N_{m-1}\|_{\rho_m} \cdot \|Z_{m-1}\|_{\rho_m})^{-1} \|N_{m-1}\|_{\rho_m} \le$$
$$\le \|N_{m-1}\|_{\rho_m} + 2\|N_{m-1}\|_{\rho_m}^2 \|Z_{m-1}\|_{\rho_m} \le$$
$$\le \|N_{m-1}\|_{\rho_m} + 2\|N_{m-1}\|_{\rho_m}^2 [2\|DK_{m-1}\|_{\rho_{m-1}} c_{m-1}^{(9)} + \frac{1}{2} c_{m-1}^{(9)}] \gamma^{-2} \delta_{m-1}^{-2(\sigma+1)} \|e_{m-1}\|_{\rho_{m-1}} \le$$
$$\le 2(\|N_{m-1}\|_{\rho_{m-1}} + 1)^2 [2\|DK_{m-1}\|_{\rho_{m-1}} c_{m-1}^{(9)} + c_{m-1}^{(9)}] \gamma^{-2} \delta_{m-1}^{-2(\sigma+1)} \|e_{m-1}\|_{\rho_{m-1}} \le$$
$$\le c_{m-1}^{(10)} \gamma^{-2} \delta_{m-1}^{-2(\sigma+1)} \|e_{m-1}\|_{\rho_{m-1}}.$$

Accordingly we obtain 1).

On the other hand, for any $t \in [0,\ 1]$

$$(1 - t \| N_{m-1} \|_{\rho_m} \cdot \| Z_{m-1} \|_{\rho_m})^{-1} \leq 2. \qquad (5.19)$$

Let's estimate the norm $\| \Delta S_{m-1}^0 \|_{\rho_m}$ of

$$\Delta S_{m-1}^0(\theta) = S_{m-1}^0(\theta) - S_{m-1}^0(\theta),\ \theta \in U_\rho$$

for $m \geq 1$. Let's introduce simplicity symbols:

$$B_{m-1}(\theta) = B(K_{m-1}(\theta)), \quad \Delta B_{m-1}(\theta) = B(K_m(\theta)) - B(K_{m-1}(\theta)),$$

$$DB_{m-1}(\theta) = DB(K_{m-1}(\theta)), \quad \Delta DB_{m-1}(\theta) = DB(K_m(\theta)) - DB(K_{m-1}(\theta)),$$

$$\nabla H_{m-1}(\theta) = \nabla H(K_{m-1}(\theta)),\ \Delta \nabla H_{m-1}(\theta) = \nabla H(K_m(\theta)) - \nabla H(K_{m-1}(\theta)).$$

Then

$$\begin{aligned}
\Delta S_{m-1}^0 &= S_m^0 - S_{m-1}^0 = \\
&= N_m DK_m^T \{ A_m B_m - B_m A_m - DB_m B_m \nabla H_m + \\
&\quad + B_m DK_m N_m DK_m [A_m + A_m^T] \} DK_m N_m - \\
&\quad - N_{m-1} DK_{m-1}^T \{ A_{m-1} B_{m-1} - B_{m-1} A_{m-1} - DB_{m-1} B_{m-1} \nabla H_{m-1} + \\
&\quad + B_{m-1} DK_{m-1} N_{m-1} DK_{m-1}^T [A_{m-1} + A_{m-1}^T] \} DK_{m-1} N_{m-1} = \\
&= (N_{m-1} + \Delta N_{m-1})(DK_{m-1}^T + D\Delta K_{m-1}^T) \cdot \\
&\quad \cdot \{ (A_{m-1} + \Delta A_{m-1})(B_{m-1} + \Delta B_{m-1}) - \\
&\quad - (B_{m-1} + \Delta B_{m-1})(A_{m-1} + \Delta A_{m-1}) - \\
&\quad - (DB_{m-1} + D\Delta B_{m-1})(B_{m-1} + \Delta B_{m-1})(\nabla H_{m-1} + \Delta \nabla H_{m-1}) + \\
&\quad + (B_{m-1} + \Delta B_{m-1})(DK_{m-1} + D\Delta K_{m-1})(N_{m-1} + \Delta N_{m-1})(DK_{m-1}^T + D\Delta K_{m-1}^T) \cdot \\
&\quad \cdot [(A_{m-1} + \Delta A_{m-1}) + (A_{m-1}^T + \Delta A_{m-1}^T)] \} (DK_{m-1} + D\Delta K_{m-1})(N_{m-1} + \Delta N_{m-1}) - \\
&\quad - N_{m-1} DK_{m-1}^T \{ A_{m-1} B_{m-1} - B_{m-1} A_{m-1} - DB_{m-1} B_{m-1} \nabla H_{m-1} + \\
&\quad + B_{m-1} DK_{m-1} N_{m-1} DK_{m-1}^T [A_{m-1} + A_{m-1}^T] \} DK_{m-1} N_{m-1}.
\end{aligned}$$

holds. Therefore $\| \Delta S_{m-1}^0 \|_{\rho_m}$ is bounded to upper direction by a polynomial with positive coefficients which variables are

$$\| N_{m-1} \|_{\rho_m},\ \| \Delta N_{m-1} \|_{\rho_m},$$
$$\| DK_{m-1} \|_{\rho_m},\ \| D\Delta K_{m-1} \|_{\rho_m},$$
$$\| A_{m-1} \|_{\rho_m},\ \| \Delta A_{m-1} \|_{\rho_m},$$
$$\| B_{m-1} \|_{\rho_m},\ \| \Delta B_{m-1} \|_{\rho_m},\ \| DB_{m-1} \|_{\rho_m},\ \| D\Delta B_{m-1} \|_{\rho_m},$$
$$\| \nabla H_{m-1} \|_{\rho_m},\ \| \nabla \Delta H_{m-1} \|_{\rho_m}.$$

Hence in order to estimate $\| \Delta S_{m-1}^0 \|_{\rho_m}$ by a polynomial with positive coefficients which variables are

$$\| N_{m-1} \|_{\rho_{m-1}}, \| DK_{m-1} \|_{\rho_{m-1}}, \| A_{m-1} \|_{\rho_{m-1}}, \| B \|_{C^k, \mathcal{B}_r}, \| H \|_{C^k, \mathcal{B}_r}, \tag{5.20}$$

we have to estimate quantities

$$\| \Delta N_{m-1} \|_{\rho_m}, \ \| D\Delta K_{m-1} \|_{\rho_m}, \ \| \Delta A_{m-1} \|_{\rho_m}, \| \Delta B_{m-1} \|_{\rho_m}, \ \| D\Delta B_{m-1} \|_{\rho_m}, \ \| \nabla \Delta H_{m-1} \|_{\rho_m}$$

by the quantities which are in (5.20).

$\| D\Delta K_{m-1} \|_{\rho_m}$ is estimated as

$$\| D\Delta K_{m-1} \|_{\rho_{m-1}-3\delta_{m-1}} \le c_{m-1}^{(8)} \gamma^{-2} \delta_{m-1}^{-(2\sigma+1)} \| e_{m-1} \|_{\rho_{m-1}}$$

by (5.10). Let's estimate $\| \Delta N_{m-1} \|_{\rho_m}$. We denote the Banach space of all continuous linear maps from Banach space $\mathbf{E}$ to Banach space $\mathbf{F}$ by $L(\mathbf{E}, \mathbf{E})$ and denote the set of all topological linear isomorphisms from $\mathbf{E}$ to $\mathbf{F}$ denote by $IL(\mathbf{E}, \mathbf{E})$. Then $IL(\mathbf{E}, \mathbf{E})$ is open set in $L(\mathbf{E}, \mathbf{E})$. If we define

$$\Phi(X) = X^{-1}, \ (X \in IL(\mathbf{R}^n, \ \mathbf{R}^n)),$$

then we have

$$D\Phi(X)\mathrm{H} = -X^{-1}\mathrm{H}X^{-1}, \ (\forall \mathrm{H} \in L(\mathbf{R}^n, \ \mathbf{R}^n)).$$

Accordingly

$$\| D\Phi(X) \| = \sup_{\|\mathrm{H}\|=1} \| X^{-1}\mathrm{H}X^{-1} \| \le \| X^{-1} \|^2$$

holds. Then from the integral-typed Taylor's formula we have

$$\Phi(Y) - \Phi(X) = \int_0^1 D\Phi(X + t(Y-X))(Y-X)dt =$$
$$= -\int_0^1 [X + t(Y-X)]^{-1}(Y-X)[X + t(Y-X)]^{-1} dt.$$

Hence

$$N_m - N_{m-1} = (DK_m{}^T DK_m)^{-1} - (DK_{m-1}{}^T DK_{m-1})^{-1} =$$
$$= -\int_0^1 E(t)^{-1}(DK_m{}^T DK_m - DK_{m-1}{}^T DK_{m-1})E(t)^{-1} dt =$$
$$= -\int_0^1 E(t)^{-1} Z_{m-1} E(t)^{-1} dt$$

holds, where

$$E(t) = DK_{m-1}{}^T DK_{m-1} + t(DK_m{}^T DK_m - DK_{m-1}{}^T DK_{m-1}) =$$
$$= Y_{m-1} + t Z_{m-1} = Y_{m-1}(I_n + t N_{m-1} Z_{m-1}).$$

If we consider (5.19), then we have

$$\| E(t)^{-1} \|_{\rho_m} = \| (I_n + tN_{m-1}Z_{m-1})^{-1} N_{m-1} \|_{\rho_m} \leq (1 - t \| N_{m-1}Z_{m-1} \|_{\rho_m})^{-1} \| N_{m-1} \|_{\rho_m} \leq 2 \| N_{m-1} \|_{\rho_m}.$$

Hence

$$\| \Delta N_{m-1} \|_{\rho_m} = \| N_m - N_{m-1} \|_{\rho_m} \leq 4 \| N_{m-1} \|_{\rho_m}^2 \| Z_{m-1} \|_{\rho_m} \leq$$
$$\leq 4 \| N_{m-1} \|_{\rho_m}^2 [2 \| DK_{m-1} \|_{\rho_{m-1}} c_{m-1}^{(9)} + \frac{1}{2} c_{m-1}^{(9)}] \gamma^{-2} \delta_{m-1}^{-2(\sigma+1)} \| e_{m-1} \|_{\rho_{m-1}} \leq .$$
$$\leq c_{m-1}^{(10)} \gamma^{-2} \delta_{m-1}^{-2(\sigma+1)} \| e_{m-1} \|_{\rho_{m-1}}.$$

Now we estimate $\| \Delta A_{m-1} \|_{\rho_m}$. If $\xi \in \mathbf{R}^{2n}$ is the element satisfying $\|\xi\| = 1$, then

$$\| \Delta A_{m-1}(\theta)\xi \| \leq \| A_m(\theta)\xi - A_{m-1}(\theta)\xi \| =$$
$$= \| DB(K_m(\theta))\xi \nabla H(K_m(\theta)) + B(K_m(\theta))D\nabla H(K_m(\theta))\xi -$$
$$- DB(K_{m-1}(\theta))\xi \nabla H(K_{m-1}(\theta)) - B(K_{m-1}(\theta))D\nabla H(K_{m-1}(\theta))\xi \| \leq$$
$$\leq \| DB(K_m(\theta))\xi \nabla H(K_m(\theta)) - DB(K_{m-1}(\theta))\xi \nabla H(K_{m-1}(\theta)) \| +$$
$$\| B(K_m(\theta))D\nabla H(K_m(\theta))\xi - B(K_{m-1}(\theta))D\nabla H(K_{m-1}(\theta))\xi \|.$$

From mean value theorem we have

$$\| DB(K_m(\theta))\xi \nabla H(K_m(\theta)) - DB(K_{m-1}(\theta))\xi \nabla H(K_{m-1}(\theta)) \| \leq$$
$$\leq \| DB(K_m(\theta))\xi \nabla H(K_m(\theta)) - DB(K_{m-1}(\theta))\xi \nabla H(K_m(\theta)) \| +$$
$$+ \| DB(K_{m-1}(\theta))\xi \nabla H(K_m(\theta)) - DB(K_{m-1}(\theta))\xi \nabla H(K_{m-1}(\theta)) \| \leq$$
$$\leq \| DB(K_m(\theta)) - DB(K_{m-1}(\theta)) \| \cdot \| \nabla H(K_m(\theta)) \| +$$
$$+ \| DB(K_{m-1}(\theta)) \| \cdot \| \nabla H(K_m(\theta)) - \nabla H(K_{m-1}(\theta)) \| \leq$$
$$\leq \| B \|_{C^2, \mathcal{B}_r} |H|_{C^1, \mathcal{B}_r} \| \Delta K_{m-1} \|_{\rho_m} + \| B \|_{C^1, \mathcal{B}_r} |H|_{C^2, \mathcal{B}_r} \| \Delta K_{m-1} \|_{\rho_m} \leq$$
$$\leq 2 \| B \|_{C^2, \mathcal{B}_r} |H|_{C^2, \mathcal{B}_r} \| \Delta K_{m-1} \|_{\rho_m} \leq$$
$$\leq 2 \| B \|_{C^2, \mathcal{B}_r} |H|_{C^2, \mathcal{B}_r} c_{m-1}^{(8)} \gamma^{-2} \delta_{m-1}^{-2\sigma} \| e_{m-1} \|_{\rho_{m-1}}.$$

Similarly we obtain

$$\| B(K_m(\theta))D\nabla H(K_m(\theta))\xi - B(K_{m-1}(\theta))D\nabla H(K_{m-1}(\theta))\xi \| \leq$$
$$\leq \| B \|_{C^1, \mathcal{B}_r} |H|_{C^2, \mathcal{B}_r} \| \Delta K_{m-1} \|_{\rho_m} + \| B \|_{\mathcal{B}_r} |H|_{C^3, \mathcal{B}_r} \| \Delta K_{m-1} \|_{\rho_m} \leq$$
$$\leq 2 \| B \|_{C^1, \mathcal{B}_r} |H|_{C^3, \mathcal{B}_r} \| \Delta K_{m-1} \|_{\rho_m} \leq$$
$$\leq 2 \| B \|_{C^1, \mathcal{B}_r} |H|_{C^3, \mathcal{B}_r} c_{m-1}^{(8)} \gamma^{-2} \delta_{m-1}^{-2\sigma} \| e_{m-1} \|_{\rho_{m-1}}.$$

Therefore we have

$$\left\| \Delta A_{m-1} \right\|_{\rho_m} \leq 4 \| B \|_{C^2, \mathcal{B}_r} |H|_{C^3, \mathcal{B}_r} c_{m-1}^{(8)} \gamma^{-2} \delta_{m-1}^{-2\sigma} \| e_{m-1} \|_{\rho_{m-1}}.$$

Similarly $\| \Delta B_{m-1} \|_{\rho_m}$, $\| D\Delta B_{m-1} \|_{\rho_m}$, $\| \nabla \Delta H_{m-1} \|_{\rho_m}$ are estimated as follows

$$\| \Delta B_{m-1} \|_{\rho_{m-1}} \leq \| B \|_{C^1, \mathcal{B}_r} \| \Delta K_{m-1} \|_{\rho_m} \leq \| B \|_{C^1, \mathcal{B}_r} c_{m-1}^{(8)} \gamma^{-2} \delta_{m-1}^{-2\sigma} \| e_{m-1} \|_{\rho_{m-1}},$$

$$\| D\Delta B_{m-1} \|_{\rho_{m-1}} \leq \| B \|_{C^2, \mathcal{B}_r} \| \Delta K_{m-1} \|_{\rho_m} \leq \| B \|_{C^2, \mathcal{B}_r} c_{m-1}^{(8)} \gamma^{-2} \delta_{m-1}^{-2\sigma} \| e_{m-1} \|_{\rho_{m-1}},$$

$$\|\nabla \Delta H_{m-1}\|_{\rho_m} \leq |H|_{C^2,\mathcal{B}_r} \|\Delta K_{m-1}\|_{\rho_m} \leq |H|_{C^2,\mathcal{B}_r} c_{m-1}^{(8)} \gamma^{-2} \delta_{m-1}^{-2\sigma} \|e_{m-1}\|_{\rho_{m-1}}$$

by mean value theorem.

From the above arguments there exists constant $c_{m-1}^{(11)}$ depending on

$$n, \sigma, \|N_{m-1}\|_{\rho_{m-1}}, \|DK_{m-1}\|_{\rho_{m-1}}, |<S_{m-1}^0>^{-1}|, \|B\|_{C^2,\mathcal{B}_r}, \|B^{-1}\|_{\mathcal{B}_r}, \|H\|_{C^3,\mathcal{B}_r}$$

that satisfies $c_{m-1}^{(11)} \geq c_{m-1}^{(10)}$ such that

$$\|\Delta S_{m-1}^0\|_{\rho_m} \leq c_{m-1}^{(11)} \gamma^{-2} \delta_{m-1}^{-2(\sigma+1)} \|e_{m-1}\|_{\rho_{m-1}}. \tag{5.20}$$

Since from the assumption 3) of lemma 5.2: $c_{m-1}^{(12)} \gamma^{-2} \delta_{m-1}^{-(2\sigma+1)} \|e_{m-1}\|_{\rho_{m-1}} \leq 1/2$ we have

$$\begin{aligned}
\left|<S_{m-1}^0>^{-1}\right| \cdot \left|<S_{m-1}^0>\right| &\leq \left|<S_{m-1}^0>^{-1}\right| \cdot \|S_{m-1}^0\|_{\rho_m} \leq \\
&\leq \left|<S_{m-1}^0>^{-1}\right| c_{m-1}^{(11)} \gamma^{-2} \delta_{m-1}^{-(2\sigma+1)} \|e_{m-1}\|_{\rho_{m-1}} < \\
&< 2(|<S_{m-1}^0>^{-1}|+1)^2 c_{m-1}^{(11)} \gamma^{-2} \delta_{m-1}^{-(2\sigma+1)} \|e_{m-1}\|_{\rho_{m-1}} = \\
&= c_{m-1}^{(12)} \gamma^{-2} \delta_{m-1}^{-(2\sigma+1)} \|e_{m-1}\|_{\rho_{m-1}} \leq 1/2,
\end{aligned} \tag{5.21}$$

matrix $(I_n + <S_{m-1}^0>^{-1} <\Delta S_{m-1}^0>)$ and $(I_n + N_{m-1} Z_{m-1})$ are invertible. From this fact

$$<S_m^0> = <S_{m-1}^0>[I_n + <S_{m-1}^0>^{-1} <\Delta S_{m-1}^0>]$$

is invertible and

$$<S_m^0>^{-1} = [I_n + <S_{m-1}^0>^{-1} <\Delta S_{m-1}^0>]^{-1} <S_{m-1}^0>^{-1}.$$

Hence if we consider (5.20),(5.21), then we have

$$\begin{aligned}
|<S_m^0>^{-1}| &\leq |I_n + <S_{m-1}^0>^{-1} <\Delta S_{m-1}^0>]^{-1} | \cdot |<S_{m-1}^0>^{-1}| \\
&\leq (1-|<S_{m-1}^0>^{-1}| \cdot |<\Delta S_{m-1}^0>|)^{-1} \cdot |<S_{m-1}^0>^{-1}| = \\
&= [1 + (1-|<S_{m-1}^0>^{-1}| \cdot |<\Delta S_{m-1}^0>|)^{-1} |<S_{m-1}^0>^{-1}| \cdot |<\Delta S_{m-1}^0>|] \cdot |<S_{m-1}^0>^{-1}| \leq \\
&\leq |<S_{m-1}^0>^{-1}| + 2|<S_{m-1}^0>^{-1}|^2 \cdot |<\Delta S_{m-1}^0>| \leq \\
&\leq |<S_{m-1}^0>^{-1}| + 2|<S_{m-1}^0>^{-1}|^2 c_{m-1}^{(11)} \gamma^{-2} \delta_{m-1}^{-2(\sigma+1)} \|e_{m-1}\|_{\rho_{m-1}} = \\
&= |<S_{m-1}^0>^{-1}| + c_{m-1}^{(12)} \gamma^{-2} \delta_{m-1}^{-2(\sigma+1)} \|e_{m-1}\|_{\rho_{m-1}}.
\end{aligned}$$

□

Let's put $c_{m-1} = c_{m-1}^{(12)}$, $(m \geq 1)$. We can obtain the following lemma 5.3 by lemma 5.1 and lemma 5.2:

**Lemma 5.3.** If for $m \geq 1$ the approximation solution of (2.16) $K_{m-1} \in \mathcal{P}_{\rho_{m-1}}$ satisfies

1) $DK_{m-1}(\theta)^T DK_{m-1}(\theta)$, $<S_{m-1}^0>$ are invertible  (5.22)

2) $r_{m-1} + c_{m-1} \gamma^{-2} \delta_{m-1}^{-2\sigma} \|e_{m-1}\|_{\rho_{m-1}} < r$  (5.23)

3) $c_{m-1}\gamma^{-2}\delta_{m-1}^{-(2\sigma+1)} \|e_{m-1}\|_{\rho_{m-1}} \leq 1/2$, (5.24)

then $K_m(\theta) = K_{m-1}(\theta) + DK_{m-1}(\theta)$, $(\theta \in U_{\rho_m})$ obtained by quasi-Newton method in one step, satisfies

a) $\|\Delta K_{m-1}\|_{\rho_{m-1}-2\delta_{m-1}} \leq c_{m-1}\gamma^{-2}\delta_{m-1}^{-2\sigma} \|e_{m-1}\|_{\rho_{m-1}}$ (5.25)

b) $K_m \in \mathcal{P}_{\rho_m}$ (5.26)

c) $DK_m(\theta)^T DK_m(\theta)$, $<S_m^0>$ is invertible (5.27)

d) $K_m(U_{\rho_m}) \subset \mathcal{B}_r$ (5.28)

e) $\|e_m\|_{\rho_m} \leq c_{m-1}\gamma^{-4}\delta_{m-1}^{-4\sigma} \|e_{m-1}\|_{\rho_{m-1}}^2$. (5.29)

Therefore $K_m(\theta)$, $(\theta \in U_{\rho_m})$ is approximate solution of (2.16) and $K_m \in \mathcal{ND}(\rho_m)$. □

Now we establish the condition for $\|e_0\|_{\rho_0}$ in which (5.22), (5.23) hold for any $m \geq 1$.

Constant $c_{m-1}$ depends on parameters $n$, $\sigma$, $\rho$, $\gamma$, $\|B\|_{C^2,\mathcal{B}_r}$, $\|B^{-1}\|_{\mathcal{B}_r}$, $\|H\|_{C^3,\mathcal{B}_r}$. $c_{m-1}$ also depends on $\rho_{m-1} \leq \rho_0$ and quantities

$$d_{m-1} = \|DK_{m-1}\|_{\rho_{m-1}}, \quad v_{m-1} = \|N_{m-1}\|_{\rho_{m-1}}, \quad s_{m-1} = \left|<S_{m-1}^0>^{-1}\right| \quad (m \geq 1)$$

related to the approximate solution $K_{m-1}$. Moreover we know that dependency of $c_m$ for $(d_m, v_m, s_m)$ is polynomial with positive coefficients from their definition process. Thus there exists polynomial with positive coefficients $\lambda(y_1, y_2, y_3)$ related to

$n$, $\sigma$, $\rho$, $\gamma$, $\|B\|_{C^2,\mathcal{B}_r}$, $\|B^{-1}\|_{\mathcal{B}_r}$, $\|H\|_{C^3,\mathcal{B}_r}$ such that

$$c_{m-1} = \lambda(d_{m-1}, v_{m-1}, s_{m-1}), \quad \forall m \geq 1. \tag{5.30}$$

Since we are finding the solution of (2.16) close to $K_0$, we expect that quantities $d_m, v_m, s_m$ remain in a neighbourhood of $d_0, v_0, \tau_0$ respectively. In this case, $c_{m-1}$ is bounded by

$$\lambda(d_m + \beta, v_m + \beta, s_m + \beta),$$

where $\beta$ will be conveniently determined.

**Theorem 5.1.** Let's $\{c_m\}_{m \geq 0}$ is sequence of constants that is given by (5.30). For fixed $0 < \delta_0 < \min(\rho_0/12, 1)$ we define

$$\delta_m = \delta_0 2^{-m}, \quad \varepsilon_m = \|e_m\|_{\rho_m}, \quad r_m = \|K_m - K_0\|_{\rho_m}, \quad m \geq 0.$$

Then there exists constant $c$ depending on

$$n, \ \sigma, \ r, \ \rho_0, \|N_{m-1}\|_{\rho_0}, \|DK_{m-1}\|_{\rho_0}, |<S_0^0>^{-1}|, \|B\|_{C^2, \mathcal{B}_r}, \|B^{-1}\|_{\mathcal{B}_r}, \|H\|_{C^3, \mathcal{B}_r}$$

such that: if $\|e_0\|_{\rho_0}$ satisfies

$$\kappa := 2^{4\sigma} c \gamma^{-4} \delta_0^{-4\sigma} \|e_0\|_{\rho_0} \leq 1/2 , \qquad (5.31)$$

$$\left[1 + \frac{2^{4\sigma}}{2^{2\sigma}-1}\right] c \gamma^{-2} \delta_0^{-2\sigma} \|e_0\|_{\rho_0} < r , \qquad (5.32)$$

then (5.22),(5.23),(5.24) hold for $m \geq 1$.

**Proof.** Let

$$\beta = \gamma^2 \delta_0^{2\sigma-1} 2^{-(4\sigma+1)}(1+2^{4\sigma-1}), \ c = \lambda(d_0+\beta, \ \nu_0+\beta, \ \tau_0+\beta). \qquad (5.33)$$

We take four proposition related to $m \in \mathbf{Z}_+$ as follows:

C1(m) $r_m \leq \left(1 + \dfrac{\kappa 2^{4\sigma}}{2^{2\sigma}-1}\right) c \gamma^{-2} \delta_0^{-2\sigma} \varepsilon_0$

C2(m) $\varepsilon_m \leq 2^{-4\sigma(m-1)} \kappa^{(2^m-1)} \varepsilon_0$

C3(m) $c_m \leq c$

C4(m) $c_m \gamma^{-2} \delta_m^{-(2\sigma+1)} \varepsilon_m \leq 1/2$

First step: We will prove that C1(m),C2(m),C3(m),C4(m) ( $m \geq 0$ ) mean (5.22),(5.23),(5.24) for any $m \geq 1$ under the assumptions (5.31),(5.32).

C4(m) is just (5.24).

Next we prove that for $m=1$ (5.23) holds. From $r_0 = \|K_0 - K_0\|_{\rho_0} = 0$, (5.32) and C3(0) we have

$$r_0 + c_0 \gamma^{-2} \delta_0^{-2\sigma} \varepsilon_0 = c_0 \gamma^{-2} \delta_0^{-2\sigma} \varepsilon_0 \leq c \gamma^{-2} \delta_0^{-2\sigma} \varepsilon_0 < r .$$

Therefore (5.23) holds for $m=1$. Next we prove (5.23) holds for $m \geq 2$. Let $m \geq 1$. From C1(m),C2(m),C3(m) and $\delta_m = \delta_0 2^{-m}$ we have

$$r_m + c_m \gamma^{-2} \delta_m^{-2\sigma} \|e_m\|_{\rho_m} \leq$$

$$\leq \left(1 + \frac{\kappa 2^{4\sigma}}{2^{2\sigma}-1}\right) c\gamma^{-2}\delta_0^{-2\sigma}\varepsilon_0 + c\gamma^{-2}\delta_0^{-2\sigma} 2^{2m\sigma} 2^{-4\sigma(m-1)} \kappa^{(2^m-1)} \varepsilon_0 =$$

$$= c\gamma^{-2}\delta_0^{-2\sigma}\varepsilon_0 \left\{1 + \frac{\kappa 2^{4\sigma}}{2^{2\sigma}-1} + 2^{2m\sigma} 2^{-4\sigma(m-1)} \kappa^{(2^m-1)}\right\} =$$

$$= c\gamma^{-2}\delta_0^{-2\sigma}\varepsilon_0 \left\{1 + \frac{\kappa 2^{4\sigma}}{2^{2\sigma}-1} + 2^{-2\sigma m+4\sigma} \kappa^{(2^m-1)}\right\} =$$

$$= c\gamma^{-2}\delta_0^{-2\sigma}\varepsilon_0 \left\{1 + \frac{\kappa 2^{4\sigma} + 2^{-2\sigma m+6\sigma}\kappa^{(2^m-1)} - 2^{-2\sigma m+4\sigma}\kappa^{(2^m-1)}}{2^{2\sigma}-1}\right\}.$$

Therefore if we obtain

$$\kappa 2^{4\sigma} + 2^{-2\sigma m+6\sigma}\kappa^{(2^m-1)} - 2^{-2\sigma m+4\sigma}\kappa^{(2^m-1)} \leq 2^{4\sigma}, \tag{5.34}$$

then we get (5.23) from (5.32). By the way since we have $\kappa 2^{4\sigma} \leq \frac{2^{4\sigma}}{2}$ from (5.31), if

$$2^{-2\sigma m+6\sigma}\kappa^{(2^m-1)} \leq 2^{4\sigma-1} \tag{5.35}$$

is true, then (5.34) is true. However since

$$2^{-2\sigma m+6\sigma}\kappa^{(2^m-1)} \leq 2^{-2\sigma m+6\sigma}(2^{-1})^{(2^m-1)} = 2^{-2\sigma m+6\sigma}2^{-2^m+1} = 2^{-2\sigma m+6\sigma-2^m+1},$$

if

$$-2\sigma m + 6\sigma - 2^m + 1 \leq 4\sigma - 1 \tag{5.36}$$

is true, then (5.23) is true. Here inequality

$$-2\sigma m - 2^m + 2 + 2\sigma \leq 0$$

that is equivalent to (5.36) holds and therefore (5.23) holds. By lemma 5.2 we have (5.22) from (5.23),(5.24). Accordingly if C1(m)-C4(m) are true for $m \geq 0$, then (5.22)-(5.24) hold.

Second step: Next we prove C1(m)-C4(m) for $m \geq 0$ under the assumptions (5.31),(5.32). First of all, we prove C1(0)-C4(0). Since $r_0 = \|K_0 - K_0\|_{\rho_0} = 0$, C1(0) holds. Since C2(0) implies $\varepsilon_0 \leq \varepsilon_0$, C2(0) holds trivially. Since from the definition of $c_0$, $c$ we have

$$c_0 = \lambda(d_0,\ \nu_0,\ \tau_0) \leq \lambda(d_0+\beta,\ \nu_0+\beta,\ \tau_0+\beta) = c,$$

C3(0) holds. By the assumptions $\sigma > 1$ and $0 < \gamma \leq 1$, from C3(0) and

(5.31): $\kappa = 2^{4\sigma} c \gamma^{-4} \delta_0^{-4\sigma} \|e_0\|_{\rho_0} \leq 1/2$ we have

$$c_0 \gamma^{-2} \delta_0^{-(2\sigma+1)} \varepsilon_0 \leq (c\gamma^{-4}\delta_0^{-4\sigma}\varepsilon_0)\delta_0^{2\sigma-1} \leq 2^{-4\sigma-1} < 2^{-4\sigma} < 1/2.$$

And C4(0) holds.

Next we prove C1(m)-C4(m) under the inductive assumptions C1(j)-C4(j) for $j = 1, \cdots, m-1$ (we can perform $m-1$ step of quasi-Newton method under the assumptions). First we confirm

$$2^{j-1}\sum_{s=1}^{j-1} s 2^{-s} \leq 2^j - j - 1, \quad (j \geq 2). \tag{5.37}$$

In the case of $j = 2$ since

$$2 \cdot 2^{-1} \leq 2^2 - 2 - 1,$$

(5.37) is true. Now we prove (5.37) for case of $j+1$ under the inductive assumptions (5.37) for case of $j \geq 2$. Since

$$2^j \sum_{s=1}^{j} s 2^{-s} = 2 \cdot 2^{j-1}\left[\sum_{s=1}^{j-1} s 2^{-s} + j 2^{-j}\right] \leq 2 \cdot \left[2^j - j - 1 + 2^{j-1} j 2^{-j}\right] =$$
$$= 2^{j+1} - 2(j+1) + j = 2^{j+1} - (j+1) - 1,$$

(5.37) holds in the case of $j+1$.

Now we prove that (5.31),(5.32),C1(m-1),C2(m-1),C3(m-1) mean (5.23). (5.23) holds for $m = 1$. Let's prove (5.23) for $m \geq 2$. Since we have

$$r_m + c_m \gamma^{-2} \delta_m^{-2\sigma} \|e_m\|_{\rho_m} \leq$$

$$\leq [1 + \frac{\kappa 2^{4\sigma}}{2^{2\sigma} - 1}]c\gamma^{-2}\delta_0^{-2\sigma}\varepsilon_0 + c\gamma^{-2}\delta_0^{-2\sigma} 2^{2m\sigma} 2^{-4\sigma(m-1)} \kappa^{(2^m-1)} \varepsilon_0 =$$

$$= c\gamma^{-2}\delta_0^{-2\sigma}\varepsilon_0 \left\{1 + \frac{\kappa 2^{4\sigma}}{2^{2\sigma} - 1} + 2^{2m\sigma} 2^{-4\sigma(m-1)} \kappa^{(2^m-1)}\right\} =$$

$$= c\gamma^{-2}\delta_0^{-2\sigma}\varepsilon_0 \left\{1 + \frac{\kappa 2^{4\sigma}}{2^{2\sigma} - 1} + 2^{-2\sigma m + 4\sigma} \kappa^{(2^m-1)}\right\} =$$

$$= c\gamma^{-2}\delta_0^{-2\sigma}\varepsilon_0 \left\{1 + \frac{\kappa 2^{4\sigma} + 2^{-2\sigma m + 6\sigma} \kappa^{(2^m-1)} - 2^{-2\sigma m + 4\sigma} \kappa^{(2^m-1)}}{2^{2\sigma} - 1}\right\}$$

considering $\delta_m = \delta_0 2^{-m}$, for (5.23) it is sufficient to prove

$$\kappa 2^{4\sigma} + 2^{-2\sigma m + 6\sigma} \kappa^{(2^m-1)} - 2^{-2\sigma m + 4\sigma} \kappa^{(2^m-1)} \leq 2^{4\sigma}$$

considering (5.32). Since $\kappa 2^{4\sigma} \leq \frac{2^{4\sigma}}{2}$, for the last inequality it is sufficient to prove

$$2^{-2\sigma m + 6\sigma} \kappa^{(2^m-1)} - 2^{-2\sigma m + 4\sigma} \kappa^{(2^m-1)} \leq 2^{4\sigma - 1}.$$

Now we prove $2^{-2\sigma m + 6\sigma} \kappa^{(2^m-1)} \leq 2^{4\sigma - 1}$. Since

$$2^{-2\sigma m + 6\sigma} \kappa^{(2^m-1)} \leq 2^{-2\sigma m + 6\sigma} (2^{-1})^{(2^m-1)} = 2^{-2\sigma m + 6\sigma} 2^{-2^m + 1} = 2^{-2\sigma m + 6\sigma - 2^m + 1},$$

it is sufficient to prove $-2\sigma m + 6\sigma - 2^m + 1 \leq 4\sigma - 1$ for our aim. Since the inequality $-2\sigma m + 6\sigma - 2^m + 1 \leq 4\sigma - 1$ is equivalent to $-2\sigma m - 2^m + 2 + 2\sigma \leq 0$, (5.23) holds.

On the other hand, we have $c_{m-1}\gamma^{-1}\delta_{m-1}^{-(\sigma+1)}\varepsilon_{m-1} \leq 1/2$ from C4(m-1). Therefore all of assumptions in lemma 5.1 are satisfied. From lemma 5.1 we have

$$(5.13): \|e_m\|_{\rho_m} \leq c_{m-1}^{(9)}\gamma^{-4}\delta_{m-1}^{-4\sigma}\|e_{m-1}\|_{\rho_{m-1}}^2 .$$

Now we prove

$$\varepsilon_j \leq \kappa^{2^j-1}(2^{-4\sigma j})\varepsilon_0 \leq \kappa^{2^j-1}(2^{-4\sigma(j-1)})\varepsilon_0 \quad (1 \leq j \leq m) \tag{5.38}$$

by inductive method. First we consider the case of $j=1$. Since

$$\varepsilon_1 \leq c\gamma^{-4}\delta_0^{-4\sigma}\varepsilon_0^2 = 2^{4\sigma}c\gamma^{-4}\delta_0^{-4\sigma}\varepsilon_0^2 2^{-4\sigma} = \kappa 2^{-4\sigma}\varepsilon_0,$$

(5.38) holds for $j=1$. Next in the case of $j \geq 2$ we prove (5.38) for $j$ under the inductive assumption (5.38) for $j-1$. Since if we use (5.37) considering

$(5.29): \|e_m\|_{\rho_m} \leq c_{m-1}\gamma^{-4}\delta_{m-1}^{-4\sigma}\|e_{m-1}\|_{\rho_{m-1}}^2$ and C3(m-1), then

$$\begin{aligned}
\varepsilon_j &\leq c\gamma^{-4}\delta_{j-1}^{-4\sigma}\varepsilon_{j-1}^2 \leq c\gamma^{-4}\delta_{j-1}^{-4\sigma}(c\gamma^{-4}\delta_{j-2}^{-4\sigma}\varepsilon_{j-2}^2)^2 = \\
&= c\gamma^{-4}(c\gamma^{-4})^2\delta_{j-1}^{-4\sigma}(\delta_{j-2}^{-4\sigma})^2\varepsilon_{j-2}^4 \leq c\gamma^{-4}(c\gamma^{-4})^2\delta_{j-1}^{-4\sigma}(\delta_{j-2}^{-4\sigma})^2(c\gamma^{-4}\delta_{j-3}^{-4\sigma}\varepsilon_{j-3}^2)^4 = \\
&= c\gamma^{-4}(c\gamma^{-4})^2(c\gamma^{-4})^4\delta_{j-1}^{-4\sigma}(\delta_{j-2}^{-4\sigma})^2(\delta_{j-3}^{-4\sigma})^4\varepsilon_{j-3}^8 = \\
&= (c\gamma^{-4})^{1+2+4}(\delta_{j-1}\delta_{j-2}^2\delta_{j-3}^4)^{-4\sigma}\varepsilon_{j-3}^8 \leq \\
&\leq \cdots \leq (c\gamma^{-4})^{(1+2+\cdots+2^{j-1})}(\delta_{j-1}^{2^0}\delta_{j-2}^{2^1}\cdots\delta_0^{2^{j-1}})^{-4\sigma}\varepsilon_0^{2^j} = \\
&= (c\gamma^{-4})^{(2^0+2^1+\cdots+2^{j-1})}((\delta_0 2^{-(j-1)})^{2^0}(\delta_0 2^{-(j-2)})^{2^1}\cdots(\delta_0 2^{-0})^{2^{j-1}})^{-4\sigma}\varepsilon_0^{2^j} = \\
&= (c\gamma^{-4}\delta_0^{-4\sigma})^{(2^0+2^1+\cdots+2^{j-1})}((2^{j-1})^{2^0}(2^{j-2})^{2^1}\cdots(2^0)^{2^{j-1}})^{4\sigma}\varepsilon_0^{2^j} = \\
&= (c\gamma^{-4}\delta_0^{-4\sigma})^{2^j-1}(2^{4\sigma})^{2^0(j-1)+2^1(j-2)+\cdots+2^{j-2}\cdot 1}\varepsilon_0^{2^j} \leq \\
&\leq (c\gamma^{-4}\delta_0^{-4\sigma})^{2^j-1}(2^{4\sigma})^{2^j-j-1}\varepsilon_0^{2^j} \\
&\leq \underbrace{(c\gamma^{-4}\delta_0^{-4\sigma}2^{4\sigma}\varepsilon_0)}_{\kappa}^{2^j-1}2^{-4\sigma j}\varepsilon_0 \leq \kappa^{2^j-1}(2^{-4\sigma j})\varepsilon_0.
\end{aligned}$$

And we have (5.38). From (5.38) we obtain C2(m).

On the other hand since $(5.9): \|\Delta K_{m-1}\|_{\rho_{m-1}-2\delta_{m-1}} \leq c_{m-1}^{(8)}\gamma^{-2}\delta_{m-1}^{-2\sigma}\|e_{m-1}\|_{\rho_{m-1}}$ holds from lemma 5.1, we have

$$r_m = \|K_0 - K_m\|_{\rho_m} \le \|K_0 - K_{m-1}\|_{\rho_{m-1} - 3\delta_{m-1}} + \|K_m - K_{m-1}\|_{\rho_{m-1} - 3\delta_{m-1}} \le$$

$$\le \|K_0 - K_{m-1}\|_{\rho_{m-1}} + \|\Delta K_m\|_{\rho_{m-1} - 2\delta_{m-1}} \le r_{m-1} + c\gamma^{-2}\delta_{m-1}^{-2\sigma}\|e_{m-1}\|_{\rho_{m-1}} \le$$

$$\le r_{m-2} + c\gamma^{-2}\delta_{m-2}^{-2\sigma}\|e_{m-2}\|_{\rho_{m-2}} + c\gamma^{-2}\delta_{m-1}^{-2\sigma}\|e_{m-1}\|_{\rho_{m-1}} \le \cdots \le$$

$$\le c\gamma^{-2}\delta_0^{-2\sigma}\varepsilon_0 + c\gamma^{-2}\sum_{j=1}^{m-1}\delta_j^{-2\sigma}\varepsilon_j \le$$

$$\le c\gamma^{-2}\delta_0^{-2\sigma}\varepsilon_0 + c\gamma^{-2}\sum_{j=1}^{m-1}(\delta_0 2^{-j})^{-2\sigma}\cdot 2^{-4\sigma(j-1)}\kappa^{(2^j-1)}\varepsilon_0 \le$$

$$\le c\gamma^{-2}\delta_0^{-2\sigma}\varepsilon_0 + c\gamma^{-2}\delta_0^{-2\sigma}\kappa\varepsilon_0\sum_{j=1}^{m-1}2^{2j\sigma}2^{-4\sigma(j-1)} =$$

$$= c\gamma^{-2}\delta_0^{-2\sigma}\varepsilon_0 + c\gamma^{-2}\delta_0^{-2\sigma}\kappa\varepsilon_0 2^{4\sigma}\sum_{j=1}^{m-1}2^{-2j\sigma} \le$$

$$\le c\gamma^{-2}\delta_0^{-2\sigma}\varepsilon_0\left(1 + \kappa 2^{4\sigma}\frac{2^{-2\sigma}}{1 - 2^{-2\sigma}}\right) =$$

$$= c\gamma^{-2}\delta_0^{-2\sigma}\varepsilon_0\left(1 + \frac{\kappa 2^{4\sigma}}{2^{2\sigma} - 1}\right)$$

considering (5.38). Therefore C1(m) holds.

Now we prove C3(m). Since from

$$(5.10): \|D\Delta K_{m-1}\|_{\rho_{m-1} - 3\delta_{m-1}} \le c_{m-1}^{(8)}\gamma^{-2}\delta_{m-1}^{-(2\sigma+1)}\|e_{m-1}\|_{\rho_{m-1}}$$

in lemma 5.1 we have

$$d_m = \|DK_m\|_{\rho_m} = \|D(K_{m-1} + \Delta K_{m-1})\|_{\rho_m} \le \|DK_{m-1}\|_{\rho_m} + \|D\Delta K_{m-1}\|_{\rho_m} \le$$

$$\le \|DK_{m-1}\|_{\rho_{m-1}} + \|D\Delta K_{m-1}\|_{\rho_{m-1} - 3\delta_{m-1}} \le d_{m-1} + c_{m-1}^{(8)}\gamma^{-2}\delta_{m-1}^{-(2\sigma+1)}\varepsilon_{m-1} \le$$

$$\le d_{m-1} + c_{m-1}\gamma^{-2}\delta_{m-1}^{-(2\sigma+1)}\varepsilon_{m-1},$$

$$d_m \le d_{m-1} + c_{m-1}\gamma^{-2}\delta_{m-1}^{-(2\sigma+1)}\varepsilon_{m-1}$$

holds. If we apply the last inequality repeatedly, then we have

$$d_m \le d_{m-1} + c_{m-1}\gamma^{-2}\delta_{m-1}^{-(2\sigma+1)}\varepsilon_{m-1}$$
$$\le d_{m-2} + c_{m-2}\gamma^{-2}\delta_{m-2}^{-(2\sigma+1)}\varepsilon_{m-2} + c_{m-1}\gamma^{-2}\delta_{m-1}^{-(2\sigma+1)}\varepsilon_{m-1} \le \cdots \le$$
$$\le d_0 + \sum_{j=0}^{m-1}c_j\gamma^{-2}\delta_j^{-(2\sigma+1)}\varepsilon_j = d_0 + \tau_{m-1},$$

where $\tau_{m-1} = \sum_{j=0}^{m-1}c_j\gamma^{-2}\delta_j^{-(2\sigma+1)}\varepsilon_j$. Similarly from lemma 5.1 we have

$$v_m \leq v_{m-1} + c_{m-1}\gamma^{-2}\delta_{m-1}^{-(2\sigma+1)}\varepsilon_{m-1} \leq v_0 + \tau_{m-1},$$
$$s_m \leq s_{m-1} + c_{m-1}\gamma^{-2}\delta_{m-1}^{-(2\sigma+1)}\varepsilon_{m-1} \leq s_0 + \tau_{m-1}.$$

Let's estimate $\tau_{m-1}$. From (5.37): $\varepsilon_j \leq \kappa^{2^j-1}(2^{-4\sigma(j-1)})\varepsilon_0$ and C3(j)( $j=1, \cdots, m-1$) we have

$$c_j\gamma^{-2}\delta_j^{-(2\sigma+1)}\varepsilon_j \leq c\gamma^{-2}(\delta_0 2^{-j})^{-(2\sigma+1)} \cdot \kappa^{2^j-1}\varepsilon_0 2^{-4\sigma(j-1)} =$$
$$= (c\gamma^{-2}\delta_0^{-(2\sigma+1)}\varepsilon_0)\kappa^{2^j-1}2^{4\sigma-2\sigma j+j} = (c\gamma^{-2}\delta_0^{-(2\sigma+1)}\varepsilon_0)2^{4\sigma}\kappa^{2^j-1}2^{-j(2\sigma-1)}$$

for $j=1, \cdots, m-1$. If we consider $\kappa = 2^{4\sigma}c\gamma^{-4}\delta_0^{-4\sigma}\varepsilon_0 \leq 1/2$ and $(2^{2\sigma-1}-1)^{-1} \leq 1$, then we have

$$\tau_{m-1} = \sum_{j=0}^{m-1}c_j\gamma^{-2}\delta_j^{-(2\sigma+1)}\varepsilon_j \leq$$
$$\leq c\gamma^{-2}\delta_0^{-(2\sigma+1)}\varepsilon_0 + (c\gamma^{-2}\delta_0^{-(2\sigma+1)}\varepsilon_0)2^{4\sigma}\sum_{j=1}^{m-1}\kappa^{(2^j-1)}2^{-j(2\sigma-1)} =$$
$$= c\gamma^{-2}\delta_0^{-(2\sigma+1)}\varepsilon_0(1+2^{4\sigma}\sum_{j=1}^{m-1}\kappa 2^{-j(2\sigma-1)}) \leq$$
$$\leq (c\gamma^{-4}\delta_0^{-4\sigma}\varepsilon_0 2^{4\sigma}) \cdot \gamma^2\delta_0^{2\sigma-1}2^{-4\sigma}(1+\kappa 2^{4\sigma}\sum_{j=1}^{m-1}2^{-j(2\sigma-1)}) =$$
$$= \kappa\gamma^2\delta_0^{2\sigma-1}2^{-4\sigma}(1+\kappa 2^{4\sigma}\sum_{j=1}^{m-1}2^{-j(2\sigma-1)}) \leq$$
$$\leq \kappa\gamma^2\delta_0^{2\sigma-1}2^{-4\sigma}(1+\kappa 2^{4\sigma}\frac{1}{2^{2\sigma-1}-1}) \leq$$
$$\leq \gamma^2\delta_0^{2\sigma-1}2^{-(4\sigma+1)}(1+2^{4\sigma-1}) = \beta.$$

Therefore we have

$$d_m \leq d_0 + \beta, \ v_m \leq v_0 + \beta, \ s_m \leq s_0 + \beta.$$

This means $c_m \leq c$. Thus C3(m) holds.

On the other hand, since from (5.38) we have

$$c_m\gamma^{-2}\delta_m^{-(2\sigma+1)} \cdot \varepsilon_m \leq c\gamma^{-2}(\delta_0 2^{-m})^{-(2\sigma-1)} \cdot \kappa^{2^m-1}\varepsilon_0 2^{-4\sigma(m-1)} =$$
$$= c\gamma^{-2}\delta_0^{-(2\sigma-1)}2^{m(2\sigma-1)} \cdot \kappa^{2^m-1}\varepsilon_0 2^{-4\sigma(m-1)} =$$
$$= 2^{4\sigma}c\gamma^{-4}\delta_0^{-4\sigma}\varepsilon_0 \cdot \delta_0^{2\sigma+1}2^{m(2\sigma-1)}\gamma^2 \cdot \kappa^{2^m-1}2^{-4\sigma m} =$$
$$= \kappa\delta_0^{2\sigma+1}\gamma^2\kappa^{2^m-1}2^{-m(2\sigma+1)} = \delta_0^{2\sigma+1}\gamma^2\kappa^{2^m}2^{-m(2\sigma+1)} \leq$$
$$\leq 2^{-2^m-m(2\sigma+1)} \leq 1/2,$$

C4(m) holds. □

# 6. Existence of KAM tori

## 6.1. Convergence of a type of Nash-Moser iteration

**Definition 6.1.** Let's $\rho_0 > 0$. Suppose that for any $\rho \in [0, \rho_0]$ real Banach space $(\mathbf{E}_\rho, \|\cdot\|_\rho)$ is given and they satisfy the following two conditions: For any $0 \le \rho_1 \le \rho_2 \le \rho_0$,

1) $\mathbf{E}_{\rho_2}$ is a sub Banach space of $\mathbf{E}_{\rho_1}$. Thus $\mathbf{E}_{\rho_1} \supset \mathbf{E}_{\rho_2}$.

2) For any $x \in \mathbf{E}_{\rho_1}$, $\|x\|_{\rho_1} \le \|x\|_{\rho_2}$.

Then we call $\{(\mathbf{E}_\rho, \|\cdot\|_\rho)\}_{0 \le \rho \le \rho_0}$ as scale of Banach spaces. Suppose that map $F: \mathbf{E}_{\rho_0} \to \mathbf{E}_{\rho_0}$ satisfies $F(\mathbf{E}_\rho) \subset \mathbf{E}_\rho$ for any $\rho \in [0, \rho_0]$. Then we call that map $F$ coincide with the scale of Banach spaces $\{(\mathbf{E}_\rho, \|\cdot\|_\rho)\}_{0 \le \rho \le \rho_0}$. □

Assume that $C^1$ map $F: \mathbf{E}_{\rho_0} \to \mathbf{E}_{\rho_0}$ is given, which coincide with the scale of Banach spaces $\{(\mathbf{E}_\rho, \|\cdot\|_\rho)\}_{0 \le \rho \le \rho_0}$. And assume that for any $\rho \in [0, \rho_0]$ open set $ND(\rho) \subset \mathbf{E}_\rho$ is given. Now we consider a procedure of quasi-Newton method which clarify the existence of solution $K \in ND(\rho)$ of the equation

$$F(K) = 0 \qquad (6.1)$$

for some $\rho \in (0, \rho_0]$. Let's fix $\rho \in (0, \rho_0]$ arbitrarily. Assume that an element $K \in \mathbf{E}_\rho$ is given. If $X \in ND(\rho)$ satisfies $F(X) = 0$, then from the Taylor's formula we have

$$F(K) + DF(K)(X - K) + O(\|X - K\|_\rho^2) = 0. \qquad (6.2)$$

If we linearize the equation (6.2) neglecting higher terms in that, then we have

$$DF(K)(X - K) = -F(K).$$

Therefore if a procedure to obtain the approximate solution $\Delta K$ of the linearized equation

$$DF(K)\Delta K = -F(K) \qquad (6.3)$$

is given, then $X = K + \Delta K$ would give an approximate solution of nonlinear equation (6.1). We call (6.3) as linearized equation of (6.1) at point $K$.

**Definition 6.2.** If $K \in \mathbf{E}_{\rho_0}$ satisfies $K \in ND(\rho)$ for some $\rho \in (0, \rho_0]$, then we call that $K$ is the approximate solution of (6.1) with error $e = F(K)$. We call $\varepsilon = \|e\|_\rho = \|F(K)\|_\rho \ge 0$ as error-size of approximate solution $K$ for (6.1). On the other hand, if $K \in \mathbf{E}_{\rho_0}$ satisfies $K \in ND(\rho)$ for $\rho \in (0, \rho_0]$ and $F(K) = 0$, then we call that $K$ is the solution of (6.1). Thus solution of (6.1) is just the approximate solution of (6.1) with zero error.

Let's suppose that $K \in \mathbf{E}_{\rho_0}$ satisfies $K \in ND(\rho)$ for $\rho \in (0, \rho_0]$. We call any element $\Delta K \in \mathbf{E}_{\rho_0}$ as approximate solution of (6.3) at $K$ with error $\varepsilon' = DF(K)\Delta K + F(K)$. □

**Fundamental Assumptions 6.1.** Assume that real numbers $c > 0$, $\gamma > 0$, $\sigma \ge 1$ are given. We choose $\delta_0 > 0$ as

$$0 < \delta_0 \le \min\{1, \rho_0/12\}. \qquad (6.4)$$

Next we put

$$\delta_m = \delta_0 2^{-m} \ (m \ge 1), \quad \rho_m = \rho_{m-1} - 3\delta_{m-1} \ (m \ge 1). \qquad (6.5)$$

Then

$$\rho_m = \rho_{m-1} - 3\delta_{m-1} = \rho_0 - 6\delta_0(1 - 2^{-m}) > 0$$

is obvious. And put

$$\rho_\infty := \lim_{m \to \infty} \rho_m = \rho_0 - 6\delta_0 > 0 . \tag{6.6}$$

We assume that an approximate solution $K_0 \in ND(\rho_0)$ of equation (6.1) is given. We call $K_0$ as 0-th approximate solution for (6.1).

Let's assume that $m$-th an approximate solution of (6.1) $K_m \in ND(\rho_m)$ is given. Then assume that a rule to give an approximate solution of (6.3) at point $K_m$, $\Delta K_m \in \mathbf{E}_{\rho_m}$ such that

$$K_{m+1} := K_m + \Delta K_m \in ND(\rho_{m+1}) \tag{6.7}$$

(therefore $K_{m+1} = K_m + \Delta K_m$ become approximate solution of (6.1) in this time). We put

$$e_m = F(K_m) \ (m \geq 0), \tag{6.8}$$

$$\varepsilon_m = \|e_m\|_{\rho_m} \ (m \geq 0), \tag{6.9}$$

$$\kappa = 2^{4\sigma} c \gamma^{-4} \delta_0^{-4\sigma} \|e_0\|_{\rho_0} . \tag{6.10}$$

**Theorem 6.1** We assume that:

1) There exists a sequence of non-negative real numbers $\{c_m\}_{m \geq 0}$ such that

$$c_m \leq c . \tag{6.11}$$

2) $\kappa := 2^{4\sigma} c \gamma^{-4} \delta_0^{-4\sigma} \varepsilon_0 \leq \dfrac{1}{2}$ (thus $\varepsilon_0 \leq c^{-1} 2^{-4\sigma-1} \gamma^4 \delta_0^{4\sigma}$). (6.12)

3) $\|e_m\|_{\rho_m} \leq c_{m-1} \gamma^{-4} \delta_{m-1}^{-4\sigma} \|e_{m-1}\|_{\rho_{m-1}}^2 \quad (m \geq 1)$. (6.13)

4) $\|\Delta K_m\|_{\rho_m - 2\delta_m} \leq c_m \gamma^{-2} \delta_m^{-2\sigma} \|e_m\|_{\rho_m} \quad (m \geq 0)$. (6.14)

Then $\{K_m\}$ converges to some solution $K_\infty \in \mathbf{E}_{\rho_\infty}$ of (6.1) and $K_\infty$ satisfies

$$\|K_\infty - K_0\|_{\rho_\infty} \leq c \gamma^{-2} \delta_0^{-2\sigma} \varepsilon_0 (1 + \kappa \frac{2^{4\sigma}}{2^{2\sigma} - 1}) . \tag{6.15}$$

**Proof.** First we prove

$$2^{j-1} \sum_{s=1}^{j-1} s 2^{-s} \leq 2^j - j - 1, \ (j \geq 2) . \tag{6.16}$$

Since

$$2 \cdot 2^{-1} \leq 2^2 - 2 - 1 ,$$

then (6.16) holds in the case of $j = 2$. Now we prove (6.16) holds in the case of $j+1$ under the assumption that (6.16) holds for $j \geq 2$. Since

$$2^j \sum_{s=1}^{j} s 2^{-s} = 2 \cdot 2^{j-1} \left[ \sum_{s=1}^{j-1} s 2^{-s} + j 2^{-j} \right] \leq 2 \cdot \left[ 2^j - j - 1 + 2^{j-1} j 2^{-j} \right] =$$

$$= 2^{j+1} - 2(j+1) + j = 2^{j+1} - (j+1) - 1 ,$$

(6.16) holds in the case of $j+1$. Now we prove

$$\varepsilon_j \leq \kappa^{2^j - 1} (2^{-4\sigma j}) \varepsilon_0 \leq \kappa^{2^j - 1} (2^{-4\sigma(j-1)}) \varepsilon_0 \quad (j \geq 1) \tag{6.17}$$

by inductive method. First we consider case of $j = 1$. Since

$$\varepsilon_1 \leq c \gamma^{-4} \delta_0^{-4\sigma} \varepsilon_0^2 = 2^{4\sigma} c \gamma^{-4} \delta_0^{-4\sigma} \varepsilon_0^2 2^{-4\sigma} = \kappa 2^{-4\sigma} \varepsilon_0 ,$$

(6.17) holds for $j=1$. Now we prove (6.17) holds in the case of $j$ under the assumption that (6.17) holds in the case of $j-1$ for $j \geq 2$. If we apply (6.16) considering (6.13) and (6.11), then we have

$$\varepsilon_j \leq c\gamma^{-4}\delta_{j-1}^{-4\sigma}\varepsilon_{j-1}^2 \leq c\gamma^{-4}\delta_{j-1}^{-4\sigma}(c\gamma^{-4}\delta_{j-2}^{-4\sigma}\varepsilon_{j-2}^2)^2 =$$
$$= c\gamma^{-4}(c\gamma^{-4})^2\delta_{j-1}^{-4\sigma}(\delta_{j-2}^{-4\sigma})^2\varepsilon_{j-2}^4 \leq c\gamma^{-4}(c\gamma^{-4})^2\delta_{j-1}^{-4\sigma}(\delta_{j-2}^{-4\sigma})^2(c\gamma^{-4}\delta_{j-3}^{-4\sigma}\varepsilon_{j-3}^2)^4 =$$
$$= c\gamma^{-4}(c\gamma^{-4})^2(c\gamma^{-4})^4\delta_{j-1}^{-4\sigma}(\delta_{j-2}^{-4\sigma})^2(\delta_{j-3}^{-4\sigma})^4\varepsilon_{j-3}^8 =$$
$$= (c\gamma^{-4})^{1+2+4}(\delta_{j-1}\delta_{j-2}^2\delta_{j-3}^4)^{-4\sigma}\varepsilon_{j-3}^8 \leq$$
$$\leq \cdots \leq (c\gamma^{-4})^{(1+2+\cdots+2^{j-1})}(\delta_{j-1}^{2^0}\delta_{j-2}^{2^1}\cdots\delta_1^{2^{j-2}}\delta_0^{2^{j-1}})^{-4\sigma}\varepsilon_0^{2^j} \leq$$
$$\leq (c\gamma^{-4})^{(2^0+2^1+\cdots+2^{j-1})}((\delta_0 2^{-(j-1)})^{2^0}(\delta_0 2^{-(j-2)})^{2^1}\cdots(\delta_0 2^{-1})^{2^{j-2}}(\delta_0 2^{-0})^{2^{j-1}})^{-4\sigma}\varepsilon_0^{2^j} =$$
$$= (c\gamma^{-4}\delta_0^{-4\sigma})^{(2^0+2^1+\cdots+2^{j-1})}((2^{j-1})^{2^0}(2^{j-2})^{2^1}\cdots(2^1)^{2^{j-2}}(2^0)^{2^{j-1}})^{4\sigma}\varepsilon_0^{2^j} =$$
$$= (c\gamma^{-4}\delta_0^{-4\sigma})^{(2^0+2^1+\cdots+2^{j-1})}(2^{4\sigma})^{2^0(j-1)+2^1(j-2)+\cdots+2^{j-2}\cdot 1}\varepsilon_0^{2^j} \leq$$
$$\leq (c\gamma^{-4}\delta_0^{-4\sigma})^{(2^0+2^1+\cdots+2^{j-1})}(2^{4\sigma})^{2^j-j-1}\varepsilon_0^{2^j}$$
$$\leq (\underbrace{c\gamma^{-4}\delta_0^{-4\sigma}2^{4\sigma}\varepsilon_0}_{\kappa})^{2^j-1}2^{-4\sigma j}\varepsilon_0 \leq \kappa^{2^j-1}(2^{-4\sigma j})\varepsilon_0.$$

Therefore (6.17) holds. We obtain

$$K_1 = K_0 + \Delta K_0$$
$$K_2 = K_1 + \Delta K_1 = K_0 + \Delta K_0 + \Delta K_1$$
$$\vdots$$
$$K_{m-1} = K_{m-2} + \Delta K_{m-2} = K_0 + \sum_{j=0}^{m-2}\Delta K_j$$
$$K_m = K_{m-1} + \Delta K_{m-1} = K_0 + \sum_{j=0}^{m-1}\Delta K_j.$$

From (6.14): $\|\Delta K_m\|_{\rho_m-2\delta_m} \leq c_m\gamma^{-2}\delta_m^{-2\sigma}\|e_m\|_{\rho_m}$ we have

$$\|\Delta K_j\|_{\rho_{j+1}} = \|\Delta K_j\|_{\rho_j-3\delta_j} \leq \|\Delta K_j\|_{\rho_j-2\delta_j} \leq c_j\gamma^{-2}\delta_j^{-2\sigma}\|e_j\|_{\rho_j}$$

for $j \geq 0$. Hence if we consider (6.12): $\kappa = 2^{4\sigma}c\gamma^{-4}\delta_0^{-4\sigma}\|e_0\|_{\rho_0} \leq 1/2$ and

(6.17): $\varepsilon_m \leq 2^{-4\sigma(m-1)}\kappa^{(2^m-1)}\varepsilon_0$, then we have

$$\|\Delta K_j\|_{\rho_\infty} \leq \|\Delta K_j\|_{\rho_{j+1}} \leq c_j\gamma^{-2}\delta_j^{-2\sigma}\cdot\varepsilon_j \leq c\gamma^{-2}(\delta_0 2^{-j})^{-2\sigma}\cdot 2^{-4\sigma(j-1)}\kappa^{2^j-1}\varepsilon_0 =$$
$$= c\gamma^{-2}\delta_0^{-2\sigma}\kappa^{2^j-1}\varepsilon_0 2^{4\sigma}2^{-2\sigma j} \leq c\gamma^{-2}\delta_0^{-2\sigma}\kappa\varepsilon_0 2^{4\sigma}2^{-2\sigma j}$$

(6.18)

for $j \geq 1$. Since from (6.18) and (6.14) for $m=0$: $\|\Delta K_0\|_{\rho_0-2\delta_0} \leq c_0\gamma^{-2}\delta_0^{-2\sigma}\varepsilon_0$ we have

$$\sum_{j=0}^{\infty}\|\Delta K_j\|_{\rho_\infty} \leq c\gamma^{-2}\delta_0^{-2\sigma}\varepsilon_0(1+2^{4\sigma}\kappa\sum_{j=1}^{\infty}2^{-2\sigma j}) = c\gamma^{-2}\delta_0^{-2\sigma}\varepsilon_0(1+\kappa\frac{2^{4\sigma}}{2^{2\sigma}-1}), \quad (6.19)$$

$\sum_{j=0}^{\infty} \Delta K_j$ thus $\{K_m\}$ converges in $\mathbf{E}_{\rho_\infty}$ and

$$K_\infty = K_0 + \sum_{j=0}^{\infty} \Delta K_j. \tag{6.20}$$

From (6.19) and (6.20) we have

$$\|K_\infty - K_0\|_{\rho_\infty} = \left\|\sum_{j=0}^{\infty} \Delta K_j\right\|_{\rho_\infty} \leq \sum_{j=0}^{\infty} \|\Delta K_j\|_{\rho_\infty} \leq c\,\gamma^{-2}\delta_0^{-2\sigma}\varepsilon_0(1+\kappa\frac{2^{4\sigma}}{2^{2\sigma}-1}).$$

On the other hand since from (6.17): $\varepsilon_j \leq \kappa^{2^j-1}(2^{-4\sigma(j-1)})\varepsilon_0$ we have $\lim_{m\to+\infty}\varepsilon_m = 0$,

$$\|F(K_\infty)\|_{\rho_\infty} = \lim_{m\to+\infty}\|F(K_m)\|_{\rho_\infty} = \lim_{m\to+\infty}\|e_m\|_{\rho_m} = 0$$

holds. Therefore $K_\infty$ is the solution of (6.1). □

**Lemma 6.1** Suppose that the following three assumptions are satisfied:

1) Map $F:\mathbf{E}_{\rho_0} \to \mathbf{E}_{\rho_0}$ is $C^2$. Here we fix $\mathrm{H} \geq 0$ such that

$$\|D^2F(K_{m-1}+t\Delta K_{m-1})\| \leq \mathrm{H},\quad (t\in[0,\ 1]). \tag{6.21}$$

For example we can take as $\mathrm{H} = \sup_{0\leq t\leq 1}\|D^2F(K_{m-1}+t\Delta K_{m-1})\|$.

2) There exists sequence $\{c''_m\}_{m\geq 0}$ of non-negative real numbers such that

$$\|\Delta K_{m-1}\|_{\rho_{m-1}-2\delta_{m-1}} \leq c''_{m-1}\gamma^{-2}\delta_{m-1}^{-2\sigma}\|e_{m-1}\|_{\rho_{m-1}}\quad (m\geq 1). \tag{6.22}$$

3) There exists sequence $\{c'''_m\}_{m\geq 0}$ of non-negative real numbers such that

$$\|e_{m-1} + DF(K_{m-1})\Delta K_{m-1}\|_{\rho_{m-1}-2\delta_{m-1}} \leq c'''_m\gamma^{-3}\delta_{m-1}^{-(3\sigma+1)}\|e_{m-1}\|_{\rho_{m-1}}^2\quad (m\geq 1). \tag{6.23}$$

Then $K_m = K_{m-1}+\Delta K_{m-1} \in ND(\rho_m)$ satisfies

$$\|e_m\|_{\rho_m} \leq c^{(4)}_{m-1}\gamma^{-4}\delta_{m-1}^{-4\sigma}\|e_{m-1}\|_{\rho_{m-1}}^2 \tag{6.24}$$

for $c^{(4)}_m = c'''_m + c''^{\,2}_m\mathrm{H},\ (m\geq 0)$.

**Proof.** Let's suppose that $\Delta K_{m-1} \in \mathbf{E}_{\rho_{m-1}}$ is an approximate solution of $DF(K_{m-1})\Delta = -e_{m-1}$ given by fundamental assumption 6.1. If we define remainder term in the Taylor's formula as
$$\mathcal{R}(K,\ K') = F(K') - F(K) - DF(K)(K'-K),$$
then
$$F(K_m) = F(K_{m-1}) + DF(K_{m-1})(K_m - K_{m-1}) + \mathcal{R}(K_{m-1},\ K_m).$$

Therefore
$$e_m = e_{m-1} + DF(K_{m-1})\Delta K_{m-1} + \mathcal{R}(K_{m-1},\ K_m).$$

Hence we have
$$\|e_m\|_{\rho_m} \leq \|e_{m-1}+DF(K_{m-1})\Delta K_{m-1}\|_{\rho_m} + \|\mathcal{R}(K_{m-1},\ K_m)\|_{\rho_m}.$$

Since $\rho_m = \rho_{m-1}-3\delta_{m-1}$ from (6.23), we have
$$\|e_{m-1}+DF(K_{m-1})\Delta K_{m-1}\|_{\rho_m} = \|e_{m-1}+DF(K_{m-1})\Delta K_{m-1}\|_{\rho_{m-1}-3\delta_{m-1}} \leq$$
$$\leq \|e_{m-1}+DF(K_{m-1})\Delta K_{m-1}\|_{\rho_{m-1}-2\delta_{m-1}} \leq c'''_{m-1}\gamma^{-3}\delta_{m-1}^{-(3\sigma+1)}\|e_{m-1}\|_{\rho_{m-1}}^2.$$

Now let's estimate $\|\mathcal{R}(K_{m-1}, K_m)\|_{\rho_m}$. Since from the integral type Taylor's formula we have

$$F(K_m) = F(K_{m-1}) + DF(K_{m-1})\Delta K_{m-1} + \int_0^1 D^2 F(K_{m-1} + t\Delta K_{m-1})(1-t)(\Delta K_{m-1})^2 dt,$$

by (6.22) we have

$$\|\mathcal{R}(K_{m-1}, K_m)\|_{\rho_m} = \left\|\int_0^1 D^2 F(K_{m-1} + t\Delta K_{m-1})(1-t)(\Delta K_{m-1})^2 dt\right\|_{\rho_m} =$$

$$\leq \mathrm{H} \cdot \|\Delta K_{m-1}\|_{\rho_m}^2 \leq \mathrm{H}\|\Delta K_{m-1}\|_{\rho_{m-1}-2\delta_{m-1}}^2 \leq c_{m-1}''^2 \mathrm{H}\gamma^{-4}\delta_{m-1}^{-4\sigma}\|e_{m-1}\|_{\rho_{m-1}}^2.$$

Therefore we have

$$\|e_m\|_{\rho_m} \leq c_{m-1}''' \gamma^{-3}\delta_{m-1}^{-(3\sigma+1)}\|e_{m-1}\|_{\rho_{m-1}}^2 + c_{m-1}''^2 \mathrm{H}\gamma^{-4}\delta_{m-1}^{-4\sigma}\|e_{m-1}\|_{\rho_{m-1}}^2 \leq$$

$$\leq (c_{m-1}''' + c_{m-1}''^2 \mathrm{H})\gamma^{-4}\delta_{m-1}^{-4\sigma}\|e_{m-1}\|_{\rho_{m-1}}^2 = c_{m-1}^{(4)}\gamma^{-4}\delta_{m-1}^{-4\sigma}\|e_{m-1}\|_{\rho_{m-1}}^2. \square$$

**Theorem 6.2.** We suppose that

1) Map $F:\mathbf{E}_{\rho_0} \to \mathbf{E}_{\rho_0}$ is $C^2$. Here we fix $\mathrm{H} \geq 0$ such that

(6.21): $\|D^2 F(K_{m-1} + t\Delta K_{m-1})\| \leq \mathrm{H}$, $(t \in [0, 1])$.

2) There exists a sequence $\{c_m''\}_{m \geq 0}$ of non-negative real numbers such that

$$\|\Delta K_{m-1}\|_{\rho_{m-1}-2\delta_{m-1}} \leq c_{m-1}'' \gamma^{-2}\delta_{m-1}^{-2\sigma}\|e_{m-1}\|_{\rho_{m-1}} \quad (m \geq 1).$$

3) There exists a sequence $\{c_m'''\}_{m \geq 0}$ of non-negative real numbers such that

$$\|e_{m-1} + DF(K_{m-1})\Delta K_{m-1}\|_{\rho_{m-1}-2\delta_{m-1}} \leq c_m''' \gamma^{-4}\delta_{m-1}^{-(3\sigma+1)}\|e_{m-1}\|_{\rho_{m-1}}^2 \quad (m \geq 1).$$

4) If we put $c_m = c_m''' + c_m''^2 \mathrm{H}$, $(m \geq 0)$, then $c_m \leq c$.

5) $\kappa = 2^{4\sigma} c \gamma^{-4}\delta_0^{-4\sigma}\varepsilon_0 \leq \dfrac{1}{2}$.

Then $\{K_m\}$ converges to a solution $K_\infty \in \mathbf{E}_{\rho_\infty}$ such that

$$\|K_\infty - K_0\|_{\rho_\infty} \leq c\gamma^{-2}\delta_0^{-2\sigma}\varepsilon_0(1 + \kappa\frac{2^{4\sigma}}{2^{2\sigma}-1}).$$

**Proof.** From lemma 6.1 all assumptions in theorem 6.1 are satisfied for $c_m = c_m''' + c_m''^2 \mathrm{H}$, $(m \geq 0)$. Therefore from theorem 6.1 we obtain the conclusion of theorem 6.2. $\square$

### 6.2. Existence of KAM tori

**Theorem 6.3** Let $\omega \in \mathbf{R}^n$ belongs to $D_n(\gamma, \sigma)$ for some $\gamma > 0$ and $\sigma > n-1$. Let's suppose that $0 < \rho_0 \leq 12$ and $K_0 \in \mathcal{P}(\rho_0)$ are given which is non-degenerate. Let's consider generalized Hamiltonian system

(2.1): $\dfrac{dz}{dt} = B(z)\nabla H(z)$.

We suppose $H:\mathbf{U}^{2n} \to \mathbf{R}$ and $b_{ij}:\mathbf{U}^{2n} \to \mathbf{R}$ are real analytic and they can be holomorphically extended to some complex neighbourhood of the image under $K_0$ of $U_{\rho_0}$:

$$\mathcal{B}_r = \{z \in \mathbf{C}^{2n} \; ; \; \inf_{|\mathrm{Im}\theta| < \rho_0} |z - K_0(\theta)| < r\}$$

for some $r > 0$. We define error function for $K_0$ as $e_0(\theta) = J\nabla H(K_0(\theta)) - \partial_\omega K_0(\theta)$. Let

$\delta_0 = \rho_0/12$. And we put

$$d_0 = \|DK_0\|_{\rho_0}, \ \nu_0 = \|N_0\|_{\rho_0}, \ s_0 = |<S_0>^{-1}|,$$
$$\beta = \gamma^2 \delta_0^{2\sigma-1} 2^{-(4\sigma+1)}(1+2^{4\sigma-1}),$$
$$c = \lambda(d_0+\beta, \ \nu_0+\beta, \ s_0+\beta),$$
$$C = (1+2^{4\sigma-1})c.$$

If $e_0$ satisfies

$$C\gamma^{-4}\delta_0^{-4\sigma}\|e_0\|_{\rho_0} \leq 1/2, \tag{6.25}$$

$$C\gamma^{-2}\delta_0^{-2\sigma}\|e_0\|_{\rho_0} < r, \tag{6.26}$$

then there exists solution $K_\infty \in \mathcal{ND}(\rho/2)$ of

$$(2.16): \partial_\omega K(\theta) = B(K(\theta))\nabla H(K(\theta))$$

such that

$$\|K_\infty - K_0\|_{\rho_0/2} < r.$$

**Proof.** Let $K_0 \in \mathcal{ND}(\rho_0)$ is given. We put

$$\delta_m = \delta_0 2^{-m} \ (m \geq 1), \ \rho_m = \rho_{m-1} - 3\delta_{m-1} \ (m \geq 1).$$

Then

$$\rho_m = \rho_{m-1} - 3\delta_{m-1} = \rho_0 - 6\delta_0(1-2^{-m}) > 0.$$

Let $\rho_\infty := \lim_{m \to \infty} \rho_m$. Then $\rho_\infty = \rho_0 - 6\delta_0 = \rho_0/2$. Assume that (6.25),(6.26) hold for $C = (1+2^{4\sigma-1})c$.

Then

$$(1+2^{4\sigma-1})c\gamma^{-4}\delta_0^{-4\sigma}\|e_0\|_{\rho_0} \leq 1/2 \tag{6.27}$$

$$(1+2^{4\sigma-1})c\gamma^{-2}\delta_0^{-2\sigma}\|e_0\|_{\rho_0} < r \tag{6.28}$$

hold. From (6.27) we have

$$2^{4\sigma}c\gamma^{-4}\delta_0^{-4\sigma}\|e_0\|_{\rho_0} \leq 1/2 \tag{6.29}$$

And from $\sigma \geq 1$, $\dfrac{1}{2^{2\sigma}-1} \leq \dfrac{1}{2}$ holds. Therefore from (6.28) we have

$$\left(1+\frac{2^{4\sigma}}{2^{2\sigma}-1}\right)c\gamma^{-2}\delta_0^{-2\sigma}\|e_0\|_{\rho_0} < r. \tag{6.30}$$

Then from theorem 5.1

a) $r_{m-1} + c_{m-1}\gamma^{-2}\delta_{m-1}^{-2\sigma}\|e_{m-1}\|_{\rho_{m-1}} < r \ (m \geq 1)$

b) $c_{m-1}\gamma^{-2}\delta_{m-1}^{-(\sigma+1)}\|e_{m-1}\|_{\rho_{m-1}} \leq 1/2 \quad (m \geq 1)$

c) $c_m \leq c$

hold. And then from 5.3 $K_m = K_{m-1} + DK_{m-1}$ satisfies

1) $\|\Delta K_{m-1}\|_{\rho_{m-1}-2\delta_{m-1}} \leq c_m\gamma^{-2}\delta_{m-1}^{-2\sigma}\|e_{m-1}\|_{\rho_{m-1}}$,

2) $K_m \in \mathcal{P}_{\rho_m}$,

3) $DK_m(\theta)^T DK_m(\theta)$, $<S_m>$ is invertible,

4) $K_m(U_{\rho_m}) \subset \mathcal{B}_r$,

5) $\|e_m\|_{\rho_m} \le c_{m-1}\gamma^{-4}\delta_{m-1}^{-4\sigma}\|e_{m-1}\|_{\rho_{m-1}}^2$ .

Therefore from theorem 6.1 there exists solution $K_\infty \in \mathcal{ND}(\rho_0/2)$ of (2.16) such that

$$\|K_\infty - K_0\| \le c(1+\frac{\kappa 2^{4\sigma}}{2^{2\sigma}-1})\gamma^{-2}\delta_0^{-2\sigma} < r \;. \square$$

We can obtain the following local uniqueness of invariant tori just as [Llave, 2005].

**Theorem 6.4.** Let's assume that $\omega \in D(\gamma,\sigma)$ and $\delta = \rho/8$. Suppose that $K_1, K_2 \in \mathcal{ND}(\rho)$ are two solutions of (2.16) such that $K_1(U_\rho) \subset \mathcal{B}_r$ and $K_2(U_\rho) \subset \mathcal{B}_r$. For $N$ and $S^0$ in which we replace $K$ with $K_2$, there exists a constant $c > 0$, depending on

$$n,\; \sigma,\; \gamma,\; \rho,\; \rho^{-1},\; |H|_{C^2,\mathcal{B}_r},\; |B|_{\mathcal{B}_r},\; \left|B^{-1}\right|_{\mathcal{B}_r},\; \|K_2\|_{C^2,\rho},\; \|N\|_\rho,\; \left|(<S^0>)^{-1}\right|,$$

such that if $\|K_1 - K_2\|$ satisfies $c\gamma^{-2}\delta^{-2\sigma}\|K_1 - K_2\|_\rho < 1$, then there exists an initial phase $\tau \in \mathbf{R}^n$ such that $(K_1 \circ T_\tau)(\theta) = K_2(\theta)$ ($\theta \in U_{\rho/2}$) . $\square$

## References

[Ahlfors 1966] L.V.Ahlfors: Complex Analysis, McGraw-Hill,1966.

[Arnold 1963] Arnold, V. I., Proof of a theorem of A. N. Kolmogorov on the persistence of quasi-periodic motions under small perturbations of the Hamiltonian, Usp.Mat. Nauk 18:5 (1963), 13–40. English transl.: Russ. Math. Surv. 18:5 (1963).

[Arnold- Avez 1968] Arnold V.I., Avez A.: Ergodic Problems of Classical Mechanics, Benjamin-Cummings, Reading, Mass., 1968.

[Benettin 1984] G. Benettin, L. Galgani, A. Giorgili, and J.-M. Strelcyn. A proof of Kolmogorov's theorem on invariant tori using canonical transformations defined by the Lie method. Nuovo Cimento B (11), 79(2):201–223, 1984.

[Calleja 2011] R. Calleja: KAM theory for dissipative systems, McGill University, 2011.

[Cartan 1963] Henri Cartan: Theorie elémentaire des fonctions analytiquew d'une ou plusieurs variables complexes, Editions Scientifiques Hermann, Paris, 1963.

[Dieudonné 1960] Dieudonné J.: Foundations of modern analysis, Academic Press, 1960.

[Fejoz 2010] J.Fejoz; A simple proof of the invariant torus theorem, arXiv: 1105.5604v1 [math.DS], 2010, 19p.

[Fejoz 2011] J.Fejoz; A proof of the invariant torus theorem of Kolmogorov, arXiv:1102.0923v2 [math.DS], 2011, 4p.

[Hairer, 2006] Hairer E., Lubich C., Wanner G.: Geometric Numerical Integration 2ed, Springer, 2006.

[Haro, 2004] Haro A and de la Llave R: A parameterization method for the computation of whiskers in quasi periodic maps: numerical implementation Preprint mp arc 04-350.

[Irwin 1980] M.C.Irwin: Smooth Dynamical Systems,Academic Press, 1980.

[Ito 1998] Hideyuki Ito: Ordinary differential equation and analytical mechanics,Kyouritu Syuppann,1998 (Japanese)；伊藤秀行：常微分方程式と解析力学, 共立出版,1998.

[Jon-Paek 2014] WuHwan Jong, Jin Chol Paek: Existence of invariant tori for differentiable Hamiltonian vector fields without action-angle variables, Chaos, Solitons & Fractals 68 (2014) 114–122.

[Kolmogorov 1954] Kolmogorov A.N.: On conservation of conditionally periodic motions for a small perturbations in Hamiltonian's function, Dokl. Akad. Nauk SSSR 98:4 (1954), 527–530.

[Moser 1962] J. K. Moser, On invariant curves of area-preserving mappings of an annulus, Nachr. Akad. Wiss. Göttingen, Math.-Phys. Kl. II 1 (1962), 1–20.

[Llave, 2005] de la Llave R., González A., Jorba À., Villanueva J.: KAM theory without action-angle variables, Nonlinearity 18, 855–895, (2005).

[Li-Yi 2002] Yong Li, Yingfei Yi: Persistence of invariant tori in generalized Hamiltonian systems, Ergodic Theory Dynamical Systems 22 (2002), 1233–1261.

[Li-Yi 2006] Yong Li1 and Yingfei Yi : Nekhoroshev and KAM Stabilities in Generalized Hamiltonian Systems, Journal of Dynamics and Differential Equations, Vol. 18, No. 3, 2006.

[Liu-Yihe-Huang 2005] Zh. Liu, D. Yihe and Q. Huang, Persistence of hyperbolic tori in generalized Hamiltonian systems, Northeast. Math. J. 21 (2005), 447–464.

[Liu-Zhu-Han 2006] Baifeng Liu, Wenzhuang Zhu, Yuecai Han: Persistence of lower-dimensional hyperbolic invariant tori for generalized Hamiltonian systems, J. Math. Anal., Appl. 322 (2006), 251-275.

[Luque-Villanueva 2010] A. Luque, J.Villanueva: A KAM theorem without action-angle variables for elliptic lower dimensional tori, Universitat Politècnica de Catalunya, 1-64p, 2010.

[Pöschel 2009] J. Pöschel: A Lecture on the Classical KAM Theorem,ArXiv 0908.2234v1 [math.DS], 1-33p, (2009).

[Rüssmann 1975] H. Rüssmann: On Optimal Estimates for the Solutions of Linear PDEs of First Order with Constant Coefficients on the Torus, in Moser J. eds.: Dynamical Systems. Theory and Applications, Lecture Notes in Physics Vol. 38, Springer, 1975, 598-624.

[Sardanashvily 2003] G. Sardanashvily: The quasi-periodic stability condition (the KAMtheorem) for partiallyintegrable systems,airXiv: math/0301068v1 [math.DS] 8 Jan 2003, pp.1-4.

[Zhao 1995] Zhao Xiaohua: Periodic orbits in perturbed generalized Hamiltonian systems, Mathematica Acta Scientia, 15 (4), 370-384,1995.